\RequirePackage{anyfontsize}
\documentclass{svjour3}
\usepackage[caption=false]{subfig}
\usepackage{csquotes}
\usepackage{pgfplots}
\usepackage{pdfpages}

\pgfplotsset{compat=1.10}
\usepgfplotslibrary{fillbetween}
\usetikzlibrary{patterns}

\usepackage{url}
\usepackage{mathtools}
\usepackage{bbm}
\usepackage{enumitem}
\usepackage{esvect}
\usepackage{centernot}
\usepackage{float}
\usepackage{fullpage}
\usepackage{amsfonts}
\usepackage{latexsym}
\usepackage{amssymb}
\usepackage{color}
\usepackage{hyperref}
\usepackage{graphicx}
\graphicspath{{./figures/}}
\usepackage{epstopdf}
\usepackage{algorithm,algpseudocode}

\usepackage{array, multirow, boldline}

\usepackage{siunitx}
\usepackage{mathrsfs}
\usepackage{threeparttable}
\usepackage{textcomp}
\usepackage{booktabs}
\usepackage{tikz}
\usetikzlibrary{positioning,angles,quotes}

\usepackage{makecell}

\usepackage{algorithm,algpseudocode}

\usepackage{tikz}
\usetikzlibrary{bayesnet}
\usetikzlibrary{arrows}
\usetikzlibrary{decorations.pathreplacing}

\numberwithin{equation}{section}

\allowdisplaybreaks
\def\vgap{\vspace*{.1in}}
\newcommand{\ts}{\textsuperscript}
\newtheorem{assumption}{Assumption}
\makeatletter
\@ifundefined{theproperty}{\newtheorem{property}{Property}}{}
\makeatother

\newcommand{\R}{\mathbb{R}}

\newcommand{\La}{\mathcal{L}}

\def \grad {\nabla}
\DeclareMathOperator*{\argmax}{arg\,max}
\DeclareMathOperator*{\argmin}{arg\,min}

\newcommand{\inner}[2]{\langle {#1,#2} \rangle}
\NewDocumentCommand\norm{ m O{}}{\left\lVert #1 \right\rVert_{#2}}
\NewDocumentCommand\normsq{ m O{}}{\left\lVert #1 \right\rVert^2_{#2}}

\DeclarePairedDelimiter{\ceil}{\lceil}{\rceil}

\newcommand{\tsum}{\textstyle \sum}

\newcommand{\ep}{\epsilon}
\newcommand{\lam}{\lambda}

\DeclarePairedDelimiter\floor{\lfloor}{\rfloor}

\makeatletter
\newcommand{\algmargin}{\the\ALG@thistlm}
\makeatother
\newlength{\whilewidth}
\algdef{SE}[parWHILE]{parWhile}{EndparWhile}[1]
{\parbox[t]{\dimexpr\linewidth-\algmargin}{%
    \hangindent\whilewidth\strut\algorithmicwhile\ #1\ \algorithmicdo\strut}}{\algorithmicend\ \algorithmicwhile}%
\algnewcommand{\parState}[1]{\State%
  \parbox[t]{\dimexpr\linewidth-\algmargin}{\strut #1\strut}}

\def \bigO {\mathcal{O}}

\usepackage{pifont}

\usepackage{tcolorbox}

\usepackage{empheq}

\newcommand{\jim}[1]{#1}

\definecolor{jimdarkred}{RGB}{96,0,0}

\usepackage{enumitem}
\setlist{itemsep=0pt}
\usepackage{accents}
\NewDocumentCommand\utilde{m}{\underaccent{\tilde}{#1}}

\newcommand{\addt}[1]{#1^t}
\newcommand{\addtt}[1]{#1^{t-1}}
\newcommand{\addttt}[1]{#1^{t-2}}
\newcommand{\addp}[1]{#1^{t+1}}
\newcommand{\addN}[1]{#1^N}

\newcommand{\addO}[1]{#1^0}

\NewDocumentCommand\bci{O{i}}{#1}
\newcommand{\ave}{{\text{av}}}

\NewDocumentCommand\SubAddi{m m}{#1_{\bci[#2]}}

\NewDocumentCommand\SubAddI{m m}{#1_{#2}}
\NewDocumentCommand\Addit{m O{\bci}}{\addt{\SubAddI{#1}{#2}}}
\NewDocumentCommand\Additt{m O{\bci}}{\addtt{\SubAddI{#1}{#2}}}
\NewDocumentCommand\Addittt{m O{\bci}}{\addttt{\SubAddI{#1}{#2}}}
\NewDocumentCommand\AddiN{m O{\bci}}{\addN{\SubAddI{#1}{#2}}}
\NewDocumentCommand\Addip{m O{\bci}}{\addp{\SubAddI{#1}{#2}}}
\NewDocumentCommand\AddiO{m O{\bci}}{\addO{\SubAddI{#1}{#2}}}
\NewDocumentCommand\Addav{m O{\ave}}{\SubAddI{#1}{#2}}
\NewDocumentCommand\Addavt{m O{\ave}}{\addt{\SubAddI{#1}{#2}}}
\NewDocumentCommand\Addavtt{m O{\ave}}{\addtt{\SubAddI{#1}{#2}}}
\NewDocumentCommand\Addavttt{m O{\ave}}{\addttt{\SubAddI{#1}{#2}}}
\NewDocumentCommand\AddavN{m O{\ave}}{\addN{\SubAddI{#1}{#2}}}
\NewDocumentCommand\Addavp{m O{\ave}}{\addp{\SubAddI{#1}{#2}}}
\NewDocumentCommand\AddavO{m O{\ave}}{\addO{\SubAddI{#1}{#2}}}

\def \targmin {{\textstyle\argmin}}

\NewDocumentCommand\sumi{O{m}}{\tsum_{i=1}^{#1}}
\NewDocumentCommand\nSt{O{t}}{{{S_{#1}}}}
\NewDocumentCommand\St{O{t}}{S_{#1}}

\NewDocumentCommand\sums{O{\St}}{\tsum_{s=1}^{#1}}

\NewDocumentCommand\sumt{O{N}}{\tsum_{t=1}^{#1}}
\def \sumwt {\tsum_{t=1}^N \wt}

\NewDocumentCommand\sumT{O{N}}{\tsum_{t=1}^{#1}}

\NewDocumentCommand\sumwT{O{N}}{\tsum_{t=1}^{#1} \wt}
\NewDocumentCommand\sumwNtil{O{\Ntil}}{\tsum_{t=1}^{#1} \wt}
\NewDocumentCommand\qs{O{s}}{q_{#1}}

\NewDocumentCommand\etat{O{t}}{\eta_{#1}}
\def \etatt {\etat[t-1]}

\def \etatt {\eta_{t-1}}

\NewDocumentCommand\taut{O{t}}{\tau_{#1}}
\def \tautt {\taut[t-1]}
\NewDocumentCommand\thetat{O{t}}{\theta_{#1}}

\NewDocumentCommand\wt{O{t}}{\omega_{#1}}
\NewDocumentCommand\epxt{O{t}O{x}}{\ep^{#1}_{#2}}
\NewDocumentCommand\epuxt{O{t}O{\underline{x}}}{\ep^{#1}_{#2}}
\NewDocumentCommand\eppt{O{t}}{\epxt[#1][p]}
\def \wtt {\wt[t-1]}

\def \Nep {N_\ep}
\def \Lf {L_f}

\NewDocumentCommand\Lfi{O{}}{L_{f_{#1}}}

\def \Mtil {\tilde{M}}
\def \Mbar {\bar{M}}
\NewDocumentCommand\Mtilt{O{t}}{\Mtil_{#1}}

\NewDocumentCommand\MtU{O{t}}{M_{{#1}}}
\NewDocumentCommand\bMtU{O{t}}{\bar{M}_{{#1}}}

\NewDocumentCommand\MPiu{O{U^*}}{M_{\Pi, #1}}
\NewDocumentCommand\MPi{O{i}}{M_{\Pi_{#1}}}

\NewDocumentCommand\DX{O{}}{D_{X^{#1}}}

\NewDocumentCommand\Rt{O{t}}{R_{{#1}}}
\NewDocumentCommand\gi{O{i}}{\SubAddi{g}{#1}}

\NewDocumentCommand\eptl{O{}}{\tilde{\ep}_{#1}}

\NewDocumentCommand\MXtlt{O{t}}{M_{\tilde{X}}^{#1}}

\NewDocumentCommand\MPitlt{O{t}}{M_{\tilde{\Pi}}^{#1}}
\NewDocumentCommand\fistar{O{i}}{f^*_{#1}}
\NewDocumentCommand\fundert{O{N}}{\underline{f}^{#1}}
\NewDocumentCommand\Fundert{O{N}}{\underline{F}^{#1}}
\NewDocumentCommand\gistar{O{i}}{g^*_{#1}}
\NewDocumentCommand\fitilstar{O{i}}{\tilde{f}^*_{#1}}
\NewDocumentCommand\Dfistar{O{i}}{D_{\fistar[#1]}}
\NewDocumentCommand\Ufstar{}{U_{f^*}}

\NewDocumentCommand\fistarp{O{i}}{(f^*_{#1})'}
\NewDocumentCommand\ftlit{O{t}O{i}}{\tilde{f}_{#2}^{#1}}
\NewDocumentCommand\futit{O{t}O{i}}{\utilde{f}_{#2}^{#1}}

\NewDocumentCommand\myfi{O{\bci}}{f_{#1}}

\NewDocumentCommand\vt{O{t}}{\nu^{#1}}
\NewDocumentCommand\vit{O{t}O{i}}{\nu_{#2}^{#1}}

\NewDocumentCommand\nut{O{t}}{\nu^{#1}}
\NewDocumentCommand\nui{O{i}}{\nu_{#1}}
\NewDocumentCommand\nuit{O{t}O{i}}{\nu_{#2}^{#1}}
\NewDocumentCommand\nuitt{O{t-1}O{i}}{\nu_{#2}^{#1}}
\NewDocumentCommand\nuti{O{t}O{i}}{\nu_{#2}^{#1}}
\NewDocumentCommand\nutibar{O{t}O{i}}{\bar{\nu}_{#2}^{#1}}
\NewDocumentCommand\nubarti{O{t}O{i}}{\bar{\nu}_{#2}^{#1}}
\def \nubar {\bar\nu}
\def \nuitt {\vit[t-1]}
\def \hnu {\hat{\nu}}
\def \nustar {\nu^*}

\NewDocumentCommand\Ai{O{i}} {A_{#1}}
\NewDocumentCommand\Aitr{O{i}}{A_{#1}^\top}

\def \zt {\addt{z}}
\def \hz {\hat{z}}
\def \zbar {\bar{z}}
\NewDocumentCommand\zbart{O{N}}{\zbar^{#1}}
\def \zstar {z^*}

\NewDocumentCommand\qts{O{s} O{t}}{q_{#1}^{#2}}
\NewDocumentCommand\betat{O{t}}{\beta^{#1}}
\NewDocumentCommand\betast{O{s} O{t}}{\beta_{#1}^{(#2)}}
\NewDocumentCommand\gamt{O{t}}{\gamma^{#1}}
\NewDocumentCommand\gams{O{s}}{\gamma_{#1}}
\NewDocumentCommand\gamst{O{s} O{t}}{\gamma_{#1}^{(#2)}}
\NewDocumentCommand\delst{O{s} O{t}}{\delta_{#1}^{(#2)}}
\NewDocumentCommand\delt{O{t}}{\delta^{#1}}
\NewDocumentCommand\lamst{O{s}O{t}}{\lam_{#1}^{(#2)}}
\NewDocumentCommand\lamist{O{s}O{t}}{\lam_{#1, i}^{(#2)}}
\NewDocumentCommand\htilst{O{s}O{t}}{\tilde{h}^{(#2), #1}}
\NewDocumentCommand\wst{O{s}O{t}}{\tilde{w}^{(#2)}_{#1}}
\NewDocumentCommand\yst{O{s}O{t}}{y_{#1}^{(#2)}}
\NewDocumentCommand\rhost{O{s}O{t}}{\rho^{(#2)}_{#1}}

\NewDocumentCommand\Wt{O{t}}{W_{#1}}
\NewDocumentCommand\Gamt{O{t}}{\Gamma_{#1}}
\def \lam {\lambda}
\def \lambar {\bar{\lam}}

\def \rplus {r^+}
\NewDocumentCommand\lamt{O{t}}{\lam^{#1}}
\NewDocumentCommand\lambart{O{t}}{\bar{\lam}^{#1}}
\NewDocumentCommand\lamtilt{O{t}}{\tilde{\lam}^{#1}}

\def \lamstar {\lam^*}
\NewDocumentCommand\lami{O{i}}{\lam_{#1}}
\NewDocumentCommand\lamti{O{t}O{i}}{\lam_{#2}^{#1}}

\NewDocumentCommand\Lamr{O{r}}{\Lambda_{#1}}
\def \hlam {\hat{\lam}}
\def \lamhat {\hat{\lam}}

\def \Lam {\Lambda}

\NewDocumentCommand\kar{O{r}}{\kappa_{#1}}
\NewDocumentCommand\Llam{O{\Lam}}{{L}(#1)}
\NewDocumentCommand\Llamr{O{r}}{\Llam[\Lamr[#1]]}
\NewDocumentCommand\rlamr{O{r}}{d(\Lamr[#1])}
\def \Ltil {\tilde{L}}
\def \Ltilagg {\tilde{L}_{\text{agg}}}
\def \Lagg {L_{\text{agg}}}
\def \Funder {\underline{F}}
\def \FunderPD {\underline{F}_{\text{PD}}}
\def \DeltaFP {\Delta_{\text{FP}}}
\def \DeltaPD {\Delta_{\text{PD}}}
\def \half {\tfrac{1}{2}}
\def \rtil {\tilde{r}}

\NewDocumentCommand\bpi{O{}}{\bar{\pi}_{\bci[#1]}}
\NewDocumentCommand\pit{O{t}O{}}{\pi_{#2}^{#1}}
\NewDocumentCommand\piutt{O{t}O{}} {\utilde{\pi}_{#2}^{#1}}
\NewDocumentCommand\utpit{O{t}O{}} {\utilde{\pi}_{#2}^{#1}}
\NewDocumentCommand\pitt{O{}} {\pit[t-1][#1]}
\NewDocumentCommand\pibrt{O{t}}{\{\pi^{#1}\}}
\NewDocumentCommand\hpi{O{}}{\hat{\pi}_{\bci[#1]}}
\NewDocumentCommand\pistar{O{}}{\pi^*_{\bci[#1]}}
\NewDocumentCommand\hpiN{O{}}{\hat{\pi}^N_{\bci[#1]}}

\NewDocumentCommand\pibar{O{}}{\bar{\pi}_{#1}}
\NewDocumentCommand\pibart{O{t}}{\bar{\pi}^{#1}}

\def \uy {\underline{y}}

\NewDocumentCommand\yt{O{t}}{y^{#1}}

\NewDocumentCommand\yit{O{t}O{i}}{y^{#1}_{#2}}
\NewDocumentCommand\ytlitj{O{t}O{i}O{j}}{\tilde{y}^{#1}_{#2\setminus\{#3\}}}
\NewDocumentCommand\yutlit{O{t}O{i}}{\utilde{y}^{#1}_{#2}}
\NewDocumentCommand\utyit{O{t}O{i}}{\utilde{y}^{#1}_{#2}}
\NewDocumentCommand\ytlit{O{t}O{i}}{\tilde{y}^{#1}_{#2}}
\NewDocumentCommand\uyit{O{t}O{i}}{\uy_{#2}^{#1}}

\def \blx {\pmb{x}} 
\NewDocumentCommand\blxt{O{t}}{\blx^{#1}}
\def \xstar {x^*}
\NewDocumentCommand\xt{O{t}}{x^{#1}}
\def \xtt {\addtt{x}}

\def \xunder {\underline{x}}
\NewDocumentCommand\xundert{O{t}}{\xunder^{#1}}
\NewDocumentCommand\blxundert{O{t}}{\underline{\blx}^{#1}}
\NewDocumentCommand\xunderit{O{t}O{i}}{\xunder^{#1}_{#2}}
\NewDocumentCommand\xunderavt{O{t}O{i}}{\xunder^{#1}_{\ave}}

\NewDocumentCommand\myxi{O{i}}{x_{#1}}
\NewDocumentCommand\xit{O{t}O{i}}{x^{#1}_{#2}}
\NewDocumentCommand\xitt{O{i}}{\xit[t-1][#1]}
\NewDocumentCommand\xavt{O{t}}{\xit[#1][\ave]}

\def \hx {\hat{x}}
\NewDocumentCommand\hxit{O{t}O{i}}{\hat{x}^{#1}_{#2}}
\NewDocumentCommand\hblxit{O{t}O{i}}{\hat{\pmb{x}}^{#1}_{#2}}

\NewDocumentCommand\xtilt{O{t}}{\tilde{x}^{#1}}
\NewDocumentCommand\blxtlt{O{t}}{\tilde{\blx}^{#1}}
\NewDocumentCommand\xtlit{O{t}O{i}}{\tilde{x}^{#1}_{#2}}
\NewDocumentCommand\xtlavt{O{t}}{\xtlit[#1][\ave]}

\def \xbar {\bar{x}}
\def \xbarN {\xbart[N]}
\NewDocumentCommand\bxit{O{t}O{i}}{\bar{x}^{#1}_{#2}}

\NewDocumentCommand\xbart{O{t}}{\bar{x}^{#1}}

\NewDocumentCommand\dit{O{t}O{i}}{d^{#1}_{#2}}
\NewDocumentCommand\Nit{O{t}O{i}}{\mathcal{N}^{#1}(#2)}
\NewDocumentCommand\Nitj{O{t}O{i}O{j}}{\mathcal{N}^{#1}(#2;\{#3\})}
\NewDocumentCommand\Ni{O{i}}{\mathcal{N}(#1)}
\NewDocumentCommand\Hij{O{i,j}}{H_{#1}}
\NewDocumentCommand\Wij{O{i,j}}{W_{#1}}
\NewDocumentCommand\sumNi{O{i}}{\tsum_{j \in \Ni[#1]} \Hij[j, #1]}

\NewDocumentCommand\calWt{O{N}}{\mathcal{W}_{#1}}

\def \blw {\pmb{w}}
\NewDocumentCommand\blwt{O{t}}{\blw^{#1}}
\NewDocumentCommand\ublwt{O{t}}{\underline{\blw}^{#1}}
\NewDocumentCommand\bblwt{O{t}}{\bar{\blw}^{#1}}
\NewDocumentCommand\tilwt{O{t}}{\tilde{\omega}^{#1}}

\NewDocumentCommand\Q{O{}}{Q_{#1}}

\NewDocumentCommand\delx{O{x}}{\delta_{#1}}
\NewDocumentCommand\delpit{O{t}}{\delta^{#1}_\pi}
\NewDocumentCommand\delpt{O{t}}{\delta^{#1}_p}
\NewDocumentCommand\delxt{O{t}}{\delta_x^{#1}}

\NewDocumentCommand\uit{O{t}}{u^{#1}_i}

\def \calM {\mathcal{M}}
\NewDocumentCommand\Memt{O{t}}{\calM_{#1}}
\NewDocumentCommand\Mit{O{t} O{i}}{\calM_{#2, #1}}
\NewDocumentCommand\Mpiitl{O{t} O{l} O{i}}{\calM^{\pi, #2}_{{#3},#1}}
\NewDocumentCommand\Mpiit{O{t} O{i}}{\calM^\pi_{{#2}, #1}}
\NewDocumentCommand\Mpiitcp{O{t} O{i}}{\calM^{\pi, \text{cp}}_{{#2}, #1}}
\NewDocumentCommand\Mst{O{t} O{s}}{\calM_{#2, #1}}
\NewDocumentCommand\Mitcp{O{t} O{i} O{}}{\Mit[#1][#2]^{\text{cp} #3}}
\NewDocumentCommand\Mitcm{O{t} O{i}}{\Mit[#1][#2]^{\text{comm}}}
\NewDocumentCommand\Mstcp{O{t} O{s} O{}}{\Mit[#1][#2]^{\text{cp} #3}}
\NewDocumentCommand\Mstcm{O{t} O{s} }{\Mit[#1][#2]^{\text{comm}}}
\NewDocumentCommand\Ki{O{i}}{\mathcal{K}_{#1}}
\NewDocumentCommand\Si{O{i}}{\mathcal{S}_{#1}}
\def \Span {\text{span}}
\NewDocumentCommand\prox{O{}}{\text{prox}_{#1}}
\NewDocumentCommand\Kl{O{l}} {K_{#1}}
\NewDocumentCommand\Mt{O{t}}{\calM_{#1}}
\def \calK {\mathcal{K}}
\def \calF {{\mathcal{F}}}

\def \Fstar {F_*}

\def \skipdisplay {\setlength{\abovedisplayskip}{7pt}
\setlength{\belowdisplayskip}{7pt}
\setlength{\abovedisplayshortskip}{7pt}
\setlength{\belowdisplayshortskip}{7pt}}

\NewDocumentCommand\Lfp{O{p}}{L_{f, #1}}
\def \Ntil {\tilde{N}}

\def \Nep {N_\ep}
\def \Cep {C_\ep}
\def \Htil {\tilde{H}}
\NewDocumentCommand\Htilt{O{l}}{\Htil_{#1}}
\def \calR {\mathcal{R}}
\NewDocumentCommand\calRt{O{t}}{\calR_{#1}}
\journalname{Mathematical Programming}
\renewcommand{\subclassname}{\bfseries Mathematics Subject Classification (2020)\enspace}
\newcommand{\bbr}{\mathbb{R}}
\def\endproof{{\ \hfill\hbox{%
      \vrule width1.0ex height1.0ex
    }\parfillskip 0pt}\par}

\title{ Solving Convex Smooth Function Constrained Optimization Is Almost As Easy As Unconstrained Optimization
 \thanks{This work is partially supported by the ONR grant N00014-20-1-2089. 
}
 \thanks{This paper contains an almost parameter-free extension to the original version that was submitted to \textit{Mathematical Programming} in October 2022.}
}
\author{Zhe Zhang \and Guanghui Lan}
\institute{Zhe Zhang \at
    Edwardson School of Industrial Engineering,
    Purdue University, West Lafayette, IN 47907
    \email{zhan5111@purdue.edu}
    \and
    Guanghui Lan \at
    H. Milton Stewart School of Industrial \& Systems Engineering,
    Georgia Institute of Technology, Atlanta, GA 30332
    \email{george.lan@isye.gatech.edu}}
\date{August 25, 2026}
 
\begin{document}

%


\maketitle

\begin{abstract}

    While Nesterov's Accelerated Gradient Descent (AGD) efficiently solves constrained problems when the constraint set $X \subseteq \bbr^n$ is simple and easy to project onto, it remains an open question whether function-constrained problems $\min_{x \in X} \{F(x) : g(x) \leq 0\}$ can be solved as efficiently as unconstrained problems in terms of oracle complexity. We provide an affirmative answer by proposing the Accelerated Constrained Gradient Descent (ACGD) method, a single-loop algorithm that modifies AGD by replacing the descent step with a constrained descent step, adding only a few linear constraints to the prox mapping. ACGD achieves nearly the same oracle complexity as minimizing the optimal Lagrangian function (with the multiplier fixed at its optimal value). We establish matching lower bounds, demonstrating these complexity results are unimprovable. For large-scale problems with many constraints, we introduce ACGD-S, which replaces the computationally demanding constrained descent step with basic matrix-vector multiplications, maintaining optimal oracle and computation complexities. Together, these methods provide a nearly complete characterization of the hardness of smooth function-constrained optimization. Under the additional assumption that $X$ is bounded, we also propose adaptive versions that achieve an order-optimal oracle complexity without requiring knowledge of either the aggregate smoothness constant or the diameter of $X$; only the strong convexity modulus is required in the strongly convex case. We present encouraging numerical results demonstrating their efficiency.

\end{abstract}
\section{Introduction}\label{sec:intro}
Consider convex smooth function-constrained optimization of the form
\begin{equation}\label{eq:prob}
    \begin{split}
        \min_{x \in X}\  & \{F(x) :=f(x) + u(x)\} \\
        s.t.\            & g(x) \leq 0,
    \end{split}
\end{equation}
where both $f:\R^n \rightarrow \R$ and $g := [g_1, g_2, \ldots, g_m]^\top:\R^n \rightarrow \R^m$ are blackbox convex functions with Lipschitz continuous gradients, and the domain $X$ is convex and closed. The regularization function $u(x)$ is assumed to be both simple \cite{nemirovski2001lectures} and $\alpha$-strongly convex for some $\alpha\geq 0$.
Referred to as composite objectives \cite{nesterov2013gradient}, $F$ covers black-box optimization as a special case. For example, an $\alpha$-strongly convex objective function $\tilde{f}$ can be decomposed into $f(x):=\tilde{f}(x)-\frac{\alpha}{2}\normsq{x}$ and $u(x):=\frac{\alpha}{2}\normsq{x}$.
Problems of this type find a wide range of applications in, for example, the Neyman-Pearson classification problem \cite{rigollet2011neyman}, the fairness-constrained classification \cite{zafar2017fairness}, and the risk-constrained portfolio optimization \cite{gandy2005portfolio}. Since the dimensions $n$ and $m$ are large in many applications, we focus on first-order methods to find an approximate solution. Specifically, $(f, g)$ is assumed to be accessible only via a black-box first-order oracle, which returns $(f(x), g(x); \grad f(x), \grad g(x))$ when queried at some $x \in \R^n$, and the goal is to find an $(\ep; \ep/c)$-optimal solution (or, in short, an $\ep$-solution):
\begin{equation}\label{eq:opt-criteria}
    F(x^N) - F(\xstar) \leq \ep \text{ and} \norm{[g(x^N)]_+} \leq \ep/c,\end{equation}
where the scaling constant $c \geq 1$ represents the modeler's preference for constraint violation relative to sub-optimality.

This paper intends to develop fast methods and to understand the requisite computation cost. To ensure practical efficiency, special attention is paid to simple methods involving only oracle evaluations, projections onto $X$, and basic vector operations like matrix-vector multiplication. 
In particular, we consider two kinds of cost involved in solving \eqref{eq:prob}: a) the oracle complexity, i.e., the total number of queries to the first-order oracle, and b) the computation complexity, i.e., the total number of matrix-vector multiplications (each costing at most $O(mn)$ FLOPs). During implementation, these complexities translate to different computation burdens, so the dominating cost depends on the context. For instance,
if the constraint $g$ is complicated, say a large finite-sum function in the fairness-constrained problem, the gradient evaluation could be the bottleneck. On the other hand, if both $n$ and $m$ are large and $g$ is an affine function, the matrix-vector multiplication might be the bottleneck. Thus the ideal optimization method should excel in both directions.

\begin{table}
    \caption{Ideal Complexity for Solving \eqref{eq:prob}}\label{tb:ideal-complexity}
    \centering
    \jim{
        \begin{tabular}{|l| l | l|}
            \hline
            Case                           & Oracle Complexity      & Computation Complexity \\
            \hline
            Non-strongly Convex $\alpha=0$ & $\Theta(1/\sqrt{\ep})$ & $\Theta(1/\ep)$        \\
            \hline
            Strongly Convex $\alpha > 0$   & $\Theta(\log(1/\ep))$  & $\Theta(1/\sqrt{\ep})$ \\
            \hline
        \end{tabular}
    }
\end{table}

Since convex optimization has been studied extensively, e.g. \cite{nesterov2003introductory}, we can conjecture the best possible complexities for \eqref{eq:prob} by reducing it to simpler cases for which the optimal results are available. Recall from \cite{nemirovsky1983problem,nesterov1983method}, the number of oracle evaluations required to find an $\ep$-optimal solution to a smooth problem without function constraints is $\Theta(1/\sqrt{\ep})$ in general, and $\Theta(\sqrt{\kappa}\log(1/\ep))$ if the problem is also strongly convex, where $\kappa$ denotes the condition number.
One way to reduce \eqref{eq:prob}  is to consider the Lagrange dual formulation:
\begin{equation}\label{eq:La}
    \min_{x\in X} \max_{\lambda \in \R^m_+} \{\La(x; \lambda):=f(x) + u(x) + \lambda^\top (g(x))\}.\end{equation}
If the multiplier $\lam$ is fixed to the optimal dual multiplier $\lamstar$, the optimal Lagrangian function $\La(x; \lamstar)$ is smooth and does not have any function constraint. Since some minimizer to $\La(x; \lamstar)$ is the solution to \eqref{eq:prob}, in the ideal case, optimizing \eqref{eq:prob} should have the same oracle complexity as optimizing $\La(x; \lamstar)$, i.e., $\bigO(1/\sqrt{\ep})$ for the non-strongly convex case ($\alpha =0$) and   $\bigO(\sqrt{\kappa}\log(1/\ep))$ for the strongly convex case ($\alpha > 0$), where $\kappa$ is the condition number associated with $\La(x; \lamstar)$. On the other hand, inspired by \cite{nesterov2003introductory} and \cite{nemirovsky1991optimality}, Xu and Ouyang constructed {novel} large-scale linearly constrained quadratic programs in \cite{ouyang2021lower} to show the tight computation complexity for the linearly constrained smooth problem is $\Theta(1/\ep)$ if $\alpha =0 $, and $\Theta(1/\sqrt{\ep})$ if  $\alpha > 0$.  Since the nonlinear constraint function includes the affine function as a special case,  the computation complexity for \eqref{eq:prob} should be at least as expensive as the affine case. These conjectured complexities, summarized in Table \ref{tb:ideal-complexity}, lead naturally to the research question:
\vgap
\begin{center}
    \fbox{%
        \parbox{0.8\textwidth}{%
            Can we solve function-constrained problems with the same oracle complexities as those without function constraints, while maintaining the same computation complexities as solving linearly constrained problems?
        }%
    }
\end{center}
\vgap

However, despite much research effort on the subject from different directions, the question remains open. Broadly speaking, the current results can be grouped according to reformulations (see the summary in Table \ref{tb:lit_complexity}). The methods based on the Lagrangian formulation \eqref{eq:La} are usually simple to implement and have a low per-iteration cost. For instance, the ConEx method in \cite{boob2022stochastic} and the APD method in \cite{aybat2021primal} are single-loop algorithms with optimal computation complexities, but their oracle complexities are worse than the ideal ones by an order of magnitude. The current best method based on the augmented Lagrangian reformulation \cite{xu2020first}, a three-loop algorithm constructed from the inexact augmented Lagrangian method, the ellipsoid method, and the accelerated gradient descent (AGD) method, is closer to the ideal oracle complexities. However, the proposed method in \cite{xu2020first} scales poorly with the number of constraints. In fact, the method is advantageous to the APD method only when $m \leq 5$ in the numerical experiments \cite{xu2020first}.

Another line of research \cite{lin2018level,nesterov2018lectures}, called the level-set method, reformulates \eqref{eq:prob} to a root-finding problem associated with a certain mini-max problem. These methods typically require all constraint functions to share the same strong convexity modulus as the objective function. Such an assumption can be quite restrictive because it is violated when there exists one affine constraint. Assuming the uniform strong-convexity condition holds, the method proposed by Nesterov in Section 2.3.5 of \cite{nesterov2018lectures} can achieve an $\bigO(\log(1/\ep))$ oracle complexity when $\alpha >0$. However, the method might be too computationally demanding to implement for the large-scale setting because it requires the exact solution to a quadratic program (QP) in each iteration, and the exact solution to a quadratic-constrained quadratic program (QCQP) from time to time. Lin et al.\ relax the expensive computation oracle assumption in \cite{lin2018level}, but the oracle complexity also becomes worse by an order of magnitude. To sum up, there exist two major deficiencies: a) current methods fail to match the ideal oracle complexity, and b) methods that are close to the ideal oracle complexity are impractical in the large-scale setting.

\begin{table}
    \caption{Complexities for Smooth Constrained Optimization \eqref{eq:prob}}\label{tb:lit_complexity}
    \centering
    \begin{threeparttable}
        \begin{tabular}{|l| r r |r r| l |}
            \hline                                         & \multicolumn{2}{c|}{Strongly Convex $\alpha > 0$}    & \multicolumn{2}{c|}{Convex $\alpha =0$}              & Strong Oracle \&                                                                                 \\
            \cline{1-1}\cline{2-2}\cline{3-3}\cline{4-4}\cline{5-5}
            Method                                         & Oracle                                               & Computation                                          & Oracle                                             & Computation        & Assumption             \\
            \hline
            Level Set\tnote{1} \cite{nesterov2018lectures} &                                                      &                                                      & $\bigO(\tfrac{1}{\sqrt{\ep}}\log(\tfrac{1}{\ep}))$ & \-                 & QCQP \& QP Oracle      \\
            \hline
            Level Set\tnote{1} \cite{lin2018level}         &                                                      &                                                      & $\bigO(1/\ep)$                                     & $\bigO(1/\ep)$     &                        \\
            \hline
            iALM \cite{xu2020first}                        & $\bigO(m\log^3(\tfrac{1}{\ep}))$                     &
            $\bigO(m\log^3(\tfrac{1}{\ep}))$               & $\bigO(\tfrac{m}{\sqrt{\ep}}\log^3(\tfrac{1}{\ep}))$ & $\bigO(\tfrac{m}{\sqrt{\ep}}\log^3(\tfrac{1}{\ep}))$ & $m=\bigO(1)$ constraints                                                                         \\
            \hline
            APD \cite{aybat2021primal}                     & $\bigO({1}/{\sqrt{\ep}})$                            & $\bigO({1}/{\sqrt{\ep}})$                            & $\bigO({1}/{\ep})$                                 & $\bigO({1}/{\ep})$ &                        \\
            \hline
            ConEx \cite{boob2022stochastic}                & $\bigO({1}/{\sqrt{\ep}})$                            & $\bigO({1}/{\sqrt{\ep}})$                            & $\bigO({1}/{\ep})$                                 & $\bigO({1}/{\ep})$ &                        \\
            \hline
            ACGD [*]                                       & $\bigO(\log(1/\ep))$                                 &                                                      & $\bigO(1/\sqrt{\ep})$                              &                    & {\makecell{QP Oracle}} \\
            \hline
            ACGD-S [*]                                     & $\bigO(\log(1/\ep))$                                 & $\bigO(1/\sqrt{\ep})$                                & $\bigO(1/\sqrt{\ep})$                              & $\bigO(1/\ep)$     &                        \\
            \hline \hline
        \end{tabular}
        \begin{tablenotes}
            \item[1] If, in addition to $f$, all constraint functions, $g_1,\ldots,g_m$, are $\alpha$-strongly convex, the oracle complexity is $\bigO(\log(1/\ep))$ for the level-set method in \cite{nesterov2018lectures}, and is  $\bigO(1/\sqrt{\ep})$ for \cite{lin2018level}.
        \end{tablenotes}
    \end{threeparttable}
\end{table}

In this paper, we provide an affirmative answer to the research question by proposing efficient algorithms to simultaneously achieve both ideal complexities. An important observation for our development is that the ideal oracle complexity and the ideal computation complexity, shown in Table \ref{tb:ideal-complexity}, have different orders of magnitude. This explains why the single-loop algorithms which carry out $\bigO(1)$ matrix-vector operations per oracle evaluation, for example, the APD method \cite{aybat2021primal} and the ConEx method \cite{boob2022stochastic}, can only achieve the ideal computation complexity, but not the ideal oracle complexity. This also explains why the more complicated algorithms with better oracle complexities, for example, the level set method in \cite{nesterov2018lectures} and the iALM method in \cite{xu2020first}, require stronger computation oracles than basic matrix-vector operations. Following this observation, we first ask if the ideal oracle complexity in Table \ref{tb:ideal-complexity} is attainable with any strong computation oracle. The question leads us to the oracle-efficient Accelerated Constrained Gradient Descent (ACGD).

The ACGD method is based on the Lagrangian formulation \eqref{eq:La}. The method can be motivated by how Nesterov adapts the AGD method to solve a minimax problem of the form $\min_{x}\max_{i}\{f_1(x), ..., f_m(x)\}$ in Section 2.3 of \cite{nesterov2018lectures}.  Since the Lagrangian function is also a minimax problem, it also satisfies the following max-type smoothness condition (see Lemma 2.3.1 in \cite{nesterov2018lectures}):
\begin{equation}
    \max_{\lam \in \R^m_+} \La(x, \lam) - \max_{\lam \in \R^m_+} \{l_f(x;\xbar) + \jim{\sumi \lami l_{g_i}(x; \xbar) + u(x)}\} \leq \tfrac{L(\R^m_+)}{2} \normsq{x -\xbar},
\end{equation}
where $l_f$ and $l_g$ denote the linearizations of $f$ and $g$ at $\xbar$, respectively, and the aggregate smoothness constant $L(\R^m_+)$ denotes the maximum, over $\lambda \in \R^m_+$, of the Lipschitz smoothness constants of \jim{the smooth component $f(x)+\lambda^\top g(x)$ of $\La(x;\lambda)$.}
To obtain the desired oracle complexity, Nesterov proposes to modify the AGD method by replacing the descent step with a max-type descent step (see (2.3.12) in \cite{nesterov2018lectures} or the DRAO method in \cite{lan2022optimal}), which in our case becomes
\begin{align}
    \xt \leftarrow \argmin_{x \in X} \max_{\lam \in \R^m_+}  l_f(x; \xundert) + \sum_{i=1}^m \lami (l_{g_i}(x; \xundert)) + u(x) + \tfrac{\etat}{2} \normsq{x - \xtt}, \label{eq:prox-linear}
\end{align}
where the linearization center $\xundert$ is a certain convex combination of $\{\xt[0],\ldots, \xtt\}$. Since the maximization of $\lambda$ is over $\R^m_+$, \eqref{eq:prox-linear} is equivalent to the following \jim{linearly constrained prox-mapping}:
$$\xt \leftarrow \jim{\argmin}_{x \in X} \{l_f(x; \xundert)  + u(x) + \tfrac{\etat}{2} \normsq{x - \xtt} \text{s.t.}\   l_{g_i}(x; \xundert) \leq 0\ \forall i \in [m]\}.
$$
\jim{When $X$ is polyhedral and $u$ is quadratic, this constrained prox-mapping is a quadratic program. For general $X$ and $u$, its computational cost depends on the corresponding prox oracle.}
We call it the constrained descent step and the so-modified AGD method the Accelerated \textit{Constrained} Gradient Descent method.
Just like the AGD method, the stepsize parameter $\etat$ should be selected to be proportional to the aggregate smoothness constant $L(\R^m_+).$ Note, however, that one difficulty with this approach is $L(\R^m_+) = \infty$, so $\etat$ is infinite and the convergence rate is $\bigO(\infty/N^2)$, i.e., the method may not even converge. By taking a primal-dual perspective, we show $\etat$ only needs to be proportional to  the smoothness constant of $\La(x, \lamstar)$, denoted by $L(\lamstar)$. As a result, $\etat$ remains finite and the oracle complexity of the proposed ACGD method almost matches the optimal oracle complexity for optimizing $\La(x;\lamstar)$.
\jim{Importantly, when expressed in terms of the aggregate smoothness constant $L(\Lambda_c)$, this oracle complexity is optimal for solving the original function-constrained problem \eqref{eq:prob}, as demonstrated by the matching lower oracle-complexity bounds in Section \ref{sec:lower}.}

To enhance its efficiency for the large-scale setting with many constraints, we use the sliding technique \cite{lan2016gradient,lan2022optimal} to extend the ACGD method to the ACGD with Sliding (ACGD-S) method. Given a linearization center $\xundert$, instead of solving \eqref{eq:prox-linear} to optimality in each iteration, the ACGD-S method solves the bilinear saddle point inexactly by performing only a finite number of  $\lambda$-prox mappings and $x$-prox mappings (see \eqref{eq:X-proj-def}), the most expensive operation during which is matrix-vector multiplication. Particularly, after the $t$\ts{th} oracle evaluation, i.e., $\{\grad g(\xundert), \grad f(\xundert); f(\xundert), g(\xundert)\}$, the $x$-prox mapping and $\lambda$-prox mapping are repeated only $\bigO(t)$ times if $\alpha = 0$, and $\bigO(\sqrt{\theta^t})$ times for some $\theta > 1$ if $\alpha > 0$. The proposed ACGD-S method achieves the optimal computation complexity (matching the lower bound for linearly constrained problem \cite{ouyang2021lower}), while maintaining the same oracle complexity as the ACGD method. Therefore the ACGD-S method provides an almost complete characterization of both the computation and the oracle complexity for solving a smooth function-constrained problem. Moreover, the intricate step-size choice to achieve both the optimal $\bigO(\log(1/\ep))$ oracle evaluations (outer loops) and the optimal $\bigO(1/\sqrt{ \ep})$ matrix-vector multiplications (inner iterations) appears to be new to the sliding technique \cite{lan2016gradient,lan2021graph,lan2022optimal}, so it is of independent interest.

A practical limitation of ACGD and ACGD-S is that selecting their stepsizes requires the Lipschitz smoothness constant of $\La(x; \lamstar)$ (see \eqref{eq:La}), which, in turn, depends on the unknown optimal multiplier $\lamstar$.
\jim{Recent advances in parameter-free optimization that appeared after the initial release of this manuscript \cite{deng2026uniformly,lan2023optimal,lin2025adaptive} prompted us to investigate whether this dependence on unknown problem parameters could be eliminated.
Under the additional assumption that $X$ is bounded, we introduce adaptive variants that automatically tune their stepsizes without requiring knowledge of either the aggregate smoothness constant or the diameter $D_X$. The non-strongly convex variants require no problem parameters, whereas the strongly convex variants require only the strong convexity modulus.}
\jim{Adaptive ACGD roughly matches the oracle complexity of its non-adaptive counterpart order-wise, while adaptive ACGD-S achieves the same oracle complexity and incurs an additional multiplicative $\log$ factor in computational complexity due to restarts.}
The key innovation is a set of novel termination certificates that converge rapidly regardless of parameter misspecification, thereby enabling online detection of misspecified parameters. Numerical studies in Section \ref{sec:numerical} demonstrate their efficiency on large-scale, high-dimensional problems with many constraints.

\jim{Technically, the trilinear reformulation, algorithms, and analysis are inspired by the DRAO framework of~\cite{lan2022optimal}, which was developed for distributed risk-averse optimization. In adapting this framework to the constrained optimization setting, our key contributions are handling the unbounded dual feasible region $\lambda \in \R^m_+$ and developing a new termination certificate that enables an adaptive, nearly parameter-free implementation.}
The rest of the paper is organized as follows. Section \ref{sec:ANLA} proposes the ACGD method, and Section \ref{sec:lower} develops matching lower bounds. Section \ref{sec:ANLA-S} extends it to a computationally efficient ACGD-S method, Section \ref{sec:numerical} presents numerical results, and Section \ref{sec:conclusion} provides concluding remarks.



\subsection{Notations \& Assumptions}
The following notation and assumptions are used throughout the paper.
\begin{itemize}

    \item The set of optimal solutions to \eqref{eq:prob}, $X^*$, is nonempty, and $\xstar$ is an arbitrary optimal solution. $F_*$ denotes the optimal objective, $F(\xstar)$. 
    \item \jim{Unless otherwise stated, our complexity bounds concern the nontrivial accuracy regime in which the initial point does not already satisfy the stated termination criterion.  Any additional accuracy restrictions required by a particular result are stated explicitly.}
    \item \jim{We use the convention that $0/0=0$ unless otherwise specified .}
    \item The Fenchel conjugate (see \cite{Beck2017First}) of a convex function g(x) is defined as $g^*(\pi) := \max_{x \in \R^n} \inner{x}{\pi} - g(x).$
    \item $U_h$ denotes the Bregman distance function (see \cite{Beck2017First}) generated by a convex function $h$, i.e., $U_h(\pi; \pibar):= h(\pi) - h(\pibar) - \inner{h'(\pibar)}{\pi - \pibar}$ where $h'(\pibar)$ is some fixed subgradient in $\partial h(\pibar).$ If $g$ is vector-valued, i.e.,  $g(x):=[g_1(x), g_2(x),\ldots,g_m(x)]$, $U_g$ is vector-valued with its $i$\ts{th} component being the Bregman distance function generated by $g_i$, namely, $U_{g_i}$.
    \item We refer to the following computation as either a prox-mapping or a projection:
          \begin{equation}\label{eq:prox-def}
              \hat w \leftarrow \argmin_{w \in W} \inner{y}{w} + h(w) + \tau U(w;\bar w),\end{equation}
          where the vector $y$ represents some ``descent direction'' (the gradient for example), and
          $h(w)$ is a simple convex function \cite{LanBook}. $U$ is the Bregman distance function, $\bar w$ is a prox center, and $\tau$ is a stepsize parameter. Together they ensure the output $\hat w$ is close to $\bar w$. In particular, the following will be referred to as the $X$-projection:
          \begin{equation}\label{eq:X-proj-def}
              \hat x \leftarrow \argmin_{x \in X} \inner{y}{x} + u(x) + \tfrac{\tau}{2} \normsq{x - \xbar}.\end{equation}
          \jim{We assume that $X$ and $u$ are sufficiently simple for the $X$-projection in \eqref{eq:X-proj-def} to be computed efficiently. In the computation-complexity analysis of ACGD-S, we further assume that each such projection requires $\bigO(n)$ FLOPs and is therefore no more expensive than a matrix-vector multiplication.}
    \item \jim{We refer to an operation of the following form, which includes additional affine constraints, as a constrained prox-mapping:
              \begin{equation}\label{eq:constrained-prox-def}
                  \begin{aligned}
                      \hat x \leftarrow \argmin_{x \in X}\quad & \inner{y}{x} + u(x) + \tfrac{\tau}{2}\normsq{x-\xbar} \\
                      \text{s.t.}\quad                         & \inner{a_i}{x} \leq b_i,\quad \forall i\in[m].
                  \end{aligned}
              \end{equation}}
          \jim{The constrained prox-mapping reduces to a quadratic program when $X$ is polyhedral and $u$ is quadratic; for general $X$ and $u$, the cost of evaluating this mapping depends on the corresponding prox oracle.}
    \item \jim{A differentiable function $h$ is $L$-smooth if $\norm{\grad h(x)-\grad h(y)}\le L\,\norm{x-y}$ for all $x,y\in\R^n$. We use $L_f$ to denote the Lipschitz smoothness constant of the objective function $f$. Since a Lipschitz smoothness constant need not be minimal, we adopt without loss of generality the convention $L_f\geq\alpha$, replacing any valid constant $L_f$ by $\max\{L_f,\alpha\}$ if necessary. For the vector-valued constraint function $g(x)=[g_1(x),\dots,g_m(x)]$, let $L_{g_i}$ denote the Lipschitz smoothness constant of $g_i$, and define the aggregate Lipschitz smoothness constant of $g$ by
          \begin{equation}\label{eq:Lg-def}
              L_g = \| L_{g_1}, \dots, L_{g_m} \|_2.
          \end{equation}}
    \item \jim{For fixed $\lambda$, the Lipschitz smoothness constant of $\La(x;\lambda)$ with respect to $x$ refers to that of its smooth component $f(x)+\lambda^\top g(x)$; the proximal term $u(x)$ is excluded.}
    \item { We use $B^m(x;r)$ to denote the $m$-dimensional ball centered at $x$ with radius $r$, i.e., $B^m(x;r) = \{y \in \R^m: \norm{y-x} \leq r\}$. We also use $B^m_+(x;r)$ to denote the $m$-dimensional ball centered at $x$ with radius $r$ and non-negative elements, i.e., $B^m_+(x;r) = \{y \in \R^m: \norm{y-x} \leq r, y \geq 0\}$.}
\end{itemize}

\section{The Accelerated Constrained Gradient Descent Method}\label{sec:ANLA}

We present in this section the ACGD method. Specifically, Subsection \ref{subsec:ACGD} introduces a novel primal-dual perspective to motivate the ACGD method, Subsection \ref{subsec: ACGD-Results} presents the convergence results, \jim{and Subsection \ref{subsec:ANLA-conv-analysis} furnishes the detailed proof}. Subsection \ref{subsec:ANLA-search} proposes the guess-and-check scheme to look for the problem parameter $L(\lambda^*)$, and Subsection \ref{subsec:a-ANLA-conv-analysis} contains the detailed proofs.

\subsection{The ACGD method}\label{subsec:ACGD}
This subsection presents the primal-dual perspective to motivate the ACGD method. Such a perspective is important for understanding the finite stepsize parameter $\etat$ discussed in Section \ref{sec:intro}. Specifically, we first introduce a novel nested Lagrangian function which reformulates \eqref{eq:prob} as a $\min-\max-\max$ trilinear saddle point problem. To search for the saddle point, we propose a primal-dual type method similar to \cite{lan2022optimal,zhang2019efficient}. The ACGD method then simply follows from rewriting the proposed method in the primal form.

First, we need to assume the existence of a KKT point to the Lagrangian function defined in \eqref{eq:La} throughout this paper:
\begin{assumption}
    There exists a $\lamstar \in \R^m_+$ and $\xstar \in X$ such that ${-\left[\grad f(\xstar) + \sumi \lamstar_i\grad g_i(\xstar)+ u' \right]}\in N_X(\xstar)$ for some $u' \in \partial u(\xstar)$, $g(\xstar) \leq 0$, and $\lamstar_i g_i(\xstar) = 0\ \forall i$, where $N_X$ denotes the normal cone to $X.$\\
\end{assumption}
Note  $\xstar$ and $\lamstar$ could be interpreted as an arbitrary element in $X^*$ and $\Lam^*:= \argmax_{\lam \in \R^m_+} \min_{x\in X} \La(x; \lambda)$ for the rest of the paper,
because any $(\xbar, \bar{\lam})$ with  $\xbar \in X^*$ and $\lambar \in \Lam^* $ also constitutes a KKT point (see Proposition 3.4.1 in \cite{bertsekas2009convex}).

To motivate the nested Lagrangian function, consider the simplified case where the optimal dual multiplier $\lamstar$ is known. Fix $\lambda$ to $\lamstar$, an optimal solution $\xstar$ can be found by solving the following simplified problem under certain regularity conditions:
\begin{equation}\label{eq:unconstrained}
    \min_{x \in X} f(x) + \inner{\lamstar}{g(x)} +  u(x).\end{equation}
One useful framework for designing an optimal algorithm for the constrained problem in \eqref{eq:unconstrained} is to consider a bilinear saddle-point reformulation
\begin{equation}\label{eq:fenchel}
    \min_{x\in X} \max_{  \nu \in V, \pi \in \Pi} \{ \inner{x}{\pi} - f^*(\pi) + u(x) + \inner{\lamstar}{\nu x - g^*(\nu)} \},\end{equation}
where $f^*$ is the Fenchel conjugate function to $f$ and $\Pi$ is its domain, namely,  $\{\pi \in \R^n: f^*(\pi) < \infty \}$, and $g^*:=[g_1^*, \ldots, g_m^*]$ is the vector-valued Fenchel conjugate to $g$ and $V$ is its domain, namely, $\{v \in \R^{m\times n}: g_i^*(v_i) < \infty\ \forall i \in [m] \}.$
\jim{This trilinear reformulation would be important for the design of the sliding algorithm later in Section~\ref{sec:ANLA-S}.}
Moreover, since the dual variables $\pi$ and $\nu$ are associated with a common primal variable $x$, we use the following notation to denote the joint domain $[V, \Pi]$:
\begin{equation}\label{eq:joint-domain}
    [V, \Pi] = \{(\grad g(x), \grad f(x)): x\in \R^n\}.\end{equation}

However, $\lamstar$ is unknown in practice. So we propose to consider the following nested Lagrangian reformulation which combines \eqref{eq:La} and \eqref{eq:fenchel}:
\begin{equation}\label{eq:nLa}
    \min_{x\in X} \max_{\lambda \in \R^m_+, \pi \in \Pi, \nu \in V} \{\La(x; \lambda, \pi, \nu):= \inner{x}{\pi} - f^*(\pi) + u(x) + \inner{\lambda}{\nu x - g^*(\nu)}\},\end{equation}
where $f^*$, $g^*$, $\Pi$ and $V$ are defined in the same way as \eqref{eq:fenchel}. Notice a common notation $\La$ is used for the nested Lagrangian and the ordinary Lagrangian, but the exact meaning should be clear from the context. 
Let $Z$ denote the joint domain, $Z:= X \times \R^m_+ \times \Pi \times V$ and $\zstar:= (\xstar; \lamstar, \nustar :=\grad g(\xstar), \pistar:=\grad f(\xstar))$. It is not hard to see that $\zstar$ is a saddle point to \eqref{eq:nLa}, and a useful criterion to measure the optimality of an iterate $z^t=(\xt; \lamt, \vt, \pit ) \in Z$ is to compare it to some reference point $z \in Z$ in the following gap function:
\begin{equation}\label{eq:Q_def}
    Q(\zt; z):= \La(\xt; \lam, \nu, \pi) - \La(x; \lamt, \vt, \pit).
\end{equation}
Indeed, the saddle point $\zstar$ satisfies $Q(\zstar; z) \leq 0\ \forall z \in Z$.

A crucial observation for our development is that for convergence to an $\ep$-solution, it is sufficient to consider only the reference $\lam$'s inside a certain bounded set rather than the positive orthant.
The next lemma shows the $Q$ function still provides upper bounds for both the feasibility violation and the optimality gap.
\begin{lemma}\label{lm:fun-feas-from-Q}
    Let $\zt=(\xt; \lamt, \nut, \pit) \in Z$ be given and let $\Lamr$ denote a certain set of reference $\lambda$'s,
    \begin{equation}\label{eq:ref_lambs}
        \Lamr =\{0\} \cup \{\lambda^* + \lambda: \lambda \in B^m_+(0;r)\}.\end{equation}
    If $\max_{\lam \in \Lamr, (\nu,\pi) \in [V, \Pi]}Q(\zt; (\xstar; \lam, \nu, \pi)) \leq \ep$, we have $F(\xt) - F(\xstar) \leq \ep$ and $ \norm{[g(\xt)]_+} \leq \ep/r$.
\end{lemma}

\begin{proof}
    Fixing $\hpi = \grad f(\xt)$ and $\hnu = \grad g(\xt)$, the given condition implies that $Q(\zt; (\xstar; \lam, \hnu, \hpi)) \leq \ep\ \forall \lam \in \Lamr$. It then follows from the conjugate duality  \jim{$f(x^t)+ f^*(\hat{\pi})=\inner{x^t}{\hpi}$ and $g(x^t)+g^*(\hnu)=\inner{\hnu}{x^t}$} (see Section 4.1 and 4.2 of \cite{Beck2017First}) that
    $$f(\xt) + \inner{\lam}{g(\xt)} + u(\xt) - [f(\xstar) + \inner{\lamt}{g(\xstar)} + u(\xstar)] \leq \ep\ \forall \lam \in \Lamr.$$
    Since $\xstar$ is feasible, i.e., $g(\xstar)\leq 0$, we have $\inner{\lamt}{g(\xstar)} \leq 0$. Then taking $\lam = 0$ leads to
    $$F(\xt)- \Fstar = f(\xt) + u(\xt) - f(\xstar) - u(\xstar)\leq f(\xt) + u(\xt)  - [f(\xstar) + \inner{\lamt}{g(\xstar)} + u(\xstar)] \leq \ep.$$
    Next since $(\xstar; \lamstar)$ is a saddle point to \eqref{eq:La}, we have
    $0 \leq F(\xt) + \inner{\lamstar}{g(\xt)} - [F(\xstar) + \inner{\lamt}{g(\xstar)}]$. \jim{There are two cases to consider. If $\norm{[g(\xt)]_+} = 0$, then $g(\xt)\leq 0$ and we are done. In the following we assume $\norm{[g(\xt)]_+} > 0$.} Setting $\hlam = \lamstar+ r [g(\xt)]_+/ \norm{[g(\xt)]_+} \in \Lamr$, we get
    \begin{align*}
        r\norm{[g(\xt)]_+} & \leq F(\xt) +  \inner{\lamstar+ r[g(\xt)]_+/ \norm{[g(\xt)]_+}}{g(\xt)} - [F(\xstar) + \inner{\lamt}{g(\xstar)}] \\
                           & \leq Q(\zt; (\xstar; \hlam, \hnu, \hpi)) \leq \ep.
    \end{align*}
\end{proof}
We remark that a lower bound to $F(\xt) - \Fstar \geq -\epsilon \norm{\lamstar}/r$ can also be derived in similar fashion, see \cite{lan2013iteration,yang2022data}.

Now, let us move on to consider minimizing $Q$ to find a saddle point to \eqref{eq:nLa}. An essential feature of the nested Lagrangian function \eqref{eq:nLa} is the trilinear term $\inner{\lam}{\nu x - g^*(\nu)}$. The problem cannot be simplified to a $\min-\max$ saddle point problem by combining the $\nu$ and $\lam$ into a single dual block, because \jim{jointly maximizing a bilinear term without exploiting the problem structure is computationally intractable}.
Similar problems have been studied by Lan and Zhang in \cite{lan2022optimal,zhang2019efficient,zhang2020optimal} and the key to handle the $\min-\max-\max$ trilinear structure is to decompose the $Q$ function into the sub-gap functions and optimize each sub-gap function sequentially. Specifically, the following decomposition of \eqref{eq:Q_def} into sub-gap functions associated with $x$, $\lam$, $\nu$ and $\pi$ is useful.
\begin{equation}\label{eq:Q_decomp}
    Q(\zt; z) = Q_x(\zt; z) + Q_\lam(\zt; z) + Q_\pi(\zt; z) + Q_\nu(\zt; z)
\end{equation}
where
\begin{subequations}\label{eq:subeqns}
    \begin{align}
        Q_\pi(\zt; z)     & := \La(\xt; \lam, \nu, \pi) - \La(\xt; \lam, \nu, \pit) = \inner{\xt}{\pi} - f^*(\pi) \boxed{- [\inner{\xt}{\pit} - f^*(\pit)]}, \label{eq:Qpi}                                                   \\
        Q_\nu(\zt; z)     & := \La(\xt; \lam, \nu, \pit) - \La(\xt; \lam, \nut, \pit) = \sumi \lami \left[\inner{\nui}{\xt} - g^*_i(\nui)\right] \boxed{- \sumi \lami[i] [\inner{\nuit}{\xt} - g_i^*(\nuit)]}, \label{eq:Qnu} \\
        Q_\lambda(\zt; z) & := \La(\xt; \lam, \nut, \pit) - \La(\xt; \lamt, \nut, \pit) = \inner{\lam}{\nut \xt - g^*(\nut)} \boxed{- \inner{\lamt}{\nut \xt - g^*(\nut)}}, \label{eq:Qlam}                                   \\
        Q_x(\zt; z)       & := \La(\xt; \lamt, \nut, \pit) - \La(x; \lamt, \nut, \pit) = \boxed{\inner{\pit + \sumi \lamti \nuti}{ \xt} + u(\xt)} - \inner{\pit + \sumi \lamti \nuti}{ x} - u(x). \label{eq:Qx}
    \end{align}
\end{subequations}
Following the DRAO method in \cite{lan2022optimal}, at each iteration, we reduce the boxed terms associated with $Q_\pi$, $Q_\nu$, and the pair $(Q_x,Q_\lambda)$ using the respective prox-mappings \jim{\eqref{ls:anla-pi}, \eqref{ls:anla-nu}, and \eqref{ls:anla-x-lambda}}:
\begin{subequations}\label{ls:anla}
    \begin{align}
        \xtilt       & \leftarrow \xtt + \thetat (\xtt - \xt[t-2]); \label{ls:anla-extrapolation}                                                                                                 \\
        \pit         & \leftarrow \argmax_{\pi  \in \Pi} \inner{\pi}{\xtilt} - f^*(\pi) - \taut U_{f^*}(\pi; \pit[t-1]); \label{ls:anla-pi}                                                       \\
        \nuti        & \leftarrow  \argmax_{\nui \in V_i} \inner{\nui}{\xtilt} - g_i^*(\nui) - \taut U_{g_i^*}(\nui; \nuit[t-1])\ \forall i \in [m]; \label{ls:anla-nu}                           \\
        (\xt, \lamt) & \leftarrow \text{arg}\min_{x \in X} \max_{\lam \in \R^m_+} \inner{\pit}{x} + \inner{\lam}{\nut x - g^*(\nut)} + u(x) + \etat \normsq{x - \xtt}/2. \label{ls:anla-x-lambda}
    \end{align}
\end{subequations}
In the above listing, $\thetat$, $\taut$ and $\etat$ are non-negative stepsize parameters, and $U_{g_i^*}$ and $U_{f^*}$ are Bregman distance functions generated by $g_i^*$ and $f^*$ respectively. Particularly, with the yet  not available $\xt$ replaced by some proxy $\xtilt$, the above $\pit$ update corresponds to the minimization of the variable $\pit$ in \eqref{eq:Qpi} subject to a prox term $\taut U_{f^*}(\pi; \pit[t-1])$.
Similarly, the $\nut$ update corresponds to the minimization of the variable $\nut$ in \eqref{eq:Qnu} subject to a prox term $\taut \sumi \lami U_{g_i^*}(\nui; \nuit[t-1])$. Since the summation weights $\{\lami\}$ are non-negative and the maximizations are separable, the $\nut$ update is written equivalently as individual $\vit$ updates in \eqref{ls:anla}. Finally, since both the maximization of $\lamt$ in \eqref{eq:Qlam} and the minimization of $\xt$ in \eqref{eq:Qx} do not require any oracle information related to either $f$ or $g$, we evaluate them simultaneously subject to a prox term $\etat \normsq{x - \xtt}/2$. This leads to the joint $(\xt, \lamt)$ update in \eqref{ls:anla}.
We remark that $\xtilt$, the proxy for $\xt$, is chosen as the momentum extrapolation from $\xtt$. Such a choice helps us to achieve acceleration ( e.g. see Section 3.4 in \cite{LanBook} for the connection).

The implementable version of \eqref{ls:anla} is shown in Algorithm \ref{alg:ANLA}. It employs two additional simplifications. First, we initialize the dual variables to primal gradients, i.e.,  $\pit[0] = \grad f(\xt[0]) $ and $\nut[0] = \grad g(\xt[0])$. With $U_{g^*}$ and $U_{f^*}$ selected as prox-functions, we can show recursively that  $\pit$ and $\vt$ in \eqref{ls:anla} are the same as the gradients at some averaged point (see Lemma 2 in \cite{zhang2020optimal}). Thus the $\nuit$ and $\pit$ computation in \eqref{ls:anla} simplifies to gradient evaluations in Line 4 of Algorithm \ref{alg:ANLA}. Second,  the $(\xt, \lamt)$-saddle point problem in \eqref{ls:anla} is formulated as \jim{a constrained prox-mapping} in Line 5 of Algorithm \ref{alg:ANLA}. Since  $\nut = \grad g(\xundert)$ in Algorithm \ref{alg:ANLA} implies the relation $g_i(\xundert) + g_i^*(\nuit) = \inner{\xundert}{\nuit}$ (see Theorem 4.20 in \cite{Beck2017First}), we have $g^*(\nut)=\nut\xundert - g(\xundert)$. Interestingly, other than the additional linear constraint associated with $g$ in the descent step (Line 5), Algorithm \ref{alg:ANLA} is the same as Nesterov's AGD method \cite{nesterov1983method}. Therefore we name it the Accelerated \textit{Constrained} Gradient Descent method.

\jim{Before presenting the convergence analysis, we note that the ACGD updates can alternatively be obtained by applying a prox-linear step \cite{nesterov2003introductory} to the penalty function associated with \eqref{eq:prob},
    \[
        P(x;r):=f(x)+u(x)+ r\sum_{i=1}^{m}[g_i(x)]_+, \text{ with $r\to+\infty$.}
    \]
    However, as with the Lagrangian perspective discussed in the introduction, directly adopting this viewpoint would make the convergence analysis considerably more difficult.}

\begin{minipage}{.95\textwidth}
    \begin{algorithm}[H]
        \caption{\textbf{A}ccelerated \textbf{C}onstrained \textbf{G}radient \textbf{D}escent Method}
        \label{alg:ANLA}
        \begin{algorithmic}[1]
            \Require $\xt[-1]= \xundert[0] = \xt[0]\in X$,  stepsizes $\{\thetat\}$, $\{\etat\}$, $\{\taut\}$, and  weights $\{\wt\}$.
            \State Set $\pit[0] = \grad f(\xt[0]) $, $\nut[0] = \grad g(\xt[0])$.
            \For{$t =  1,2,3 ... N$}
            \parState{Set $\xundert \leftarrow (\taut \xundert[t-1] + \xtilt)/(1 + \taut)$ where  $\xtilt  = \xtt + \thetat (\xtt - \xt[t-2])$.}
            \parState{Set $\pit \leftarrow \grad f(\xundert)$ and $\nut \leftarrow \grad g(\xundert)$.}
            \parState{Solve $\xt \leftarrow \argmin_{x \in X} \{\inner{\pit}{x} +u(x) + \etat \normsq{x - \xtt}/2\ s.t.\ \nut (x - \xundert) + g(\xundert) \leq 0\} .$} 
            \EndFor
            \State \Return $\xbart[N]:= \tsum_{t=1}^N \wt \xt / (\sumt \wt). $
        \end{algorithmic}
    \end{algorithm}
    \vgap
\end{minipage}
\subsection{Convergence Analysis of the ACGD Method}\label{subsec: ACGD-Results}
Next, we present convergence results for the ACGD method. The detailed analysis is deferred to Subsection \ref{subsec:ANLA-conv-analysis}. The next proposition states some conditions required for the $Q$ function in \eqref{eq:Q_def} to converge.
Since these conditions depend on the set of reference multiplier $\lambda$'s under consideration (see Lemma \ref{lm:fun-feas-from-Q}), it is useful to define an aggregate Lipschitz smoothness constant as a function of $\Lambda$:
\jim{
    \begin{equation}\label{eq:agg-L}
        L(\Lambda):= \sup_{\lam \in \Lambda} L_\lambda, \text{ where } L_\lambda \norm{\xbar - \hx} \geq \norm{\grad f(\xbar) + \sumi \lambda_i \grad g_i(\xbar) - \grad f(\hx) - \sumi \lambda_i \grad g_i(\hx)} \quad \forall \xbar, \hx \in \R^n .
    \end{equation}
}
\begin{proposition}\label{pr:Q-convergence}
    Let a set of reference multipliers $\Lambda \subset \R^m_+$ be given and let the aggregate smoothness constant $\Llam$ be defined in \eqref{eq:agg-L}. Consider the iterates ${\zt:=(\xt; \lamt, \nut, \pit)}$ generated by Algorithm \ref{alg:ANLA}, where $\lamt$ is defined in \eqref{ls:anla}. Suppose the following conditions are satisfied by the stepsizes together with some non-negative weights $\wt \geq 0$ for all $t \geq 2$:
    \begin{align}
         & \wt \etat \leq \wtt (\etatt + \alpha), \label{stp-con:eta-recur}                                                                      \\
         & \wt \taut \leq \wtt (\tautt + 1), \label{stp-con:Q-recur}                                                                             \\
         & \etat[N](\taut[N] + 1) \geq \Llam, \etatt \taut \geq \theta_t\Llam \text{ with }  \thetat :=\wtt /\wt \leq 1\label{stp-con:Q-cancel}.
    \end{align}
    Consider the ergodic iterate $\zbart[N]= (\xbarN; \lambart[N], \nubarti[N][], \pibart[N])$ specified according to
    \begin{equation}\label{eq:ergo_solns}
        \begin{split}
             & \xbarN := \sumt \wt \xt / (\sumwt),\ \lambart[N]:= \sumwt \lamt/ (\sumwt),\ \pibart[N]:= \sumwt \pit / (\sumwt), \\
            \vspace{0.2cm}                                                                                                      \\
             & \nutibar[N]:=  \begin{cases}\sumwt \lamti \nuit /(\sumwt \lamti) & \text{o.w.}                       \\
             \grad g_i(x^0)                       & \text{if } \lamti = 0\ \forall t.\end{cases}
        \end{split}
    \end{equation}
    Then the following convergence bound is valid for any reference point $z=(x; \lam,  \nu, \pi) \in X \times \Lam \times [V, \Pi]$ with $[V, \Pi]$ defined in \eqref{eq:joint-domain}:
    \begin{equation}\label{eq:Q-conv}
        \begin{split}
            (\sumwt) & Q(\zbart, z) + \wt[N](\etat[N] + \alpha) \normsq{\xt[N]-x}/2                                                                                      \\
            \leq     & {\tsum_{t=1}^N \wt[t]\left(\La(\xt[t]; \lam, \nu, \pi) - \La(x;\lam^t,\nu^t,\pi^t)\right)+ \frac{\wt[N](\etat[N] + \alpha)}{2} \normsq{\xt[N]-x}} \\
            \leq     & \wt[1]\etat[1] \normsq{\xt[0] - x}/2 + \wt[1]\taut[1](\jim{U_{f^*}}(\pi; \pit[0]) + \inner{\lambda}{\jim{U_{g^*}}(\nu; \nut[0])}).
        \end{split}
    \end{equation}
\end{proposition}

Proposition \ref{pr:Q-convergence}, together with Lemma \ref{lm:fun-feas-from-Q}, show that it is possible to select a finite stepsize $\etat$ related to $L(\Lamr)$ in place of the infinite $L(\R^m_+)$.
To provide convergence guarantees for both the feasibility violation and the optimality gap associated with the function-constrained problem in  \eqref{eq:prob},
Lemma \ref{lm:fun-feas-from-Q} states that we should only consider the $Q$ gap function defined with respect to $\Lamr$, some small neighborhood of reference $\lambda$'s around $\lamstar$.  Proposition \ref{pr:Q-convergence} shows such an $Q$ function converges when the stepsize choice satisfies certain conditions related to $L(\Lamr)$, rather than $L(\R^m_+)$. Thus a finite $\etat$ proportional to $L(\Lamr)$ is sufficient for our purpose, and the next theorem states the result more precisely.


\begin{theorem}\label{thm:ANLA}
    Let a smooth constrained optimization problem \eqref{eq:prob} be given and let $\Llamr$ be defined in  \eqref{eq:ref_lambs} and \eqref{eq:agg-L}. Denote its condition number by $\kar= \Llamr/\alpha$ (We set $\kar = \infty$ if $\alpha=0$).  Suppose the solution iterates $\{\xt\}$ are generated by Algorithm \ref{alg:ANLA} with the following stepsizes for $t \geq 1$
    \begin{equation}\label{stp-choice:ANLA}
        \taut = \min\{\tfrac{t-1}{2}, \sqrt{\kar}\},\ \etat = \tfrac{\Llamr}{\taut[t+1]},\ \thetat = \tfrac{\taut}{\taut[t-1] + 1},\ \wt = \begin{cases} \wt[t-1] /\thetat & \text{if } t  \geq 2, \\
              1                 & \text{if } t = 1.
        \end{cases}
    \end{equation}
    Then the ergodic average $\xbarN$ solution satisfies
    \begin{equation}\label{eq:thm-non-strong}
        \max\{F(\xbarN) - \Fstar, r\norm{[g(\xbarN)]_+}\} \leq \frac{2\Llamr}{N(N+1)} \normsq{\xt[0] - \xstar}.
    \end{equation}
    Moreover, in the strongly convex case with $\alpha >0 $, $\xbarN$ also satisfies
    \begin{align}
         & \max\{F(\xbarN) - \Fstar, r\norm{[g(\xbarN)]_+}\} \leq \frac{ \sqrt{\Llamr \alpha} \normsq{\xt[0] - \xstar}}{(1+ 1/\sqrt{\kar})^{N-4} - 1}, \label{eq:thm-strong} \\
         & \normsq{\xbarN - \xstar} \leq \frac{2 \sqrt{\kar} \normsq{\xt[0] - \xstar}}{(1+ 1/\sqrt{\kar})^{N-4} - 1}. \label{eq:thm-strong-dist}
    \end{align}
\end{theorem}

As a consequence of the preceding theorem, we can derive upper bounds on the number of iterations required by the ACGD method to find an $(\ep; \ep/c)$-solution (see \eqref{eq:opt-criteria}). The next corollary focuses on the non-strongly convex case, i.e., $\alpha=0$.
\begin{corollary}\label{cor:ANLA-complexity-non-strong}
    Suppose $\{\xt\}$ are generated by Algorithm \ref{alg:ANLA} using the stepsizes choice in \eqref{stp-choice:ANLA}, with $\Llamr=\Llamr[c]$ (see \eqref{eq:agg-L}). The ergodic average solution $\xbarN$ is an $(\ep, \ep/c)$-solution if
    $$N \geq \sqrt{\tfrac{2\Llamr[c]}{\ep}} \norm{\xt[0]- \xstar}.$$
\end{corollary}

For the strongly convex problem with $\alpha > 0$, it is advantageous to choose a small $r$ to ensure a small aggregate Lipschitz smoothness constant $\Llamr$ and hence a small condition number $\kar$. We set it to $\Llamr[1]$ in the next corollary.
\begin{corollary}\label{cor:ANLA-complexity-strong}
    Suppose $\{\xt\}$ are generated by Algorithm \ref{alg:ANLA} using the stepsizes choice in \eqref{stp-choice:ANLA}, with $\Llamr=\Llamr[1]$ (see \eqref{eq:agg-L}). Then the ergodic average solution $\xbarN$ is an $(\ep, \ep/c)$-solution if
    $$N \geq  \min\{\sqrt{\tfrac{2\max\{c, 1\} \Llamr[1]}{\ep}} \norm{\xt[0] - \xstar}, [\sqrt{\tfrac{2\Llamr[1]}{\alpha}} + 1]\log[\tfrac{\max\{c,1\}\sqrt{\Llamr[1]\alpha}\normsq{\xt[0] - \xstar}}{\ep} + 1] + 4\}.$$
\end{corollary}
Note that for a small $\lamstar$, i.e., $\norm{\lamstar}$ being significantly less than $1$, it may be worthwhile to choose $\Llamr = \Llamr[\norm{\lamstar}]$ to further reduce the condition number from $\kar[1] = \Llamr[1]/\alpha$ to $\kar[\norm{\lamstar}] = \Llamr[\norm{\lamstar}]/\alpha$.
\vspace{1mm}

Since each iteration of the ACGD method requires only one gradient evaluation, the preceding two corollaries establish the desired $\bigO(1/\sqrt{\ep})$ and $\bigO(\log(1/\ep))$ oracle complexities for the non-strongly convex and the strongly convex problems, respectively.

\subsection{Convergence Analysis of the ACGD Method}\label{subsec:ANLA-conv-analysis}
We present the detailed proofs of convergence results for Algorithm \ref{alg:ANLA}. \\

\textbf{Proof of Proposition \ref{pr:Q-convergence}}
Let $Q_x$, $Q_\lam$, $Q_\nu$, and $Q_\pi$ be defined in \eqref{eq:Q_decomp}. It is useful to view the updates in Algorithm \ref{alg:ANLA} from the perspective of prox-mappings in \eqref{ls:anla}. First, let's consider $Q_\nu$. We have from the definition of $\xtilt$ that
\begin{align*}
    \inner{\nui -\nuit }{(\xtilt-\xt)}= & -\inner{\nui - \nuit}{(\xt -\xtt)} + \thetat \inner{\nui -\nuitt}{(\xtt -\xt[t-2])} \\
                                        & +\thetat \inner{\nuitt -\nuit}{(\xtt -\xt[t-2])}.
\end{align*}
Since $g^*_i$ has a strong convexity modulus 1 with respect to $U_{g_i^*}$, the definition of $\nuit$ prox-mapping in \eqref{ls:anla} implies a three-point inequality (see Lemma 3.5 in \cite{LanBook}):
\begin{align*}
    \inner{\nui - \nuit}{ \xt} & + \gistar(\nuit) - \gistar(\nui) + \inner{\nui -\nuit }{\xtilt-\xt}                                       \\
                               & \leq \taut U_{g_i^*}(\nui; \nuitt) - (\taut + 1) U_{g_i^*}(\nui; \nuit) - \taut U_{g_i^*}(\nuit; \nuitt).
\end{align*}
So, combining the above two relations,
taking the $\wt$ weighted sum of the resulting inequalities and using the conditions $\wtt=\wt\thetat$ and $\wt\taut \leq \wtt (\tautt  + 1)$, we obtain
\begin{align*}
    \sumt \wt & (\inner{\nui - \nuit}{ \xt} + \gistar(\nuit) - \gistar(\nui))                                          \\
    \leq      & - (\wt[N](\taut[N] + 1) U_{g_i^*}(\nui; \nuit[N]) -\wt[N]\inner{\nui - \nuit[N]}{ (\xt[N] -\xt[N-1])}) \\
              & -\tsum_{t=2}^{N} [\wt \taut U_{g_i^*}(\nuit; \nuitt) + \wtt \inner{\nuitt -\nuit}{(\xtt -\xt[t-2])} ]  \\
              & + \wt[1]\taut[1] U_{g_i^*}(\nui; \nuit[0]).
\end{align*}
A $\lami$-weighted sum of the above inequality leads to the desired $Q_\nu$ convergence bound given by
\begin{align*}
    \begin{split}
        \sumt \wt & Q_{\nu}(\zt;z)                                                                                                                \\
        \leq      & - (\wt[N](\taut[N] + 1) \sumi \lami U_{g_i^*}(\nui; \nuit[N]) -\wt[N]\inner{\sumi \lami (\nui - \nuit[N])}{\xt[N] -\xt[N-1]}) \\
                  & -\tsum_{t=2}^N [\wt \taut\sumi \lami U_{g_i^*}(\nuit; \nuitt) + \wtt \inner{\sumi \lami (\nuitt -\nuit)}{\xtt -\xt[t-2]} ]    \\
                  & + \wt[1]\taut[1] \sumi \lami U_{g_i^*}(\nui; \nuit[0]) .
    \end{split}
\end{align*}
A $Q_\pi$ bound can be derived similarly. Taken together, we get
\begin{align}\label{pr1:Qnu-Qpi}
    \begin{split}
         & \sumt \wt[t] [Q_{\nu}(\zt;z) + Q_\pi(\zt; z)]                                                                                                                                       \\
         & \ \leq - (\wt[N](\taut[N] + 1) [\sumi \lami U_{g_i^*}(\nui; \nuit) + \Ufstar(\pi; \pit)] -\wt[N]\inner{\pi - \pit[N] + \sumi \lami (\nui - \nuit[N])}{\xt[N] -\xt[N-1]})            \\
         & \quad \ -\tsum_{t=2}^N \{\wt \taut[t][\sumi \lami U_{g_i^*}(\nuit; \nuitt) + \Ufstar(\pit; \pitt)] + \wtt \inner{\pit[t-1] - \pit + \sumi \lami (\nuitt -\nuit)}{\xtt -\xt[t-2]} \} \\
         & \quad \ + \wt[1]\taut[1] [\sumi \lami U_{g_i^*}(\nui; \nuit[0]) + \Ufstar(\pi; \pit[0])] .
    \end{split}
\end{align}
Since $(\nut, \pit) = (\grad g(\xundert), \grad f(\xundert))\ \forall t$, applying Lemma \ref{lm:agg_strong_cvxity} with $\tilde{\lam}= [\lam; 1]$ and $\tilde{g} = [g; f]$ implies that
$$\sumi \lami U_{g_i^*}(\nuit; \nuitt) + \Ufstar(\pit; \pitt) \geq \tfrac{1}{2\Llam} \normsq{\pit - \pitt + \sumi \lami (\nuit - \nuitt)}.$$
Similarly, with $(\nu, \pi) \in [V, \Pi]$ (see \eqref{eq:joint-domain}), we get
$$\sumi \lami U_{g_i^*}(\nu; \nuit[N]) + \Ufstar(\pi; \pit[N]) \geq \tfrac{1}{2\Llam} \normsq{\pi - \pit[N] + \sumi \lami (\nui - \nuit[N])}.$$
Thus applying the Young's inequality to \eqref{pr1:Qnu-Qpi} leads to
\begin{align}\label{eq:pr-pf-Qnu}
    \begin{split}
        \sumt & \wt[t] [Q_{\nu}(\zt;z) + Q_\pi(\zt; z)]                                                                                                              \\
        \leq  & \tfrac{\wt[N]\Llam }{2(\taut[N] + 1)} \normsq{\xt[N] - \xt[N-1]} + \tsum_{t=1}^{N-1} \tfrac{\wt\thetat[t+1] \Llam }{2\taut[t+1]} \normsq{\xt - \xtt}
        + \wt[1]\taut[1] [\sumi \lami U_{g_i^*}(\nui; \nuit[0]) + \Ufstar(\pi; \pit[0])] .
    \end{split}
\end{align}

Now let us move onto $Q_x$ and $Q_\lam$. Fix $\lamt$. The $\xt$-prox mapping implies a three-point inequality (Lemma 3.5 in \cite{LanBook}):
$$\inner{\pit + \sumi \lamti\nuti}{\xt - x} + u(\xt) - u(x) +  \tfrac{\etat + \alpha}{2} \normsq{\xt - x} \leq \tfrac{\etat}{2} \normsq{\xtt - x} - \tfrac{\etat}{2} \normsq{\xt - \xtt}.$$
Fix $\xt$. The optimality of $\lamt$ implies that:
$$\inner{\lam - \lamt}{\nut \xt - g^*(\nut)} \leq 0.$$
So we have
$$Q_x(\zt; z) + Q_\lam(\zt; z) +  \tfrac{\etat + \alpha}{2} \normsq{\xt - x} \leq \tfrac{\etat}{2} \normsq{\xtt - x} - \tfrac{\etat}{2} \normsq{\xt - \xtt}.$$
Summing across iterations with weight $\wt$ and using \eqref{stp-con:eta-recur}, we get
\begin{equation}\label{eq:pr-pf-Qx}
    \sumwt[t] [Q_x(\zt; z) + Q_\lam(\zt; z)] +  \tfrac{\wt[N](\etat[N] + \alpha)}{2} \normsq{\xt - x} \leq \tfrac{\wt[1]\etat[1]}{2} \normsq{\xt[0] - x} - \sumt \tfrac{\wt\etat}{2} \normsq{\xt - \xtt}.
\end{equation}
Utilizing the stepsize assumption \eqref{stp-con:Q-cancel}, we can add it to \eqref{eq:pr-pf-Qnu} to obtain a convergence bound for the $Q$ function
\begin{equation}\label{eq:pr-pf-Q}
    \sumwt Q(\zt; z) + \tfrac{\wt[N](\etat[N] + \alpha)}{2} \normsq{\xt - x} \leq \tfrac{\wt[1]\etat[1]}{2} \normsq{\xt[0] - x} + \wt[1]\taut[1] [\sumi \lami U_{g_i^*}(\nui; \nuit[0]) + \Ufstar(\pi; \pit[0])].
\end{equation}
Moreover, the Jensen's inequality implies that
\begin{equation}\label{pr1:jensen}
    (\sumwt)\La(\xbart[N]; \lam,\nu, \pi) \leq \sumwt \La(\xt; \lam, \nu, \pi),\end{equation}
\begin{align*}
    \sumwt  & \La(x; \lamt, \nut, \pit) \jim{\leq} (\sumwt)[\inner{x}{\pibart[N]} - f^*(\pibart[N])+ u(x)] + \sumi(\sumwt\lamti)[ \inner{\nubarti[N]}{x}-g_i^*(\nubarti[N])] \\
    \jim{=} & (\sumwt)\La(x; \lambart[N], \nubarti[N][],\pibart[N] ).
\end{align*}
Thus, we get $(\sumwt )Q(\zbart[N]; z) \leq \sumwt Q(\zt; z)$, and the desired inequality in \eqref{eq:Q-conv} follows from \eqref{eq:pr-pf-Q}.
\endproof
\vgap
Next, the proof of Theorem \ref{thm:ANLA} is a direct application of Proposition \ref{pr:Q-convergence}. The analysis is complicated by the switch from a diminishing stepsize to a constant stepsize in \eqref{stp-choice:ANLA}.\\

\textbf{Proof of Theorem \ref{thm:ANLA}}
We apply Proposition \ref{pr:Q-convergence} to obtain the results. First, we verify that the requirements in \eqref{stp-con:Q-cancel} and \eqref{stp-con:Q-recur} are satisfied by the stepsize choice in \eqref{stp-choice:ANLA}. Since $\thetat=\taut/(\taut[t-1]+1)$ and $\etat = \Llamr/\taut[t+1]$, all requirements other than \eqref{stp-con:eta-recur} hold automatically.  Now
let $T = \ceil{2 \sqrt{\kar}} + 1$ denote the first iteration at which we switch to $\taut=\sqrt{\kar}.$ Since the other iterations and the case with $\alpha = 0$ are straightforward to check, we focus on iteration $T-1$ and $T$, and assume $\alpha >0$. For $t=T$, we have
\begin{align*}
    \wt[T]\etat[T] = \wt[T-1] \etat[T] \tfrac{\taut[T-1] + 1}{\taut[T]} = \wt[T-1] \tfrac{\Llamr}{\sqrt{\kar}} \tfrac{\ceil{2\sqrt{\kar} } + 1}{2\sqrt{\kar}} = \wt[T-1] \alpha \tfrac{2\sqrt{\kar} + 2}{2} \leq \wt[T-1] (\alpha\sqrt{\kar} + \alpha) \leq \wt[T-1] (\etat[T-1] + \alpha).
\end{align*}
For $t=T-1$, we have
\begin{align*}
    \wt[T-1]\etat[T-1] = \wt[T-2] \etat[T-1] \tfrac{\taut[T-2] + 1}{\taut[T-1]} = \wt[T-2] \tfrac{2\Llamr}{2\sqrt{\kar}} \tfrac{\ceil{2\sqrt{\kar}} -1}{\ceil{2\sqrt{\kar}} -2} \leq \wt[T-2] \tfrac{2 \Llamr}{\ceil{2\sqrt{\kar}} -2} = \wt[T-2] \etat[T-2] \leq \wt[T-2](\etat[T-2] + \alpha).
\end{align*}
Thus the requirements in Proposition \ref{pr:Q-convergence} are satisfied, and we have
\begin{equation}\label{thm-pf:ANLA-Q}
    Q(\zbart, z)  \leq \Llamr \normsq{\xt[0] - x}/(\sumwt)\ \forall z=(x; \lam, \nu, \pi)\in X \times \Lamr \times [V, \Pi].
\end{equation}

We provide a lower bound to $(\sumwt)$. It is useful to show
$$\wt \geq \max\{t, (1 + 1/\sqrt{\kar})^{t -5}\}.$$
Since $\wt= 1/(\prod_{t=2}^{N} \thetat)$, the fact $\wt \geq t\ \forall t \geq 1$ is straightforward. Regarding the second lower bound, let us  first consider $t \leq T = \ceil{2 \sqrt{\kar}} +1.$ The algebraic fact in Lemma \ref{lm:al_fact} implies that
$\wt[t+2] \geq t \geq (1 + 1/\sqrt{\kar})^{t-3} \forall t \leq 2 \sqrt{\kar}$, so
$\wt \geq (1 + 1/\sqrt{\kar})^{t-5} \forall t \leq T = \ceil{2 \sqrt{\kar}} + 1$. For $t\geq T+1$, the relation is also valid since $1/\thetat = (1 + 1/\sqrt{\kar}).$ Therefore we get $\wt \geq \max\{t, (1 + 1/\sqrt{\kar})^{t -5}\} \forall t \geq 1$. Using these two lower bounds, it is easy to derive
\begin{equation}\label{eq:wt_lower_bound}
    \sumwt \geq \max\{N(N+1)/2, \sqrt{\kar}[(1+ 1/\sqrt{\kar})^{N-4} -1]\}.\end{equation}
Substituting the preceding lower bound into \eqref{thm-pf:ANLA-Q} and applying Lemma \ref{lm:fun-feas-from-Q} lead us to the convergence results in \eqref{eq:thm-non-strong} and \eqref{eq:thm-strong}.

Now we deduce the convergence bound for $\norm{\xt - \xstar}$. Choosing $\hz = (\xstar; \lamstar, \grad g(\xbarN), \grad f(\xbart[N]))$ to be the reference point $z$ leads to
\begin{equation}\label{thm-pf:ANLA-Q-dist}
    Q(\zbart; \hz) \leq \Llamr \normsq{\xt[0] - \xstar}/(\sumwt).\end{equation}
Moreover, since $\xstar$ minimizes the optimal Lagrangian, we get  $\inner{\grad f(\xstar) + (\lamstar)^\top \grad g(\xstar) + u'(\xstar)}{\xbarN - \xstar} \geq 0$ for all $u'(\xstar) \in \partial u(\xstar).$ So we have
\begin{align}\label{thmpf:Q-strong}
    \begin{split}
        Q(\zbart; \hz) & \geq f(\xbarN) + (\lamstar)^\top g(\xbarN) + u(\xbarN) - [f(\xstar) + (\lamstar)^\top g(\xstar) + u(\xstar)]                                                  \\
                       & \quad - \inner{\grad f(\xstar) + (\lamstar)^\top \grad g(\xstar) + u'(\xstar)}{\xbarN - \xstar}                                                               \\
                       & = f(\xbarN) - f(\xstar) - \inner{\grad f(\xstar)}{\xbarN - \xstar} + \sumi \lamstar_i [g_i(\xbarN) - g_i(\xstar)- \inner{\grad g_i(\xstar)}{\xbarN - \xstar}] \\
                       & \quad + [u(\xbarN) - u(\xstar) - \inner{u'(\xstar)}{\xbarN - \xstar}]                                                                                         \\
                       & \geq \alpha \normsq{\xbarN - \xstar}/2,\end{split}
\end{align}
where the last inequality follows from the fact that the regularization function $u(x)$ is $\alpha$-strongly convex.
Combining the above two relations, we obtain the desired inequality in \eqref{eq:thm-strong-dist}:
$$\normsq{\xbarN - \xstar} \leq 2 \kar \normsq{\xt[0] - \jim{\xstar}}/(\sumwt)\leq \tfrac{2 \sqrt{\kar} \normsq{\xt[0] - \xstar}}{(1+ 1/\sqrt{\kar})^{N-4} - 1}.$$

\endproof

\subsection{The Adaptive Search for $\Llamr$}\label{subsec:ANLA-search}
A crucial limitation of the ACGD method is that the aggregate Lipschitz smoothness constants, $\Llamr[c]$ and $\Llamr[1]$, required for the stepsize calculation in \eqref{stp-choice:ANLA}, are often unavailable in practice. To address this, we propose incorporating an adaptive restart mechanism into Algorithm \ref{alg:ANLA} to dynamically estimate these constants.
Since our goal is to find an $(\ep, \ep/c)$-solution, the following penalty-scale parameter is useful for our development:
\begin{equation}\label{eq:lamplus}
    \jim{\rplus_c := \max\{\norm{\lamstar}, c\}},
\end{equation}
where $\norm{\lamstar}$ is the norm of the optimal dual multiplier. 

To achieve adaptivity, we require a test to determine whether our estimated aggregate smoothness constant $\Ltilagg$ is misspecified and adjust it accordingly. One natural approach is to run the ACGD for a fixed number of iterations and check if both the feasibility violation and the optimality gap are converging according to the theoretical bounds, e.g., those implied by Theorem \ref{thm:ANLA}. However, the number of iterations required for such a test would match the worst-case theoretical complexity bound. Since empirical observations show that actual algorithm runs on real-world problems typically require far fewer iterations than the theoretical upper bound, this approach could be prohibitively expensive in practice.

In this subsection, we propose a verifiable termination certificate to accurately track the algorithm's progress and develop tests based on it.
To be useful, \jim{such a certificate should be readily computable and provide meaningful bounds on both the optimality gap and the feasibility violation. In this respect, it plays a role analogous to that of the gradient norm in smooth optimization \cite{lan2023optimal,zhang2025linearly}. The following property formalizes these requirements.}

\begin{property}\label{propy:verifiable-certificate-requirements}
    A gap $\Delta$ is a verifiable certificate if it satisfies the following three requirements:
    \begin{enumerate}[label=\alph*)]
        \item \textit{Verifiability:} The gap $\Delta$ can be evaluated without knowledge of any unknown problem parameters.
        \item \textit{Convergence:}  The gap $\Delta$ converges to $0$ at the same rate that the ACGD method would achieve as if $\Ltilagg$ were a valid aggregate smoothness constant.
        \item \textit{Optimality \& Feasibility:} If $\Ltilagg$ is well specified, the gap $\Delta$ provides bounds of the same order on both the optimality gap and the feasibility violation.
    \end{enumerate}
\end{property}
We now describe how to use the gap $\Delta$ to construct a test for the ACGD method. We run the ACGD method as if  $\Ltilagg$ were a valid aggregate Lipschitz smoothness constant. By Property \ref{propy:verifiable-certificate-requirements}.a), we can calculate both $\Delta_t$ and the feasibility violation $\delta_t$ at each iteration.
Importantly, Property \ref{propy:verifiable-certificate-requirements}.b) implies that $\Delta_t$ converges efficiently regardless of whether the estimate $\Ltilagg$ is valid.
Thus, if $\delta_t$ remains of the same order as $\Delta_t$, both the optimality gap and the feasibility violation converge with the desired rate.
Otherwise, if $\delta_t$ does not decrease in proportion to $\Delta_t$, Property \ref{propy:verifiable-certificate-requirements}.c) implies that $\Ltilagg$ is misspecified, in which case we could restart the ACGD method with a larger estimate, e.g., $2\Ltilagg$.
Since the empirical $\Delta_t$ is often much smaller than the value implied by the worst-case theoretical upper bound, this adaptive test enables us to detect an erroneous $\Ltilagg$ much earlier, ensuring that the adaptive parameter-free method remains competitive with the ACGD method implemented with the true $\Llamr[2\jim{\rplus_c}]$ (see Section \ref{sec:numerical} for numerical results).

Towards this end, we propose the following FP-gap as a verifiable certificate.
\begin{definition}\label{def:linear-cone-certificate}
    Given an evaluation point $\xbar$ and some lower linear approximation $l_f$ to $f$ and the lower linear approximations $l_{g_1}, l_{g_2}, \ldots, l_{g_m}$ to $g_1, g_2, \ldots, g_m$, the FP-gap parameterized by $r$, denoted $\DeltaFP(r)$, is defined as:
    \begin{equation}\label{eq:FP-gap}
        \begin{split}
             & \DeltaFP(r) := F(\xbar) + r\norm{[g(\xbar)]_+} -\Funder                                                     \\
             & \text{ where } \Funder := \min_{x \in X} l_f(x) + u(x), \text{ s.t. } l_{g_i}(x) \leq 0, \forall i \in [m].
        \end{split}
    \end{equation}
\end{definition}
A few remarks are in order regarding this definition.
First, the lower bound $\Funder$ on $\Fstar$ is calculated as the minimum of \jim{$l_f+u$ over the relaxed feasible region} $\mathcal{G}:=\{x\in X: l_{g_i}(x) \leq 0\ \forall i \}$, which is a relaxation of the true feasible region $\{x\in X:g_i(x)\leq 0\ \forall i\}$. \jim{In particular, when $X$ is polyhedral, $\mathcal{G}$ is a polyhedron, which motivates the name \textbf{F}easible \textbf{P}olyhedron gap. Notice that the FP-gap reduces to the usual primal-dual gap \cite{LanBook} when there are no function constraints.}
Second, the parameter $r$ specifies the dual radius $\Lambda \subset B^m_+(0;r)$ associated with the constraints and is related to the estimate $\Ltilagg$. \jim{Third, to ensure that $\Funder$ is finite, we impose the following boundedness assumption on $X$ throughout Subsections~\ref{subsec:ANLA-search}, \ref{subsec:a-ANLA-conv-analysis}, and~\ref{sec:adaptive-ACGD-S}.

    \begin{assumption}[Bounded Domain]\label{ass:bounded_domain}
        The set $X$ is nonempty, closed, bounded, and convex, with finite diameter
        \[
            D_X := \sup_{x,y \in X} \norm{x-y} < \infty.
        \]
        The value of $D_X$ need not be known by the algorithm.
    \end{assumption}
}

Now let us see why such a FP-gap corresponds to a verifiable certificate. \jim{Since $\Funder$ can be calculated by evaluating a constrained prox-mapping of the form \eqref{eq:constrained-prox-def},} the value $\DeltaFP(r)$ can be calculated without the input of any unknown problem parameters, thus it satisfies Property \ref{propy:verifiable-certificate-requirements}.a).
The following lemma demonstrates that the FP-gap satisfies Property \ref{propy:verifiable-certificate-requirements}.c) if the parameter $r$ is sufficiently large (see \eqref{eq:a-ANLA-r-lower-bound}).

\begin{lemma}\label{lm:FP-gap-to-suboptimality-feasibility}
    For the evaluation point $\xbar$, the FP-gap yields the following upper bounds on feasibility violation and optimality gap:
    \begin{enumerate}
        \item[a)] $F(\xbar) - \Fstar\leq \DeltaFP(r)$ for any $r\geq 0$.
        \item[b)] If $r \geq 2\jim{\rplus_c}$, then
              \begin{equation}\label{eq:FP-gap-to-feasibility-violation}
                  \norm{[g(\xbar)]_+} \leq 2\DeltaFP(r)/r \jim{\leq \DeltaFP(r)/c}.
              \end{equation}
    \end{enumerate}
\end{lemma}
\begin{proof}
    Part a) follows directly from the definition of the FP-gap. For part b), if $r \geq 2\jim{\rplus_c}$, we have:
    \begin{align*}
        F(\xbar) + r\norm{[g(\xbar)]_+} & \geq F(\xbar) + \norm{\lamstar}\norm{[g(\xbar)]_+} + \half r  \norm{[g(\xbar)]_+}                                                                \\
                                        & \geq F(\xbar) + (\lamstar)^\top g(\xbar) + \half r \norm{[g(\xbar)]_+}                                                                           \\
                                        & \stackrel{(a)}{\geq} F(\xstar) + (\lamstar)^\top g(\xstar) \jim{+ \half r \norm{[g(\xbar)]_+}}\geq \Funder + \half r \norm{[g(\xbar)]_+}.
    \end{align*}
    Here, (a) uses the fact that the minimal value of the optimal Lagrangian is attained at $\xstar$. Rearranging the terms gives:
    \begin{align*}
        \norm{[g(\xbar)]_+} & \leq \frac{2(F(\xbar) + r\norm{[g(\xbar)]_+} - \Funder)}{ r}       \\
                            & \leq \frac{2\DeltaFP(r)}{r}\jim{\leq\frac{\DeltaFP(r)}{c}}.
    \end{align*}
\end{proof}

Next, to establish Property \ref{propy:verifiable-certificate-requirements}.c), we provide the exact construction of $l_f$ and $l_{g_i}$ based on the iterates generated by the ACGD method. With those choices, we show that the convergence of $\DeltaFP(r)$ is implied by the convergence of the $Q$-gap associated with the dual variables $\Lambda\subset B^m_+(0;r)$ in the ACGD method.

\begin{lemma}\label{lm:FP-gap-from-Q}
    Given some primal dual iterates $\{(\xt[j],\lam^j, \pi^j, \nu^j)\}_{j=1}^t$ and some nonnegative weights $\{\wt\}_{j=1}^N$, let $\xbar^t$ be the weighted average of these iterates, i.e., $\xbar^t = \tsum_{j=1}^t \wt[j] \xt[j] /\tsum_{j=1}^t \wt[j]$. If the $Q$-gap associated with them is bounded by some $\Delta(r)$ for $\hat{\pi}^t = \grad f(\xbar^t)$, $\hat{\nu}^t = \grad g(\xbar^t)$,
    \begin{equation*}
        \max_{x\in X, \lambda \in B^m_+(0;r)}  \frac{1}{\tsum_{j=1}^t \wt[j]}\tsum_{j=1}^t \wt[j]\left(\La(\xt[j]; \lam, \hat{\nu}^t, \hat{\pi}^t) - \La(x;\lam^j,\nu^j,\pi^j)\right)\leq \Delta(r).
    \end{equation*}
    Then the following lower linear approximation functions $l_f$ and $l_g$ corresponds to a FP-gap certificate associated with $\xbar^t$ with $\DeltaFP(r) \leq \Delta(r)$,
    \begin{equation}\label{eq:FP-certificate-from-Q}
        \begin{split}
             & l_f(x) := \left(\tsum_{j=1}^t \wt[j] [\inner{\pi^j}{x} - f^*(\pi^j)] \right)/(\tsum_{j=1}^t \wt[j]),                                 \\
             & l_{g_i}(x):=\left(\tsum_{j=1}^t \wt[j] \lam^j_i [\inner{\nu_i^j}{x} - g^*_i(\nu_i^j)]  \right) / \left(\tsum_{j=1}^t \wt[j] \right).
        \end{split}
    \end{equation}

\end{lemma}
\begin{proof}
    With $\hat{\pi}^t = \grad f(\xbart[t])$, $\hat{\nu}^t = \grad g(\xbart[t])$, and $\lamhat = r \frac{[g(\xbar^t)]_+}{\norm{[g(\xbar^t)]_+}}$
    \begin{align*}
        F(\xbart[t]) & +  r\norm{[g(\xbar^t)]_+} = F(\xbar^t) + \lamhat^\top g(\xbar^t)                                                                                        \\
                     & \leq   \max_{\lambda \in B^m_+(0;r)}  F(\xbar^t) + \lambda^\top g(\xbar^t)                                                                              \\
                     & \stackrel{(a)}{=} \max_{\lambda \in B^m_+(0;r)} \La(\xbart[t]; \lambda, \hat{\nu}^t, \hat{\pi}^t )                                                      \\
                     & \stackrel{(b)}{\leq}  \max_{\lambda \in B^m_+(0;r)} \frac{1}{\tsum_{j=1}^t \wt[j]} \tsum_{j=1}^t \wt[j] \La(\xt[j]; \lambda, \hat{\nu}^t, \hat{\pi}^t).\end{align*}
    Here (a) follows from the conjugate duality relationship (see Sections 4.1 and 4.2 of \cite{Beck2017First}), and (b) follows from the convexity of the nested Lagrangian function with respect to the primal variable $x$.
    For the lower bound, since $\wt[j]$ and $\lam_i^j$ are all non-negative, we have
    \begin{align*}
        \underline{F} & = \begin{cases} \min_{x \in X } & l_f(x) + u(x)\\ \text{s.t.} & l_{g_i}(x) \leq 0\quad \forall i \in [m]\end{cases} = \begin{cases} \min_{x \in X } & \frac{1}{\tsum_{j=1}^t \wt[j]} \tsum_{j=1}^t \wt[j] [\inner{\pi^j}{x} - f^*(\pi^j)] + u(x)\\ \text{s.t.} & \frac{1}{\tsum_{j=1}^t \wt[j]  }\tsum_{j=1}^t \wt[j] \lam^j_i [\inner{\nu_i^j}{x} - g^*_i(\nu_i^j)] \leq 0\quad \forall i \in [m]\end{cases} \\
                      & = \begin{cases} \min_{x \in X } & \frac{1}{\tsum_{j=1}^t \wt[j]} \tsum_{j=1}^t \wt[j] [\inner{\pi^j}{x} - f^*(\pi^j)] + u(x)\\ \text{s.t.} & \frac{1}{\tsum_{j=1}^t \wt[j]  }\tsum_{j=1}^t \wt[j] \lam^j_i [\inner{\nu_i^j}{x} - g^*_i(\nu_i^j)] \leq 0\quad \forall i \in [m]\end{cases}                                                                                                                     \\
                      & \stackrel{(c)}{\geq} \min_{x \in X} \frac{1}{\tsum_{j=1}^t \wt[j]} \left[ \tsum_{j=1}^t \wt[j] \left( \inner{\pi^j}{x} - f^*(\pi^j)+u(x)\right)  + \tsum_{j=1}^t \wt[j] \sum_{i=1}^m  \lam^j_i [\inner{\nu_i^j}{x} - g^*_i(\nu_i^j)] \right]                                                                                                                                                                                                                                   \\
                      & {\geq} \min_{x \in X} \frac{1}{\tsum_{j=1}^t \wt[j]} \tsum_{j=1}^t \wt[j] \left[\inner{\pi^j}{x} - f^*(\pi^j) +u(x)+ \sum_{i=1}^m \lam^j_i [\inner{\nu_i^j}{x} - g^*_i(\nu_i^j)]\right]                                                                                                                                                                                                                                                                                        \\
                      & \geq \min_{x\in X} \frac{1}{\tsum_{j=1}^t \wt[j]} \tsum_{j=1}^t \wt[j] \La(x; \lam^j, \nu^j, \pi^j),
    \end{align*}
    \jim{Here (c) is the Lagrangian relaxation obtained with the all-ones multiplier vector $\mathbf{1}\in\mathbb{R}_+^m$.}
    Taken together, we get $\DeltaFP(r) = F(\xbart[t]) + r \norm{[g(\xbar^t)]_+} - \underline{F} \leq \Delta(r)$.
\end{proof}

\jim{Two remarks about the FP-gap certificate are in order. First, updating the affine minorants $l_{g_i}$ requires access to the dual multipliers associated with the constrained prox-mapping in Line~\ref{aACGD-line:x-projection}. As shown in Section~\ref{sec:ANLA-S}, these multipliers arise naturally from the sliding updates. Moreover, the aggregation in~\eqref{eq:FP-certificate-from-Q} can be updated recursively, thereby avoiding the need to store the full history of iterates.} Second, the value of the parameter $r$ plays a crucial role in our development. Lemma \ref{lm:FP-gap-to-suboptimality-feasibility} requires a large $r$ to establish the connection to feasibility violation, while Lemma \ref{lm:FP-gap-from-Q} and later Proposition \ref{pr:a-Q-convergence} shows that $\DeltaFP(r)$ convergence could be guaranteed only for $r\leq \tilde{r}$ for some $\tilde{r}$ related to the estimate $\Ltilagg$. Specifically, at the $t^\text{th}$ iteration, the upper bound $\tilde{r}_t$, derived from the condition $\Ltilagg \ge \tilde{L}_{f,t} + r \tilde{L}_{g,t}$, is given by
\begin{equation}\label{eq:a-ANLA-r-lower-bound}
    \jim{\tilde{r}_t:=\begin{cases}
            \dfrac{\Ltilagg-\tilde{L}_{f,t}}{\tilde{L}_{g,t}}, & \tilde{L}_{g,t}>0,                                          \\
            +\infty,                                           & \tilde{L}_{g,t}=0 \text{ and } \Ltilagg\geq\tilde{L}_{f,t}, \\
            0,                                                 & \tilde{L}_{g,t}=0 \text{ and } \Ltilagg<\tilde{L}_{f,t}.
        \end{cases}}
\end{equation}
where $\tilde{L}_{f,t}$ and $\tilde{L}_{g,t}$ denote certain empirical estimates of the Lipschitz smoothness constans based on iterates generated up to iteration $t$.
\jim{When $\tilde r_t=+\infty$, expressions involving $\tilde r_t$ in the certificate and the restart test are interpreted by taking the limit as $r\to+\infty$.}




Using this $\tilde{r}_t$ as the parameter $r$ in the FP-gap (cf. Definition \ref{def:linear-cone-certificate}) leads to the adaptive ACGD method shown in Algorithm~\ref{alg:a-ACGD}. Starting from an initial guess $\Ltilagg$, Lines~3--6 coincide with the original ACGD update. Lines~7--11 carry out the adaptive restart: Line~7 updates the running estimates of $\Ltil_f$ and $\Ltil_g$; Line~8 computes the ratio $\tilde r$ using the estimated $\Ltil_f$ and $\Ltil_g$ (see \eqref{eq:a-ANLA-r-lower-bound}); Line~9 evaluates the $FP$-gap; and Lines~10--11 trigger a restart whenever the feasibility violation fails to match the $FP$-gap convergence.

In particular, the theorem below shows the complexity of the adaptive ACGD method.
\begin{theorem}\label{thm:a-ANLA}
    \jim{For a prescribed scaling constant $c\geq1$,} consider a smooth constrained optimization problem of form \eqref{eq:prob} \jim{with a bounded domain $X$ satisfying Assumption \ref{ass:bounded_domain}}. Provided with initial smoothness estimates $\Ltil_f=\alpha$ and $\Ltil_g=0$, an aggregate estimate $\alpha\leq\tilde{L}_{\text{agg},0} \leq \jim{\Lagg^c:=L_f+4\rplus_c L_g}$ with $\tilde{L}_{\text{agg},0}>0$, and the strong convexity modulus $\alpha$, the adaptive ACGD method in Algorithm \ref{alg:a-ACGD} finds an \jim{$(\ep, \ep/c)$-optimal solution} within the following $T_\epsilon$ iterations for any $\ep > 0$.
    \begin{enumerate}
        \item[a)] In the non-strongly convex case with $\alpha=0$, the total number of adaptive ACGD iterations is bounded by
              \begin{equation}\label{eq:a-ANLA-complexity-non-strong}
                  T_\epsilon =\frac{7\sqrt{\jim{\Lagg^c}}\DX}{\sqrt{\ep}} + \ceil{\log_2(\jim{\Lagg^c}/\tilde{L}_{\text{agg},0})}.\end{equation}
        \item[b)] In the strongly convex case with $\alpha>0$,  the total number of adaptive ACGD iterations is bounded by
              \begin{align}\label{eq:a-ANLA-complexity-strong}
                  T_\epsilon =\left(\sqrt{\frac{2\jim{\Lagg^c}}{\alpha}} + \ceil{\log_2(\frac{\jim{\Lagg^c}}{\tilde{L}_{\text{agg},0}})}\right)\left(\jim{\log(1+\frac{2\jim{\Lagg^c} \DX^2}{\ep})} + 1\right) +4\ceil{\log_2(\frac{\jim{\Lagg^c}}{\tilde{L}_{\text{agg},0}})}
              \end{align}
    \end{enumerate}
\end{theorem}

A few comments are in order regarding the result. First, compared to the ACGD method implemented with the true $L(\Lambda_c)$, the above convergence bound achieves the same order of oracle complexity. The adaptive ACGD method incurs some additional additive $\log$ terms, the constant dependence changes from $\norm{\xt[0]-\xstar}$ to $\DX$, and the aggregate smoothness constant changes from $L(\Lambda_c)$ to the more conservative $\Lagg^c$.
Second, in practice, it might be advisable to start with some small and erroneous estimate of $\tilde{L}_{\text{agg},0}$ at the beginning so that the adaptive method could revise it up to an appropriate one. Because the algorithm convergence depends only on the empirical smoothness constant, rather than the global constants, the adaptive method could potentially converge with a larger stepsize than the optimal choice in theory, which translates into a performance improvement, as observed in Section \ref{sec:numerical}.
\jim{Third, evaluating the FP-gap at each iteration requires solving an additional LP in the non-strongly convex case or an additional QP in the strongly convex case. This extra solve can be costly in large-scale settings, particularly when both $n$ and $m$ are large. In such settings, one can instead use the adaptive ACGD-S method in Algorithm~\ref{alg:adaptive-ACGD-S}, which replaces the FP-gap with the PD-gap introduced in Subsection~\ref{sec:adaptive-ACGD-S}. Evaluating the PD-gap requires only a simple linear or quadratic minimization over $X$.} 


\begin{minipage}{.95\textwidth}
    \begin{algorithm}[H]
        \caption{Adaptive ACGD Method}
        \label{alg:a-ACGD}
        \begin{algorithmic}[1]

            \Require $\xt[-1]= \xundert[0] = \xt[0]\in X$, an aggregate Lipschitz smoothness estimate $\Ltilagg\geq\alpha$, smoothness estimates $\Ltil_f\geq\alpha$ and $\Ltil_g\geq0$, \jim{the scaling constant $c\geq1$}, and the strong convexity modulus $\alpha$.
            \State Set $\pit[0] = \grad f(\xt[0]) $, $\nut[0] = \grad g(\xt[0])$
            \parState{Compute stepsize parameter $\wt$, $\etat$, $\thetat$, and $\taut$ according to \eqref{stp-choice:ANLA} using $\Ltilagg$.}
            \For{$t =  1,2,3 ... $}
            \parState{Set $\xundert \leftarrow (\taut \xundert[t-1] + \xtilt)/(1 + \taut)$ where  $\xtilt  = \xtt + \thetat (\xtt - \xt[t-2])$.}
            \parState{\label{aACGD-line:x-projection}Set $\pit \leftarrow \grad f(\xundert)$ and $\nut \leftarrow \grad g(\xundert)$.}
            \parState{Solve $\xt \leftarrow \argmin_{x \in X} \{\inner{\pit}{x} +u(x) + \etat \normsq{x - \xtt}/2\ s.t.\ \nut (x - \xundert) + g(\xundert) \leq 0\}$ and \jim{record the dual multipliers as  $\lambda^t$.}} 
            \parState{\label{aACDG-line:update-Ltil} Compute the average solution $\xbar^t:= \tsum_{j=1}^t \wt[j] \xt[j] / (\tsum_{j=1}^t \wt[j])$ and update the Lipschitz smoothness constant estimates:
                $$\Ltil_f\leftarrow \max\{\Ltil_f,  \tfrac{1}{2}\tfrac{\normsq{\grad f(\xunder^t)-\grad f(\xunder^{t-1})}}{f(\xunder^{t-1})-f(\xunder^t)-\inner{\grad f(\xunder^{t})}{\xunder^{t-1}-\xunder^{t}}}, \jim{\tfrac{1}{2}\tfrac{\normsq{\grad f(\xbar^t)-\grad f(\xunder^{t})}}{f(\xunder^t)-f(\xbar^t)-\inner{\grad f(\xbar^t)}{\xunder^{t}-\xbar^{t}}}}\},$$
                \begin{align*}\Ltil_g\leftarrow \max\left\{\Ltil_g, \tfrac{1}{2} \norm{[\tfrac{\jim{\normsq{\grad g_1(\xunder^t)-\grad g_1(\xunder^{t-1})}}}{g_1(\xunder^{t-1})-g_1(\xunder^t)-\inner{\grad g_1(\xunder^{t})}{\xunder^{t-1}-\xunder^{t}}},...,\tfrac{\jim{\normsq{\grad g_m(\xunder^t)-\grad g_m(\xunder^{t-1})}}}{g_m(\xunder^{t-1})-g_m(\xunder^t)-\inner{\grad g_m(\xunder^{t})}{\xunder^{t-1}-\xunder^{t}}}]},\right. \\
                    \left.\jim{\tfrac{1}{2} \norm{[\tfrac{\normsq{\grad g_1(\xbar^t)-\grad g_1(\xunder^{t})}}{g_1(\xunder^t)-g_1(\xbar^t)-\inner{\grad g_1(\xbar^t)}{\xunder^{t}-\xbar^{t}}},...,\tfrac{\normsq{\grad g_m(\xbar^t)-\grad g_m(\xunder^{t})}}{g_m(\xunder^t)-g_m(\xbar^t)-\inner{\grad g_m(\xbar^t)}{\xunder^{t}-\xbar^{t}}}]}}\right\}.
                \end{align*}
            }
            \parState{Calculate $\jim{\tilde{r}:=\begin{cases}(\Ltilagg-\Ltil_f)/\Ltil_g, & \Ltil_g>0,\\ +\infty, & \Ltil_g=0\text{ and }\Ltilagg\geq\Ltil_f,\\ 0, & \Ltil_g=0\text{ and }\Ltilagg<\Ltil_f.\end{cases}}$}
            \parState{Calculate FP-gap $\DeltaFP(\tilde{r}):=F(\xbar^t) + \tilde{r} \norm{[g(\xbar^t)]_+} - \Funder$ where $\Funder$ is minimal value of the lower linear approximation functions in \eqref{eq:FP-gap} attained by the lower linear approximation functions in \eqref{eq:FP-certificate-from-Q}.}
            \If{\label{aACDG-line:restart-trigger}$\tilde{r} <\jim{2c}$ or $\norm{[g(\xbar^t)]_+} > 2\DeltaFP(\tilde{r})/\tilde{r}$}
            \State Restart the Adaptive ACGD method with inputs $\xt[0]$, $2\Ltilagg$, $\Ltil_f$, $\Ltil_g$, \jim{$c$}, and $\alpha$.
            \EndIf
            \EndFor
        \end{algorithmic}
    \end{algorithm}
    \vgap
\end{minipage}


\subsection{Convergence Analysis of the Adaptive ACGD Method}\label{subsec:a-ANLA-conv-analysis}

We present the detailed proofs of convergence results for Algorithm \ref{alg:a-ACGD}. Being the counterpart to Proposition \ref{pr:Q-convergence}, the following proposition shows the convergence of the $Q$-gap function. However, the stepsize choices depend on estimated Lipschitz smoothness constants, $\Ltil_f$ and $\Ltil_g$, rather than the true ones, $L_f$ and $L_g$.
\begin{proposition}\label{pr:a-Q-convergence}
    Consider the primal dual iterates ${\zt:=(\xt; \lamt, \nut, \pit)}$ generated during any given run of Algorithm \ref{alg:a-ACGD}. Let $\Ltil_{f,t}$, $\Ltil_{g,t}$ and $\rtil_t$ denote the values of $\Ltil_f$, $\Ltil_g$ and $\tilde{r}$ at the $t$-th iteration, respectively. For any iteration $t$ strictly preceding the first iteration, if any, at which the restart condition in Line~\ref{aACDG-line:restart-trigger} is satisfied, we have the following convergence bound associated with the averaged solution $\xbar^t$ and the corresponding FP-certificate (see \eqref{eq:FP-certificate-from-Q} and Definition \ref{def:linear-cone-certificate}):
    \jim{
        \begin{equation}\label{eq:a-Q-conv}
            \begin{split}
                \sum_{j=1}^t \wt[j] \DeltaFP(\rtil_t) \leq \max_{x \in X}\frac{\etat[1]\wt[1]}{2} \normsq{\xt[0] - x}
            \end{split}
        \end{equation}
    }
\end{proposition}
\begin{proof} Since Lemma \ref{lm:FP-gap-from-Q} shows that the $FP$-gap is a lower bound to the $Q$-gap parameterized by same $\rtil_t$, it suffices to show the following relation holds for $\hat{\nu}^t=\grad g(\xbar^t)$, $\hat{\pi}^t=\grad f(\xbar^t)$ and any $x \in X$.
    \begin{align*}
        \max_{\lam \in B^m_+(0;\rtil_t)} &
        \tsum_{j=1}^t \wt[j]\left(\La(\xt[j]; \lam, \hat{\nu}^t, \hat{\pi}^t) - \La(x;\lam^j,\nu^j,\pi^j)\right)+ \frac{\wt[t](\etat[t] + \alpha)}{2} \normsq{\xt[t]-x} \\
        \leq                             & \max_{x \in X}\frac{\etat[1]\wt[1]}{2} \normsq{\xt[0] - x}.
    \end{align*}
    However, this follows similar to the proof of Proposition \ref{pr:Q-convergence}. The key new ingredient is to utilize the estimated Lipschitz smoothness constant $\Ltil_{f,t}$ and $\Ltil_{g,t}$ as opposed to the theoretical aggregate smoothness constant $L(\Lambda_r)$. Specifically, we need to
    use the stepsize choice $\etat[j-1]\tau_j \geq \thetat[j] (\Ltil_{f,t}+\rtil_t \Ltil_{g,t})$ to provide the following inequality for any $\lam \in B^m_+(0;\rtil_t)$, and any $j \leq t$
    \begin{align*}
        \tau_j [\tsum_{i=1}^m \lam_i U_{\gistar[i]}(\nu_i^j; \nu_i^{j-1})+U_{f^*}(\pi^j;\pi^{j-1})] & +\inner{\theta_j (\xt[j-1]-\xt[j-2])}{\pi^{j-1} -\pi^{j}+\tsum_{i=1}^m \lam_i(\nu^{j-1}_i -\nu_i^{j})} \\
        \geq                                                                                        & -\frac{\theta_j\eta_{j-1}}{2} \normsq{\xt[j-1]-\xt[j-2]}.
    \end{align*}
    Let us first focus on the terms associated with $\pi$. Since $\pi^t = \grad f(\xunder^t)$ and $\pi^{t-1} = \grad f(\xunder^{t-1})$, the algebraic identity between Bregman distance functions generated by $f$ and its Fenchel conjugate $f^*$ satisfies
    \begin{align*}
        U_{f^*}(\pi^j;\pi^{j-1})=\jim{U_f(\xunder^{j-1}; \xunder^j)} \stackrel{(a)}{\geq} \frac{1}{2 \Ltil_{f,j}} \normsq{\grad f(\xunder^j)-\grad f(\xunder^{j-1})} \\
        \stackrel{(b)}{\geq} \frac{1}{2 \Ltil_{f,t}} \normsq{\grad f(\xunder^j)-\grad f(\xunder^{j-1})}.
    \end{align*}
    Here (a) follows from the definition of $\Ltil_{f,j}$ in Line \ref{aACDG-line:update-Ltil} of Algorithm \ref{alg:a-ACGD}, and (b) follows $\Ltil_{f,t}$ being monotonically non-decreasing. Thus it follows from the Young's inequality that
    $$\tau_j U_{f^*}(\pi^j;\pi^{j-1})+ \inner{\theta_j (\xt[j-1]-\xt[j-2])}{\pi^{j-1} -\pi^{j}}\geq -\frac{\theta^2_j \Ltil_{f,t}}{2 \jim{\tau_j}} \normsq{\xt[j-1]-\xt[j-2]}.$$
    Similarly, if we let $\Ltil_{g,j,i}:= \frac{1}{2} \normsq{\grad g_i(\xunder^j)-\grad g_i(\xunder^{j-1})}/[g_i(\xunder^{j-1})-g_i(\xunder^j)-\inner{\grad g_i(\xunder^j)}{\xunder^{j-1}-\xunder^j}]$, we have
    $$\tau_j U_{g_i^*}(\nu_i^j;\nu_i^{j-1})+ \inner{\theta_j (\xt[j-1]-\xt[j-2])}{\nu_i^{j-1} -\nu_i^{j}}\geq -\frac{\theta^2_j \Ltil_{g,t,i}}{2\tau_j} \normsq{\xt[j-1]-\xt[j-2]}.$$
    Adding up the above inequalities with weight $\lambda_i$, we have
    \begin{align*}
        \tau_j \tsum_{i=1}^m \lambda_i U_{g_i^*}(\nu_i^j;\nu_i^{j-1}) & + \inner{\theta_j (\xt[j-1]-\xt[j-2])}{\tsum_{i=1}^m \lambda_i(\nu_i^{j-1} -\nu_i^{j})}\geq -\frac{\theta^2_j \tsum_{i=1}^m \lambda_i \Ltil_{g,t,i}}{2\tau_j} \normsq{\xt[j-1]-\xt[j-2]} \\
                                                                      & \geq -\frac{\theta^2_j \norm{\lam}\Ltil_{g,t}}{2\tau_j} \normsq{\xt[j-1]-\xt[j-2]}\geq -\frac{\theta^2_j \rtil_t\Ltil_{g,t}}{2\tau_j} \normsq{\xt[j-1]-\xt[j-2]}.
    \end{align*}
    Here the second last inequality follows from the Holder's inequality, and the last inequality follows from the fact that $\norm{\lam} \leq \rtil_t$. Summing them up, we have
    \begin{align*}
        \tau_j [\tsum_{i=1}^m \lam_i U_{\gistar[i]}(\nu_i^j; \nu_i^{j-1})+U_{f^*}(\pi^j;\pi^{j-1})] & +\inner{\theta_j (\xt[j-1]-\xt[j-2])}{\pi^{j-1} -\pi^{j}+\tsum_{i=1}^m \lam_i\jim{(\nu_i^{j-1}-\nu^j_i)}} \\
        \geq                                                                                        & -\frac{\theta_j}{2} \frac{\theta_j(\Ltil_{f,t}+\rtil_t \Ltil_{g,t})}{\tau_j} \normsq{\xt[j-1]-\xt[j-2]}.
    \end{align*}
    Since the stepsize is chosen according to \eqref{stp-choice:ANLA} with $\Ltilagg$ in place of $L(\Lamr[r])$, we have
    $$\eta_{j-1}\tau_{j}\geq \Ltilagg \theta_j=\theta_j(\Ltil_{f,t}+\rtil_t \Ltil_{g,t}),$$
    where the last equality follows from the definition of $\rtil_t$.
    Then the desired inequality follows from the fact that $\etat[j-1] \geq \thetat[j] (\Ltil_{f,t}+\rtil_t \Ltil_{g,t}) /\tau_j$.

    Additionally, we also have to show that the stepsize choice $\etat(\taut +1) \geq  \Ltilagg =(\Ltil_{f,t}+\rtil_t \Ltil_{g,t})$ ensures
    \begin{align*}
        (\tau_t + 1) [\tsum_{i=1}^m \lam_i U_{\gistar[i]}(\jim{\hat{\nu}_i^t}; \nu_i^{t})+U_{f^*}(\hat{\pi}^t;\pi^{t})] & +\inner{ (\xt-\xt[t-1])}{\hat{\pi}^t -\pi^{t}+\tsum_{i=1}^m \lam_i(\hat{\nu}^t_i -\nu_i^{t})} \\
        \geq                                                                                                            & -\frac{\eta_t}{2} \normsq{\xt-\xt[t-1]}.
    \end{align*}
    \jim{Indeed, the conjugate-Bregman identity and the second empirical ratios in
        Line~\ref{aACDG-line:update-Ltil} give
        \begin{align*}
            U_{f^*}(\hat{\pi}^t;\pi^t)
             & =U_f(\xunder^t;\xbar^t)
            \geq \frac{\normsq{\hat{\pi}^t-\pi^t}}{2\Ltil_{f,t}}, \\
            U_{g_i^*}(\hat{\nu}_i^t;\nu_i^t)
             & =U_{g_i}(\xunder^t;\xbar^t)
            \geq \frac{\normsq{\hat{\nu}_i^t-\nu_i^t}}{2\Ltil_{g,t,i}},
            \qquad i\in[m],
        \end{align*}
        where $\Ltil_{g,t,i}$ denotes the corresponding component of the empirical
        constraint-smoothness estimate. Applying Young's inequality to these bounds
        and using $\norm{\lambda}\leq\rtil_t$ gives
        \begin{align*}
             & (\tau_t+1)\left[\tsum_{i=1}^m\lam_iU_{g_i^*}(\hat{\nu}_i^t;\nu_i^t)+U_{f^*}(\hat{\pi}^t;\pi^t)\right]
            +\inner{\xt-\xt[t-1]}{\hat{\pi}^t-\pi^t+\tsum_{i=1}^m\lam_i(\hat{\nu}_i^t-\nu_i^t)}                      \\
             & \quad\geq-\frac{\Ltil_{f,t}+\rtil_t\Ltil_{g,t}}{2(\tau_t+1)}\normsq{\xt-\xt[t-1]}
            \geq-\frac{\eta_t}{2}\normsq{\xt-\xt[t-1]},
        \end{align*}
        where the last inequality follows because
        $\eta_t(\tau_t+1)\geq
            \Ltil_{f,t}+\rtil_t\Ltil_{g,t}$.}

\end{proof}

\begin{proof}\textbf{of Theorem \ref{thm:a-ANLA}}
    We first consider the non-strongly convex case with $\alpha=0$. The following facts are useful for the proof. \jim{First, the co-coercivity inequality for convex smooth functions and the definitions in Line~\ref{aACDG-line:update-Ltil} imply
        $\Ltil_{f,t}\leq L_f$ and $\Ltil_{g,t}\leq L_g$. Suppose that the guessed aggregate smoothness constant satisfies $\Ltilagg\geq\Lagg^c=L_f+4\rplus_cL_g$. If $\Ltil_{g,t}>0$, then \eqref{eq:a-ANLA-r-lower-bound} gives
        \[
            \rtil_t
            =\frac{\Ltilagg-\Ltil_{f,t}}{\Ltil_{g,t}}
            \geq\frac{4\rplus_cL_g}{\Ltil_{g,t}}
            \geq4\rplus_c\geq2\rplus_c.
        \]
        If $\Ltil_{g,t}=0$, then $\Ltilagg\geq L_f\geq\Ltil_{f,t}$, so \eqref{eq:a-ANLA-r-lower-bound} gives $\rtil_t=+\infty$. Thus $\rtil_t\geq2\rplus_c\geq2c$, and Lemma~\ref{lm:FP-gap-to-suboptimality-feasibility}.b) ensures that neither condition in Line~\ref{aACDG-line:restart-trigger} can trigger a further restart.}
    Second, with a fixed $\Ltilagg$ during a run of the adaptive ACGD method, Proposition \ref{pr:a-Q-convergence} implies, at every iteration strictly before the next restart (or throughout the run if no restart occurs), that
    $$\DeltaFP(\rtil_t) \leq \max_{x\in X}\frac{2\Ltilagg}{2(t)(t+1)} \normsq{\xt[0]-x}.$$
    \jim{Whenever no restart is triggered, $\rtil_t\geq2c$ and Line~\ref{aACDG-line:restart-trigger} guarantees
        $\norm{[g(\xbar^t)]_+}\leq2\DeltaFP(\rtil_t)/\rtil_t\leq\DeltaFP(\rtil_t)/c$. Together with Lemma~\ref{lm:FP-gap-to-suboptimality-feasibility}.a), this shows that $\DeltaFP(\rtil_t)\leq\ep$ yields an $(\ep,\ep/c)$-optimal solution.} Thus, the total number of iterations required for either a restart or an \jim{$(\ep,\ep/c)$-optimal solution} is upper bounded by $\sqrt{\frac{2\Ltilagg}{\ep}}\DX + 1.$
    Since the guess $\Ltilagg$ is doubled every time a restart is triggered, taken together, the total number of iterations required for generating an \jim{$(\ep,\ep/c)$-optimal solution} is upper bounded by
    $$T_\epsilon =\frac{7\sqrt{\jim{\Lagg^c}}\DX}{\sqrt{\ep}} + \ceil{\log_2(\jim{\Lagg^c}/\tilde{L}_{\text{agg},0})}.$$
    Similarly, in the strongly convex case with $\alpha>0$, the total number of iterations required for either a restart or an \jim{$(\ep,\ep/c)$-optimal solution} is upper bounded by
    $$\left(\sqrt{\frac{2\Ltilagg}{\alpha}} + 1\right)\left(\jim{\log(1+\frac{2\Ltilagg \DX^2}{\ep})} + 1\right) +4.$$
    Thus the total number of iterations required for either a restart or an \jim{$(\ep,\ep/c)$-optimal solution} is upper bounded by
    $$T_\epsilon =\left(\sqrt{\frac{2\jim{\Lagg^c}}{\alpha}} + \ceil{\log_2(\frac{\jim{\Lagg^c}}{\tilde{L}_{\text{agg},0}})}\right)\left(\jim{\log(1+\frac{2\jim{\Lagg^c} \DX^2}{\ep})} + 1\right) +4\ceil{\log_2(\frac{\jim{\Lagg^c}}{\tilde{L}_{\text{agg},0}})}.$$

\end{proof}

\color{black}

\section{Lower Oracle Complexity Bound}\label{sec:lower}

In this section, we present the lower oracle complexity bounds, that is, the minimum number of queries to the FO oracle required to find an $(\ep; \ep/c)$-optimal solution. These results illustrate the optimality of the ACGD method. We assume, for the sake of simplicity, that  $u(x)= \alpha\normsq{x}/2$ (see \eqref{eq:prob}) and that $X$ is the Euclidean ball  $B(0; \jim{R_X})$ centered at the origin with radius $\jim{R_X}\in \bar{\R}$.

Similar to Nesterov's lower complexity computation model in \cite{nesterov2003introductory}, we consider the class of all first-order methods, $\calF$, verifying a linear-span update requirement. Given a (finite) memory of reachable points $\Mt[t-1]$ after the $t-1$\ts{th} query to FO, the updated memory $\Mt$ after evaluating the FO oracle at some $y \in \Span(\Mt[t-1])$ needs to satisfy
\begin{equation}\label{req:fo-model}
  \Mt \subset \{x + \etat \grad f(y) + \sumi {\taut[i,t]} \grad g_i(y): x,y \in \Span(\Mt[t-1]), \etat, \taut[i,t] \in \R \}.
\end{equation}
The freedom to choose arbitrary elements from the linear span allows $\calF$ to cover many first-order algorithms. For example, the ACGD method is a special case of $\calF$ because the generated points in the ACGD method, $\Mt =\{\xt[0], \xt[1], \ldots, \xt[t]\}$, satisfy the requirement \eqref{req:fo-model}. Specifically, when the memory $\Mt[t-1] = \{\xt[0], \xt[1], \ldots, \xt[t-1]\}$, the evaluation point $\xundert$ in Line 3 is inside $\Span(\Mt[t-1])$. Moreover, the $\xt$-update in Line 5 of Algorithm \ref{alg:ANLA} can be expressed as
\begin{align*}
                                    & \xt \leftarrow \argmin_{x \in X} \{\inner{\grad f(\xundert)}{x} + u(x) + \etat \normsq{x - \xtt}/2\ s.t.\ \grad g(\xundert) (x - \xundert) + g(\xundert) \leq 0\} \\
  \stackrel{(a)}{\Leftrightarrow}\  & \xt  \leftarrow \argmin_{x \in X} \{\inner{\grad f(\xundert) + \sumi \lamti \grad g_i(\xundert)}{x}  + \alpha\norm{x}^2/2 + \etat \normsq{x - \xtt}/2  \}         \\
  \Leftrightarrow\                  & \jim{\xt \leftarrow \argmin_{x \in X} \normsq{x - \tfrac{1}{\etat + \alpha}[\etat \xt[t-1] - (\grad f(\xundert) + \sumi \lamti \grad g_i(\xundert))]}}            \\
  \stackrel{(b)}{\Leftrightarrow}\  & \xt =  \tfrac{1}{\gamma (\etat + \alpha)}[ \xt[t-1] - \frac{1}{\etat}(\grad f(\xundert) + \sumi \lamti \grad g_i(\xundert))]\  \text{for some $\gamma>0$,}
\end{align*}
where the multiplier $\lamt$ in (a) is the optimal dual solution \jim{associated with the constrained prox-mapping in the first line}, and (b) follows from $X=B(0; \jim{R_X})$. Therefore, $\xt$ is a member of the right-hand side of \eqref{req:fo-model}, $\Mt := \Mt[t-1] \cup \{\xt\}$ satisfies the update requirement \eqref{req:fo-model}, and the ACGD method is a member of $\calF$. In fact, similar arguments can be used to show that $\calF$ covers both the primal methods \cite{lin2018level,nesterov2003introductory},  and the primal-dual methods \cite{boob2022stochastic,aybat2021primal,xu2020first} in the literature.

Since the dependence on parameters of the smooth objective function is well-established, e.g. \cite{nemirovsky1983problem,nesterov2003introductory}, we will investigate the dependence of the lower complexity bounds on the parameters of the constraint function. Toward that end, we consider an affine $f$ in the objective, i.e. $\Lf = 0$, and study the dependence of the lower complexity bound on the aggregate smoothness constant (see \eqref{eq:agg-L}).

\subsection{The Strongly Convex Problem}
\jim{
  First, we study the strongly convex problem with $\alpha>0$. Since linear convergence is expected, different optimality criteria have little impact on the lower complexity bound, that is, they only have different constants inside the ``$\log$''. So, without loss of generality, we choose to focus on the convergence of $\jim{\normsq{\xt - \xstar}}$ in the next theorem.

  \begin{theorem}\label{thm:aggregate-strong-lower}
    Let $l\geq1$ and $\bar L\geq2\alpha>0$ be given. For every first-order
    method in $\calF$ with an arbitrary starting point $\xt[0]$, there exists an
    infinite-dimensional hard problem of the form \eqref{eq:prob} with
    $\norm{\lamstar}=l$, affine $f$, and aggregate smoothness constant
    $L(\Lambda_1)\leq\bar L$ \eqref{eq:agg-L} such that, for every
    $0<\ep\leq\normsq{\xt[0]-\xstar}$, at least
    \[
      \Omega\left(\sqrt{\frac{\bar L}{\alpha}}
      \log\frac{\normsq{\xt[0]-\xstar}}{\ep}\right)
    \]
    queries to the first-order oracle are required to find an iterate $\xt[N]$
    satisfying $\normsq{\xt[N]-\xstar}\leq\ep$.
  \end{theorem}

  \begin{proof}
    Consider the following hard problem parameterized by $\beta>0$. Let
    \[
      \rho:=\frac{\alpha}{\alpha+4\beta},\qquad
      \Delta:=\frac{1-\sqrt{\rho}}{1+\sqrt{\rho}},
      \qquad
      \xstar:=(\Delta,\Delta^2,\ldots,\Delta^i,\ldots)\in\R^\infty,
    \]
    and consider
    \begin{equation}\label{ls:aggregate-strong-lower}
      \begin{split}
        \min_{x\in\R^\infty}\quad
         & -\beta x_1+\frac{\alpha}{2}\normsq{x} \\
        \mathrm{s.t.}\quad
         & g(x):=h(x)-h(\xstar)\leq0,\qquad
        h(x):=\frac{\beta}{2l}
        \left[x_1^2+\sum_{i=1}^\infty(x_i-x_{i+1})^2\right].
      \end{split}
    \end{equation}
    The objective in \eqref{ls:aggregate-strong-lower} has strong convexity
    modulus $\alpha$, its smooth component $f$ is affine, and $g$ is
    $(4\beta/l)$-smooth. Moreover, $\lamstar=l$ and the above $\xstar$ satisfy
    the KKT conditions because
    \begin{align*}
      g(\xstar) & =0,         \\
      \left(\beta
      \begin{bmatrix}
        2  & -1                            \\[-1mm]
        -1 & 2  & -1                       \\
           & -1 & 2      & -1     & \ddots \\
           &    & \ddots & \ddots & \ddots
      \end{bmatrix}
      +\alpha I\right)\xstar
                & =\beta e_1.
    \end{align*}
    Thus $\norm{\lamstar}=l$. Since $\Lambda_1=\{0\}\cup[l,l+1]$ and $f$
    is affine, the aggregate smoothness constant satisfies
    \[
      L(\Lambda_1)
      \leq\frac{4(l+1)\beta}{l}.
    \]

    Without loss of generality, take $\xt[0]=0$ and $\Mt[0]=\{0\}$. The
    zero-chain structure in
    \eqref{ls:aggregate-strong-lower} and the linear-span requirement imply
    inductively that every iterate $\xt$ has nonzero entries only in its first
    $t$ coordinates. Therefore
    \begin{equation}\label{pf:aggregate-strong-tail}
      \normsq{\xt-\xstar}
      \geq\sum_{i=t+1}^\infty\Delta^{2i}
      =\Delta^{2t}\sum_{i=1}^\infty\Delta^{2i}
      =\Delta^{2t}\normsq{\xt[0]-\xstar}.
    \end{equation}

    Now choose $\beta:=l\bar L/[4(l+1)]$. Since $l\geq1$ and
    $\bar L\geq2\alpha$, we have
    $L(\Lambda_1)\leq4(l+1)\beta/l=\bar L$ and $0<\rho\leq1/2$, while
    \[
      \frac{1}{\sqrt{\rho}}=\sqrt{1+\frac{4\beta}{\alpha}}
      \geq\sqrt{\frac{\bar L}{2\alpha}}.
    \]
    Taking logarithms in \eqref{pf:aggregate-strong-tail} and using
    $\log\Delta\geq-5\sqrt{\rho}$ therefore gives
    \[
      t\geq\frac{1}{10\sqrt{\rho}}
      \log\frac{\normsq{\xt[0]-\xstar}}{\ep}
      =\Omega\left(\sqrt{\frac{\bar L}{\alpha}}
      \log\frac{\normsq{\xt[0]-\xstar}}{\ep}\right),
    \]
    which proves the result.
  \end{proof}

  \subsection{The Non-Strongly Convex Problem}
  Now we move on to consider the non-strongly convex problem with $\alpha =0$. The next theorem states the lower oracle complexity bound to find an $(\ep; \ep/c)$-optimal solution.
  \begin{theorem}\label{thm:aggregate-non-strong-lower}
    Let problem parameters $\bar L_c>0$, $\Rt[0]\geq1$, $l\geq1$, $c\geq1$, and
    $0<\ep\leq \bar L_c\Rt[0]^2/100$ be given. For a sufficiently large problem
    dimension,
    $n>2\ceil{\Rt[0]\sqrt{\bar L_c/\ep}}$, there exists a hard problem of the
    form \eqref{eq:prob} with $\norm{\lamstar}=l$,
    $\norm{\xt[0]-\xstar}\leq\Rt[0]$, and affine $f$, whose aggregate smoothness
    constant in \eqref{eq:agg-L} satisfies $L(\Lambda_c)\leq\bar L_c$, such that every first-order
    method in $\calF$ requires at least
    $\Omega(\sqrt{\bar L_c}\Rt[0]/\sqrt{\ep})$ queries to the first-order oracle
    to find an $(\ep;\ep/c)$-optimal solution.
  \end{theorem}

  \begin{proof}
    Consider the following problem parameterized by $\beta>0$, $\gamma>0$, and
    $k\in\mathbb{N}_+$:
    \begin{align}\label{ls:aggregate-lower-eg}
      \begin{split}
        \min_{x\in\R^{2k+1}}\quad & -2l\gamma\beta x_1                                           \\
        \mathrm{s.t.}\quad
                                  & g_1(x):=\beta\left[x_1^2+\sumi[2k](x_i-x_{i+1})^2+x_{2k+1}^2
        -\frac{2k+1}{2k+2}\gamma^2\right]\leq0,                                                  \\
                                  & g_2(x):=4l\beta\left[-2\gamma x_1+x_1^2
        +\sumi[2k](x_i-x_{i+1})^2+x_{2k+1}^2
        +\frac{2k+1}{2k+2}\gamma^2\right]\leq0.
      \end{split}
    \end{align}
    Without loss of generality, take $\xt[0]=0$ and $\Mt[0]=\{0\}$. Let
    $\calK_i:=\{x\in\R^{2k+1}:x_j=0\ \forall j>i\}$. Since the first-order
    information at any point in $\calK_i$ is contained in $\calK_{i+1}$, the
    linear-span requirement implies inductively that $\Mt\subset\calK_t$.
    Therefore, after $k$ oracle queries,
    \begin{align}\label{pf:aggregate-feas-lowB}
      \norm{[g(x^k)]_+}
       & \geq\min_{x\in\calK_k}g_2(x)\notag      \\
       & \geq4l\beta\gamma^2\left(-\frac{k}{k+1}
      +\frac{2k+1}{2k+2}\right)
      =\frac{4l\beta\gamma^2}{2k+2}.
    \end{align}

    It is straightforward to verify from the KKT conditions that
    \[
      \lamstar=(l,0),\qquad
      \xstar_i=\gamma\left(1-\frac{i}{2k+2}\right),\quad i\in[2k+1].
    \]
    Hence $\norm{\lamstar}=l$ and
    $\norm{\xt[0]-\xstar}\leq\gamma\sqrt{k+1}$. Moreover, $g_1$ and
    $g_2/(4l)$ are $8\beta$-smooth. For any
    $\lambda\in\Lambda_c\setminus\{0\}$, we have
    $\lambda=\lamstar+\delta$ for some $\delta\in\R_+^2$ with
    $\norm{\delta}\leq c$, and thus
    \[
      L_\lambda
      \leq8\beta(\lambda_1+4l\lambda_2)
      \leq8\beta(l+c+4lc)\leq48lc\beta.
    \]

    Now set
    \[
      K:=\Rt[0]\sqrt{\frac{\bar L_c}{\ep}},\qquad
      k:=\floor{K/5}-1,\qquad
      \beta:=\frac{\bar L_c}{48lc},\qquad
      \gamma:=\frac{\Rt[0]}{\sqrt{k+1}}.
    \]
    The assumed accuracy range ensures $K\geq10$, so $k\geq1$ and
    $k=\Omega(K)$; the dimension assumption also gives $2k+1<n$. These choices
    yield $\norm{\xt[0]-\xstar}\leq\Rt[0]$ and
    $L(\Lambda_c)\leq\bar L_c$. It follows from
    \eqref{pf:aggregate-feas-lowB} and $k+1\leq K/5$ that
    \[
      \norm{[g(x^k)]_+}
      \geq\frac{\bar L_c\Rt[0]^2}{24c(k+1)^2}
      \geq\frac{25\ep}{24c}
      >\frac{\ep}{c}.
    \]
    Therefore every method in $\calF$ requires at least
    $\Omega(\sqrt{\bar L_c}\Rt[0]/\sqrt{\ep})$ oracle queries to find an
    $(\ep;\ep/c)$-optimal solution.
  \end{proof}

  These lower bounds establish that the ACGD method in Algorithm~\ref{alg:ANLA} attains the optimal oracle complexity. A key advantage of ACGD is its ability to exploit the aggregate Lipschitz smoothness constant $L(\Lambda_c)$ of the Lagrangian, defined in \eqref{eq:agg-L}, rather than relying on the conservative upper bound $(\norm{\lamstar}+c)L_g$ based on the aggregate constraint smoothness constant $L_g$ in \eqref{eq:Lg-def}. More precisely, we have the following relationship on the hard instance
  \[
    \frac{(\norm{\lamstar}+c)L_g}{L(\Lambda_c)}
    = \bigO\left(\frac{\norm{\lamstar}}{c}+1\right).
  \]
  For the adaptive ACGD method in Algorithm~\ref{alg:a-ACGD}, the bound $\Lagg^c \leq 4(\norm{\lamstar}+c)L_g$ implies that its oracle complexity may be suboptimal by a multiplicative factor of $\sqrt{\norm{\lamstar}/c+1}$. This potential loss stems from the lack of an empirical estimate of $L(\Lambda_c)$ when $\lamstar$ is unknown. When $c\geq\norm{\lamstar}$, this factor is bounded by a constant, and adaptive ACGD is therefore order-optimal.
}

\section{The ACGD-S method}\label{sec:ANLA-S}
We extend the ACGD method to the ACGD with sliding (ACGD-S) method to handle the large-scale problem where both the problem dimension $n$ and the number of constraints $m$ are large. This section follows the same structure as Section \ref{sec:ANLA}. We first discuss the computation bottleneck in the large-scale setting. Then Subsection \ref{subsec:S-method} introduces the ACGD-S method and presents the convergence results, and Subsection \ref{subsec:S-proof} contains the detailed proofs to the convergence results.

Despite its optimal oracle complexities, the ACGD method may lack in computation efficiency for large-scale problems. The  bottleneck of Algorithm \ref{alg:ANLA} lies in  Line 5:
\begin{align}\label{eq:ANLA-QP}\begin{split}
    \xt \leftarrow \argmin_{x \in X} & \inner{\pit}{x} +u(x) + \etat \normsq{x - \xtt}/2                  \\
    s.t.                             & \ \nuit (x - \xundert) + g_i(\xundert) \leq 0\  \forall i \in [m].
  \end{split}\end{align}
\jim{This is a constrained prox-mapping of the form \eqref{eq:constrained-prox-def}. When $u$ is quadratic, it reduces to a large-scale quadratic program (QP) if $X$ is polyhedral, such as a box, and to a large-scale quadratically constrained quadratic program (QCQP) if $X$ is a Euclidean ball.}

In this section, we address the bottleneck by replacing the \jim{large-scale constrained prox-mapping} with a sequence of basic matrix-vector operations, each requiring at most $\bigO(mn)$ FLOPs.
Overall, the ACGD-S method uses a number of matrix-vector operations comparable to that required to solve a single linearly constrained problem (i.e., a problem with affine $g$), while retaining the same optimal oracle complexity as the ACGD method.
Towards that end, we need to assume the projection onto $X$ is easy, i.e., the following operation can be computed in $\bigO(n)$ FLOPs for any $\pi, \xbar \in \R^n$ and $\etat \geq 0.$
\begin{equation}\label{eq:x-projection}
  x^t \leftarrow \targmin_{x \in X} \inner{\pi}{x} + u(x) + \etat\normsq{x - \xbar}/2.\end{equation}
For instance, if $u(x) = \alpha\normsq{x}/2$, the computation simplifies to component-wise thresholding if $X$ is a box, and to vector scaling if $X$ is a Euclidean ball. If $X$ is more challenging, we can model the complicated part using the function constraints. 

\begin{algorithm}[htb]
  \caption{The ACGD-S Method}
  \label{alg:ANLA-S}
  \begin{algorithmic}[1]
    \Require $\xt[-1]= \xundert[0] = \yst[0][1] =\xt[0]\in X$,  stepsizes $\{\thetat\}$, $\{\etat\}$, $\{\taut\}$, and  weights $\{\wt\}$.
    \State Set $\pit[0] = \grad f(\xt[0]) $, $\nut[0] = \grad g(\xt[0])$, $\lamst[-1][1] = \lamst[0][1] = 0 $.
    \For{$t =  1,2,3 ... N$}
    \parState{Set $\xundert \leftarrow (\taut \xundert[t-1] + \xtilt)/(1 + \taut)$ where  $\xtilt  = \xtt + \thetat (\xtt - \xt[t-2])$.}
    \parState{Set $\pit \leftarrow \grad f(\xundert)$ and $\nut \leftarrow \grad g(\xundert)$.}
    \parState{Calculate inner loop iteration limit $\{\St\}$ stepsizes $\{\betast\}$ and $\{\gamst\}$,  and weights $\{\delst\}$.}
    \For{$s=1,2,..., \St$}
    \parState{Set $\htilst = \begin{cases} (\nut)^\top\lamst[0] + \rhost[1] (\nut[t-1])^\top (\lamst[0] - \lamst[-1]) & \text{if } s =1 \\
              (\nut)^\top\lamst[s-1] + \rhost (\nut[t])^\top (\lamst[s-1] - \lamst[s-2]) & \text{ o.w.}\end{cases}$}
    \parState{Set $\yst \leftarrow \argmin_{y \in X}\inner{\htilst + \pit}{y} + u(y) + \etat\normsq{y - \xtt}/2 +  \betast \normsq{y - \yst[s-1]}/2. $}
    \parState{Set $\lamst \leftarrow \argmax_{\lam \in \R^m_+}\inner{\lam}{\nut[t] (\yst - \xundert) + g(\xundert)} - \gamst \normsq{\lam - \lamst[s-1]}/2.$}
    \EndFor
    \parState{Set $\lamst[0][t+1] = \lamst[\St][t]$, $\lamst[-1][t+1]= \lamst[\St-1][t]$, $\yst[0][t+1] = \yst[\St][t]$}
    \parState{Set $\xt = \sums \delst \yst /(\sums \delst)$ and $\lamtilt = \sums \delst \lamst /(\sums \delst)$.}
    \EndFor
    \State \Return $\xbart[N]:= \tsum_{t=1}^N \wt \xt / (\sumt \wt). $
  \end{algorithmic}
\end{algorithm}

\subsection{The ACGD-S Method and its Convergence Results}\label{subsec:S-method}
The ACGD-S method, listed in Algorithm \ref{alg:ANLA-S}, consists of two loops. For clarity, we will call an outer iteration a phase, and an inner iteration an iteration for the rest of this section. In phase $t$, the outer loop updates happen in Lines 3, 4, and 12; they are identical to the ACGD method except that the exact solution $\xt$ to the \jim{constrained prox-mapping} in \eqref{eq:ANLA-QP} is replaced by some average of inner iterates.
The other steps, Lines 5-11,  constitute the sliding subroutine. Its goal is to solve the Lagrangian reformulation to \eqref{eq:ANLA-QP}, or the $( \xt, \lamt)$ saddle point problem in \eqref{ls:anla}, inexactly:
\begin{equation}\label{eq:slide-saddle}
  (\xt, \lamt) \leftarrow \argmin_{y \in X} \argmax_{\lam \in \R^m_+} \inner{\pit}{y} + \inner{\lam}{\jim{\nut (y - \xundert) +g(\xundert)}} + u(y) + \etat \normsq{y - \xtt}/2.\end{equation}
To avoid confusion, we use the dummy variable $y$ to emphasize it being used only in the inner loop. Specifically, Line 5 calculates the stepsize parameters and iteration number $\St$. Lines 7-9 carry out primal-dual type updates for $\St$ iterations. Line 7 computes a momentum extrapolation term $\htilst$ as a proxy for $(\nut)^\top  \lamst$. In Line 8, with the variable $(\nut)^\top \lambda$ being fixed to $\htilst$,   $\yst$ is generated by minimizing the variable $y$ in \eqref{eq:slide-saddle} subject to a prox-function $\betast \normsq{y - \yst[s-1]}$. Then in Line 9,  with the variable $y$ fixed to $\yst$, $\lamst$ is generated by  maximizing \eqref{eq:slide-saddle} subject to a prox-function $\gamst \normsq{\lambda - \lamst[s-1]}$. After that, Line 11 prepares the initialization points for the inner loop in the next phase.

We highlight three features that are essential for achieving the desired computation efficiency. First, rather than being pre-specified, the inner loop stepsize parameters and iteration limit $\St$ are calculated in an online fashion in Line 5. This allows the method to adjust dynamically to the varying difficulty of the saddle point problems \eqref{eq:slide-saddle} from different phases. Second, the last operator $\nut[t-1]$, rather than $\nut$, is used for calculating the momentum extrapolation term at the first iteration $s=1$ in Line 7. This is characteristic of the sequential dual type algorithms \cite{lan2022optimal,zhang2020optimal,zhang2019efficient} for solving the trilinear saddle point problem in \eqref{eq:nLa}. Third, two primal iterates, $\xt$ and $\yst[\St]$, are stored after each inner loop to kick-start the next one. This is common to sliding-type algorithms \cite{lan2016gradient,lan2021graph,lan2022optimal}.

Now we suggest certain stepsize choices to obtain concrete convergence rates for Algorithm \ref{alg:ANLA-S}.
The non-strongly convex case and the strongly convex case are presented in separate theorems.
\begin{theorem}\label{thm:ANLA-S-non}
  Consider a non-strongly \jim{convex} problem of form \eqref{eq:prob} with $\alpha =0$. Let the aggregate smoothness constant $\Llamr[c]$ and the reference multiplier set $\Lamr[c]$ be defined in \eqref{eq:agg-L} and \eqref{eq:ref_lambs} respectively. Suppose Algorithm \ref{alg:ANLA-S} is run with the following stepsizes. The outer-loop stepsizes are
  \begin{equation}\label{stp:S-non-outer}
    \taut = \tfrac{t-1}{2}, \etat = \tfrac{\Llamr[c]}{\taut[t+1]}, \wt = t, \thetat[t+1] = \jim{\wt / \wt[t+1]}\  \forall t \geq 1.
  \end{equation}
  With $\jim{\MtU := \begin{cases}\norm{\nut}, & \norm{\nut}>0,\\ 1, & \norm{\nut}=0,\end{cases}}$ and some $R >0 $, the inner loop parameters in phase $t$ are calculated according to
  \begin{equation}\label{stp:S-non-inner}
    \begin{split}
       & \Delta  = R/\Llamr[c],  \St = \ceil{\MtU \Delta t}, \Mtilt = \tfrac{\St}{\Delta t},                                \\
       & \rhost = \begin{cases} \Mtilt/\Mtilt[t-1] & \jim{\text{if } s=1 \text{ and } t\geq2}         \\
              1                  & \jim{\text{if } s\geq2 \text{ or } (s,t)=(1,1)}, \\
                  \end{cases} \  \betast = \Mtilt R,\ \gamst = \Mtilt/R,\ \delst = 1\ \forall s \geq 1,\end{split}
  \end{equation}

  Then $\MtU \leq \Mbar\ \forall t$, where $\Mbar$ is an upper bound of $\jim{\max\{1,\norm{\grad g(x)}\}}$ for all $x$ in some bounded neighborhood around $\xstar$. Moreover, with  $\rlamr[c] := \norm{\lamstar} + c$ denoting the radius of the reference multiplier set $\Lamr[c]$, the ergodic average solution $\xbarN$ satisfies
  \begin{equation}\label{eq:S-non-strong-conv}
    \max\{F(\xbarN) - \Fstar, c \norm{[g(\xbarN)]_+}\} \leq \tfrac{\jim{\Llamr[c]}}{N(N+1)} (3\normsq{\xt[0] - \xstar} + \frac{\rlamr[c]^2}{R^2}).
  \end{equation}
  Thus to recover an $(\ep; \ep/c)$-optimal solution, we need at most $N_\ep$ gradient-oracle evaluations and $C_\ep$ matrix-vector multiplications, where
  \begin{align}
    \begin{split}\label{eq:S-non-strong-oracle-complexity-R}
      N_\ep & = \sqrt{\tfrac{\Llamr[c]}{\ep}}\sqrt{3\normsq{\xt[0] - \xstar}+  \tfrac{\rlamr[c]^2}{R^2}} +1,                           \\
      C_\ep & = O\{ N_\ep + \jim{\frac{\Mbar R}{\Llamr[c]}} + \frac{\Mbar}{\ep}(\normsq{\xt[0] - \xstar}R + \tfrac{\rlamr[c]^2}{R})\}.
    \end{split}
  \end{align}
\end{theorem}
\vspace{-3mm}
\begin{corollary}\label{cor:S-non}
  Under the setting of Theorem \ref{thm:ANLA-S-non}, if $R = \tfrac{\rlamr[c]}{\norm{\xt[0] - \xstar}}$, the total numbers of operations required by the ACGD-S method to find an $(\ep; \ep/c)$-optimal solution are bounded by:
  \begin{itemize}
    \item $N_\ep = \bigO(\sqrt{\tfrac{ \Llamr[c]}{\ep}}\norm{\xt[0] - \xstar})$ FO-oracle evaluations.
    \item $C_\ep = \bigO(\sqrt{\tfrac{\Llamr[c]}{\ep}}\norm{\xt[0] - \xstar} + \jim{\tfrac{\rlamr[c]\Mbar}{\Llamr[c]\norm{\xt[0]-\xstar}}} + \tfrac{ \rlamr[c] \Mbar\norm{\xt[0] - \xstar}}{\ep})$ matrix-vector multiplications.
  \end{itemize}
\end{corollary}
\vgap

Three remarks are in order regarding the above results. First, the oracle complexity of the ACGD-S method matches that of the ACGD method, while its computation complexity, measured by the number of matrix-vector multiplications, matches the lower bound for solving a single linearly constrained problem \cite{ouyang2021lower}. Second, the stepsize choices in \eqref{stp:S-non-outer} and \eqref{stp:S-non-inner} only require an upper bound to $\Llamr[c]$. Here the ratio parameter $R$ (chosen as $ \frac{\rlamr[c]}{\norm{\xt[0]-\xstar}}$ in Corollary \ref{cor:S-non})  is used to trade-off the cost associated with the primal and dual diameter, as is common in primal-dual type algorithms. An optimal choice of it would lead to the optimal constant dependence, however a misspecified one would still lead to an $\bigO(1/\sqrt{\ep})$ oracle complexity and an $\bigO(1/\ep)$ computation complexity. Third, the iteration limit function $\St$ in \eqref{stp:S-non-outer}  adapts to the varying difficulty of the saddle point sub-problem \eqref{eq:slide-saddle} from different phases. Specifically, $\St$ scales in proportion both to  $\norm{\nut}$, which characterizes the hardness of the saddle point sub-problem in  \eqref{eq:slide-saddle}, and to $t$, which captures the degree of accuracy required by the outer loop.
\vgap
\begin{theorem}\label{thm:ANLA-S-str}
  Consider a strongly convex problem of form \eqref{eq:prob} with $\alpha >0$. Let the aggregate smoothness constant $\Llamr[c]$ and the reference multiplier set $\Lamr[c]$ be defined in \eqref{eq:agg-L} and \eqref{eq:ref_lambs} respectively, and let $\kar[c]:=\Llamr[c]/\alpha$ be the condition number. Suppose Algorithm \ref{alg:ANLA-S} is run with the following stepsizes. The outer-loop stepsizes are
  \begin{equation}\label{stp:S-str-outer}
    \taut = \min\{\tfrac{t-1}{2}, \sqrt{2\kar}\}, \etat = \tfrac{\Llamr[c]}{\taut[t+1]}, \text{ and } \thetat = \tfrac{\taut}{\taut[t-1] + 1} \forall t \geq 1. \quad  \wt = \begin{cases} \wt[t-1] /\thetat & \text{if } t  \geq 2, \\
              1                 & \text{if } t = 1.\end{cases}
  \end{equation}
  For the inner loops, given some diameter ratio parameter $R> 0$, let $\jim{\MtU := \begin{cases}\norm{\nut}, & \norm{\nut}>0,\\ 1, & \norm{\nut}=0,\end{cases}}$ and $\Delta = \tfrac{R^2}{\alpha\Llamr[c]}$. At the beginning, the iteration limit is set to $\St[1] = \min \{S\in \mathbb{N}_+: \tsum_{s=1}^{S} s \geq \wt[1] \MtU[1]^2 \Delta \}$, and the stepsizes for all $s \in [\St[1]]$ are set to
  \begin{equation}\label{stp:S-str-init}
    \delst[s][1] = \tfrac{s}{\Gamt[1]},\ \betast[s][1] = \tfrac{\alpha}{4} (s-1),\ \gamst[s][1] = \tfrac{4}{\alpha}\tfrac{\MtU[1]^2 \Gamt[1]}{\delst[s][1]},\text{ and } \Wt[2] = \tfrac{\delst[\St[1]][1]}{ \MtU[1]},
  \end{equation}
  where $\Gamt[1]$ is the non-negative root to $\sums[\St[1]] s = \Gamt[1]^2 (\wt[1] \MtU[1]^2 \Delta).$ Then for phase $t \geq 2$, the iteration limit $\St$ and the parameter $\Gamt \geq 0$ are specified to satisfy
  \begin{equation}\label{stp:St-strong}
    \St = \min\{S\in\mathbb{N}_+: \sums[\St] \Wt \MtU + (s-1) \geq \wt \MtU^2 \Delta\},\ \sums \Gamt (\Wt \MtU) + (s-1) = \Gamt^2 \wt \MtU^2 \Delta,
  \end{equation}
  and the stepsizes are chosen according to
  \begin{equation}\label{stp:S-str-inner}
    \delst = \Wt \MtU + \tfrac{1}{\Gamt} (s-1),\ \gamst[s] = \tfrac{4}{\alpha}\tfrac{\MtU^2 \Gamt}{\delst[s]},\ \betast = \begin{cases}\tfrac{\alpha}{4}(\Wt\Gamt\MtU) & \text{ if } s=1, \\ \tfrac{\alpha}{4}[\Wt\Gamt\MtU + (s-2)] &\text{otherwise,}\end{cases}\ \text{ and } \Wt[t+1] = \tfrac{\delst[\St]}{ \MtU}.
  \end{equation}
  Then we have $\MtU \leq \Mbar\ \forall t$, where $\Mbar$ is an upper bound for $\jim{\max\{1,\norm{\grad g(x)}\}}$ for $x$ in some bounded neighborhood around $\xstar$, and the ergodic average solution $\xbarN$ satisfies
  \begin{equation}\label{eq:S-str-conv}
    \begin{split}
       & \max\{F(\xbarN) - \Fstar, c \norm{[g(\xbarN)]_+}\} \leq \tfrac{\Llamr[c]}{\calWt} (\tfrac{2\rlamr[c]^2}{R^2} + \normsq{\xt[0] - \xstar}), \\
       & \normsq{\xbarN - \xstar} \leq \jim{\tfrac{2\Llamr[c]}{\alpha\calWt}} (\tfrac{2\rlamr[c]^2}{R^2} + \normsq{\xt[0] - \xstar}),
    \end{split}
  \end{equation}
  where $\rlamr[c]=\norm{\lamstar} + c$ and the denominator satisfies $\calWt \geq \max\{N(N+1)/2, \sqrt{2\kar[c]}[(1+ 1/\sqrt{2\kar[c]})^{N-4} -1]\}.$
\end{theorem}

\vgap
\begin{corollary}\label{cor:S-str-opt}
  Under the setting of Theorem \ref{thm:ANLA-S-str}, \jim{let
    $A_R:=\normsq{\xt[0]-\xstar}+2\rlamr[c]^2/R^2$ and suppose that
    \[
      0<\ep\leq\frac{\sqrt{\Llamr[c]\alpha}\,A_R}{\sqrt{2}}.
    \]}
  The numbers of FO-oracle evaluations $N_\ep$ and of matrix-vector multiplications $C_\ep$ required by the ACGD-S method to find an $(\ep; \ep/c)$-optimal solution are bounded by:
  \begin{align*}
    N_\ep & = \jim{\bigO\left\{
      \sqrt{\kar[c]}\log\left(
      \frac{\sqrt{\Llamr[c]\alpha}\,A_R}{\ep}\right)+1
    \right\}},                  \\
    C_\ep & = \jim{\bigO\left\{
      N_\ep+\frac{\Mbar}{\sqrt{\alpha\ep}}
      \left(R\norm{\xt[0]-\xstar}+\rlamr[c]\right)
      \right\}}.
  \end{align*}
  Moreover, \jim{if $\xt[0]\neq\xstar$ and we choose
    $R = \rlamr[c]/\norm{\xt[0]-\xstar}$, then
    \begin{align*}
      N_\ep
       & =\bigO\left\{\sqrt{\kar[c]}
      \log\left(\frac{3\sqrt{\Llamr[c]\alpha}\,
      \normsq{\xt[0]-\xstar}}{\ep}\right)+1\right\}, \\
      C_\ep
       & =\bigO\left\{
      \frac{\Mbar\rlamr[c]}{\sqrt{\alpha\ep}}+N_\ep
      \right\}.
    \end{align*}
  }
\end{corollary}
\vgap
\begin{corollary}\label{cor:S-str-close}
  Under the assumptions of Theorem \ref{thm:ANLA-S-str}, \jim{suppose
    $\xt[0]\neq\xstar$. To find an $\ep$-close solution satisfying
    $\normsq{\xbart[N] - \xstar} \leq \ep$, set $c=1$ and
    $R = \rlamr[1]/\norm{\xt[0]-\xstar}$.} The required complexities are then bounded by:
  \begin{align*}
    N_\ep & \jim{\leq \left\lceil
      (\sqrt{2\kar[1]}+1)
      \log\left(1+\frac{3\sqrt{2\kar[1]}
        \normsq{\xt[0]-\xstar}}{\ep}\right)
    \right\rceil+4},                 \\
    C_\ep & \jim{=\bigO\left\{N_\ep+
      \frac{\Mbar\rlamr[1]}
      {\norm{\xt[0]-\xstar}\sqrt{\alpha\Llamr[1]}}
      \sqrt{\max\left\{1,
        \frac{6\kar[1]\normsq{\xt[0]-\xstar}}{\ep}\right\}}
      \right\}}.
  \end{align*}
  \jim{In the nontrivial accuracy regime
    $\ep\leq6\kar[1]\normsq{\xt[0]-\xstar}$, the computation bound simplifies to
    $C_\ep=\bigO(\Mbar\rlamr[1]/(\alpha\sqrt{\ep})+N_\ep)$.}
\end{corollary}
Again, we make a few remarks regarding the results. First, to find an $\ep$-close solution, Corollary \eqref{cor:S-str-close} implies that the ACGD-S method has the same oracle complexity as the ACGD method, and has the same computation complexity as that of the lower computation complexity bound for solving a single strongly-convex linearly constrained problem \cite{ouyang2021lower}. Second, the iteration limit function $\St$ is again adaptive to the varying difficulty of the saddle-point subproblem \eqref{eq:slide-saddle} from different phases. Third, the rather complicated inner-loop stepsize choice in \eqref{stp:S-str-init} and \eqref{stp:S-str-inner} is the first among sliding algorithms, e.g. \cite{lan2016gradient,lanOuyang2016GradientSliding,lan2021graph,lan2022optimal},  to achieve both the optimal inner loop complexity of $O(1/\sqrt{\epsilon})$ and the optimal outer loop complexity of $\bigO(\sqrt{\kappa}\log(1/\ep))$  without restarting. It is unclear if the same effect is achievable with simpler stepsize choices.
\jim{Since the stepsize choices require only conservative estimates of $\alpha$ and $\Llamr$ to achieve the optimal orders of computation and oracle complexity, its implementation does not require extensive parameter tuning.}

\subsection{Convergence Analysis of the ACGD-S Method}\label{subsec:S-proof}
We first establish a generic result for the $Q$-function defined in \eqref{eq:Q_decomp}, applicable to both the non-strongly convex and strongly convex cases.

\begin{proposition}\label{pr:S-Q}
  Consider an $\alpha$-strongly convex problem of form \eqref{eq:prob}. Let a set of reference multipliers $\Lambda \in \R_+^m$ be given and let the aggregate smoothness constant $\Llam$ be defined in \eqref{eq:agg-L}. Let iterates $\zt:=\{\xt; \lamtilt, \nut, \pit\}$ be generated by Algorithm \ref{alg:ANLA-S}.  Suppose the following stepsize requirements are met.  For all $t \geq 1$, the outer-loop stepsize requirements are
  \begin{align}\label{req:ANLA-s-outer}
    \begin{split}
       & \wt \etat \leq \wtt (\etatt + \alpha/2),                                         \\
       & \wt \taut \leq \wtt (\tautt + 1),                                                \\
       & \etatt \taut \geq \theta_t\Llam \text{ with }  \thetat :=\wtt /\wt \jim{\leq 1}, \\
       & \etat[N](\taut[N] + 1) \geq \Llam.
    \end{split}
  \end{align}
  For all $t \geq 1, s \geq 1$, the intra-phase stepsize requirements are
  \begin{align}\label{req:intra-phase}
    \begin{split}
       & \delst (\betast + \alpha/2) \geq \delst[s+1] \betast[s+1],                                \\
       & \delst \gamst \geq \delst[s+1] \gamst[s+1],                                               \\
       & \gamst \betast[s+1] \geq \rhost[s+1] \norm{\nut}^2,\  \rhost[s+1] = \delst / \delst[s+1], \\
       & \gamst[\St[N]][N] (\betast[\St[N]][N] + \alpha/2) \geq  \normsq{\nut[N]}.
    \end{split}
  \end{align}
  For all $t \geq 1$, the inter-phase requirements are
  \begin{align}\label{req:inter-phase}
    \begin{split}
       & \wst[\St] (\betast[\St] + \alpha/2) \geq \wst[1][t+1] \betast[1][t+1],                                                  \\
       & \wst[\St] \gamst[\St] \geq \wst[1][t+1] \gamst[1][t+1],                                                                 \\
       & \gamst[\St][t] \betast[1][t+1] \geq \rhost[1][t+1] \jim{\norm{\nut[t]}^2},\  \rhost[1][t+1] = \wst[\St] / \wst[1][t+1],
    \end{split}
  \end{align}
  where $\wst:= \wt \delst / (\sums \delst) $ denotes the aggregate weights. Then for any reference point $z=(x; \lam,  \nu, \pi) \in X \times \Lam \times [V, \Pi]$ with $[V, \Pi]$ being defined in \eqref{eq:joint-domain}, we have
  \begin{align}\label{eq:ANLA-S-Q-conv}
    \begin{split}
      \sumwt Q(\zt; z) + \tfrac{\wt[N]}{2} (\etat[N]  + \tfrac{\alpha}{2}) \normsq{\xt[N] - x} & \leq \tfrac{\wst[1][1] \betast[1][1] + \wt[1]\etat[1]}{2} \normsq{\xt[0] - x}  + \tfrac{\wst[1][1]\gamst[1][1]}{2} \normsq{\lamst[0][1] - \lam} \\
                                                                                               & + \wt[1] \taut[1] [U_{f^*}(\pi; \pit[0]) + \lam^\top U_{g^*} (\nu; \nut[0])].
    \end{split}
  \end{align}
\end{proposition}

\begin{proof}
  We first establish a convergence bound for the inner loop within a phase. Fix $t \geq 1$. Consider the convergence of $\yst$. Since $u(y) + \etat \normsq{y-x^{t-1}}/2$ has a strong convexity modulus of $\alpha + \etat$, the $y$-prox mapping in Line 8 of Algorithm \ref{alg:ANLA-S} leads to a three point inequality (see Lemma 3.1 of \cite{LanBook}):
  {\skipdisplay
  \begin{align*}
     & \inner{\yst - x}{\htilst} + u(\yst) - u(x) +  \tfrac{\etat}{2}(\normsq{\yst - \xtt} - \normsq{x - \xtt})\notag                                             \\
     & \quad \jim{+}\tfrac{1}{2}[(\betast+\alpha+\etat) \normsq{x - \yst[s]} + \betast \normsq{\yst - \yst[s-1]} - \betast \normsq{\yst[s-1] -x }] \leq 0. \notag
  \end{align*}}
  Equivalently, we have
  \begin{align}
     & \inner{\yst - x}{\htilst[s]}  +
    \tfrac{1}{2}[(\betast + \alpha/2) \normsq{x - \yst[s]} + \betast \normsq{\yst - \yst[s-1]} - \betast \normsq{\yst[s-1] -x }]\notag                          \\
     & \quad  + u(\yst) - u(x) +\tfrac{1}{2}[\etat \normsq{\yst - \xtt} +  (\etat + \alpha/2)\normsq{\yst - x} - \etat\normsq{x - \xtt} ] \leq 0.\label{p2:Q_X}
  \end{align}
  In particular, the definition of $\htilst[s]$ in Line 7 of Algorithm \ref{alg:ANLA-S} implies
  $$
    \begin{aligned}
      \inner{\yst - x}{\htilst[s]} = & \inner{\yst - x}{\sumi \lamist[s] \nuit} - \inner{\yst - x}{\sumi (\lamist[s] - \lamist[s-1]) \nuit} \\
                                     & + \rhost[s] \inner{\yst[s] - \yst[s-1]}{\sumi (\lamist[s-1] - \lamist[s-2]) \nuit}                   \\
                                     & + \rhost[s] \inner{\yst[s-1] - x}{\sumi (\lamist[s-1] - \lamist[s-2]) \nuit} ,\ \forall s \geq 2,    \\
    \end{aligned}
  $$
  and
  $$
    \begin{aligned}
      \inner{\yst[1]- x}{\htilst[1]} = & \inner{\yst[1]- x}{\sumi \lamist[1] \nuit} - \inner{\yst[1]- x}{\sumi (\lamist[1] - \lamist[0]) \nuit} \\
                                       & + \rhost[1] \inner{\yst[1] - \yst[0]}{\sumi (\lamist[0] - \lamist[-1]) \nuitt}                         \\
                                       & +  \rhost[1] \inner{\yst[0] - x}{\sumi (\lamist[0] - \lamist[-1]) \nuitt}.
    \end{aligned}
  $$
  So, substituting them into \eqref{p2:Q_X}, summing up the resulting inequality with weight $\delst$, noting the stepsizes conditions in \eqref{req:intra-phase}, and utilizing Young's inequality,  we get
  \begin{align}\label{p2:y_conv}
    \begin{split}
      \sum_{s=1}^{\St} & \delst \left(\La(\yst; \lamst, \nut, \pit) - \La(x; \lamst, \nut, \pit) + \tfrac{[\etat \normsq{\yst - \xtt} +  (\etat + \tfrac{\alpha}{2})\normsq{\yst - x} - \etat\normsq{x - \xtt}]}{2}\right) \\
                       & + \delst[1] \rhost[1] \inner{\yst[0] - x}{\sumi (\lamist[0] - \lamist[-1]) \nuitt}-\delst[\St]  \inner{\yst[\St]- x}{\sumi (\lamist[\St] - \lamist[\St-1]) \nuit}                                 \\
      \leq             & \tsum_{s=2}^{\St} \tfrac{\delst[s-1]\gamst[s-1]}{2} \normsq{\lamst[s-1] - \lamst[s-2]} + \tfrac{\delst[1] \jim{(\rhost[1])^2} \normsq{\nut[t-1]} }{2\betast[1]} \normsq{\lamst[0] - \lamst[-1]}   \\
                       & -\tfrac{1}{2}[\delst[\St] (\betast[\St] + \alpha/2) \normsq{\yst[\St]- x}  - \delst[1] \betast[1] \norm{\yst[0] - x}^2].
    \end{split}\end{align}

  Next, consider the convergence of $\lamst$. The $\lam$-proximal mapping in Line 9 of Algorithm \ref{alg:ANLA-S} implies
  \begin{equation*}
    \La(\yst; \lam, \nut,  \pi^t) - \La(\yst; \lamst, \nut, \pi^t) + \tfrac{\gamst}{2} [\normsq{\lam- \lamst} + \normsq{\lamst- \lamst[s-1]}- \normsq{\lam- \lamst[s-1]}] \leq 0.
  \end{equation*}
  Due to the stepsize conditions in \eqref{req:intra-phase}, the $\delst$ weighted sum satisfies
  \begin{equation*}
    \begin{split}
       & \sums \delst[s][t][\La(\yst; \lam, \nut, \pi^t) - \La(\yst; \lamst, \nut, \pi^t)] + \tfrac{\gamst[\St] \delst[\St]}{2} \normsq{\lam- \lamst[\St]} \\
       & \quad + \sums \tfrac{\delst \gamst}{2}  \normsq{\lamst[s] - \lamst[s-1]} \leq \jim{\tfrac{\delst[1] \gamst[1]}{2}}\normsq{\lam- \lamst[0]}.
    \end{split}
  \end{equation*}
  Then, combining it with the $y$ convergence bound in \eqref{p2:y_conv}, we get
  \begin{align*}
    \begin{split}
      \sums & \delst \left(\La(\yst; \lam, \nut, \pi^t) - \La(x; \lamst, \nut, \pi^t)  + \tfrac{{\eta_t} \normsq{\yst - \xtt} + (\eta_t + \alpha/2) \normsq{\yst - x} - \etat \normsq{\xtt -x}}{2}\right)                                                                                 \\
            & +\delst[1] \rhost[1] \inner{\yst[0] - x}{\sumi \nuitt (\jim{\lamist[0]} - \lamist[-1])} - \delst[\St]\inner{\yst[\St] - x}{\sumi \nuit (\lamist[\St] - \lamist[\St - 1])}                                                                                                   \\
      \leq  & \tfrac{\delst[1]\gamst[1]}{2} \normsq{\lam- \lamst[0]} - \tfrac{\delst[\St] \gamst[\St]}{2} [\normsq{\lam-\lamst[\St]} + \normsq{\lamst[\St]- \lamst[\St-1]}]  +\tfrac{\delst[1] \jim{(\rhost[1])^2} \normsq{\nut[t-1]}}{2\betast[1]} \normsq{\jim{\lamst[0]} - \lamst[-1]} \\
            & -\tfrac{1}{2}[\delst[\St] (\betast[\St] + \alpha/2) \normsq{\yst[\St]- x}  - \delst[1] \betast[1] \norm{\yst[0] - x}^2].
    \end{split}
  \end{align*}
  Moreover, since $\La(\yst; \lam, \nut, \pi^t)$, $\normsq{\yst - \xtt}$ and $\normsq{\yst - x}$ are convex with respect to $\yst$ and
  $\La(x; \lamst, \nut, \pi^t)$ is linear with respect to $\lamst$, multiplying both sides by $\wt/(\sums \delst)$ and applying the Jensen's inequality leads to
  \begin{align*}
    \begin{split}
      \wt & \left(\La(\xt; \lam, \nut, \pi^t) - \La(x; \lamtilt, \nut, \pi^t)  + \tfrac{{\eta_t} \normsq{\xt - \xtt} + (\eta_t + \alpha/2) \normsq{\xt - x} - \etat \normsq{\xtt -x}}{2}\right)                                                                                       \\
          & +\wst[1] \rhost[1] \inner{\yst[0] - x}{\sumi \nuitt (\jim{\lamist[0]} - \lamist[-1])} - \wst[\St]\inner{\yst[\St] - x}{\sumi \nuit (\lamist[\St] - \lamist[\St - 1])}                                                                                                     \\
          & \leq \tfrac{\wst[1]\gamst[1]}{2} \normsq{\lam- \lamst[0]} - \tfrac{\wst[\St] \gamst[\St]}{2} [\normsq{\lam-\lamst[\St]} + \normsq{\lamst[\St]- \lamst[\St-1]}]  +\tfrac{\jim{\wst[1](\rhost[1])^2} \normsq{\nut[t-1]}}{2\betast[1]} \normsq{\jim{\lamst[0]} - \lamst[-1]} \\
          & \ -\tfrac{1}{2}[\wst[\St] (\betast[\St] + \alpha/2) \normsq{\yst[\St]- x}  - \wst[1] \betast[1] \norm{\yst[0] - x}^2],
    \end{split}
  \end{align*}
  where $\wst = \wt \delst/ (\sums \delst)$ represents the aggregate weight for the inner iterates.

  Next, we consider the inner loops from different phases. The inter-phase stepsize condition in \eqref{req:inter-phase} implies the sum of preceding inequality across $t$ satisfies
    {\skipdisplay
      \begin{equation}\label{eq:sliding-QxQp}
        \begin{split}
          \sumt & {\wt} [Q_x(\zt; z)  + Q_\lam(\zt; z)] + \sumt \tfrac{\wt}{2} \etat\normsq{\xt - \xtt} + \jim{\tfrac{\wt[N]}{2}}(\etat[N] + \alpha/2) \normsq{\xt[N] - x}             \\
                & \leq \tfrac{\wst[1][1]}{2}(\gamst[1][1]\normsq{\jim{\lamst[0][1]} - \lam} + \betast[1][1] \normsq{\yst[0][1] - x}) + \tfrac{\wt[1] \etat[1]}{2} \normsq{\xt[0] - x}.
        \end{split}
      \end{equation}}
  Observe that \eqref{eq:sliding-QxQp} is almost identical to the $Q_x$ and $Q_\lam$ inequality in \eqref{eq:pr-pf-Qx}. Thus a similar argument to Proposition \ref{pr:Q-convergence} and the outer-loop stepsize requirements in \eqref{req:ANLA-s-outer} leads to the desired convergence result in \ref{eq:ANLA-S-Q-conv}.
\end{proof}

Now we use the preceding proposition to prove the convergence of the ACGD-S method under the non-strongly convex setting.

\textbf{Proof to Theorem \ref{thm:ANLA-S-non} and Corollary \ref{cor:S-non}:}
It is straightforward to verify that the outer loop stepsize in \eqref{stp:S-non-outer} satisfies the condition \eqref{req:ANLA-s-outer}, and the adaptive inner loop stepsize in \eqref{stp:S-non-inner} satisfies both the intra-phase condition \eqref{req:intra-phase} and the inter-phase condition \eqref{req:inter-phase}. So it follows from \eqref{eq:ANLA-S-Q-conv} that
\begin{align}\label{pfeq:thm-S-non-1}
  \begin{split}
    \sumwt Q(\zt; z) + \tfrac{\wt[N]\etat[N]}{2} \normsq{\xt[N] - x} & \leq \tfrac{\wst[1][1] \betast[1][1] + \wt[1]\etat[1]}{2} \normsq{\xt[0] - x}  + \tfrac{\wst[1][1]\gamst[1][1]}{2} \normsq{\lamst[0][1] - \lam}.
  \end{split}
\end{align}
Now consider setting the reference point to $z^* = (\xstar; \lamstar, \nustar = \grad g(\xstar), \pistar= \grad f(\xstar))$ such that $Q(\zt, \zstar) \geq 0\ \forall t$. The preceding inequality implies that $$\normsq{\xt[N] - \xstar} \leq \tfrac{1}{\Llamr[c]}[({\wst[1][1] \betast[1][1] + \wt[1]\etat[1]}) \normsq{\xt[0] - \jim{\xstar}}  + {\wst[1][1]\gamst[1][1]} \normsq{\lamst[0][1] - \jim{\lamstar}}]\ \forall N \geq 1.$$  So $\xundert$, being the convex combination of $\{\xt\}$s, remains in a bounded ball around $\xstar$, and $\nut = \grad g(\xundert)$ is bounded for all $t \geq 1.$
Next, setting the ergodic average solution as $\zbart:= (\xbarN; \sumt \wt \lamtilt/ (\sumt\wt),\nubarti[N][], \pibart[N])$  with
\begin{equation}\label{eq:ergo_sol}
  \pibart[N]:= \sumwt \pit / (\sumwt),\nutibar[N]:=  \begin{cases}\sumwt \lamti \nuit /(\sumwt \lamti) & \text{o.w.}                       \\
             \grad g_i(x^0)                       & \text{if } \lamti = 0\ \forall t,\end{cases}\end{equation}
a similar application of the Jensen's inequality as \eqref{pr1:jensen} and the stepsize choice in \eqref{stp:S-non-outer} and \eqref{stp:S-non-inner} lead to
\begin{align*}
  \wst[1][1] = \wt[1] / \St[1] = 1 / (\Mtilt[1] \Delta) =\Llamr[c]/\Mtilt[1] R,\text{such that } \wst[1][1]\betast[1][1] = \Llamr[c],\ \wst[1][1]\gamst[1][1] = \Llamr[c] / R^2,
\end{align*}
and that
\begin{equation}\label{thm4:Q-ave}
  Q(\zbart[N]; (\xstar; \lam, \nu, \pi))  \leq \tfrac{\Llamr[c]}{N(N+1)} (3\normsq{\xt[0] - \xstar} + \rlamr[c]^2/R^2), \forall \lam \in  \Lamr[c], (\nu,\pi) \in [V, \Pi].  \end{equation}
The convergence in both the optimality gap and the feasibility violation in \eqref{eq:S-non-strong-conv} then follows from Lemma \ref{lm:fun-feas-from-Q}.

We now show the oracle and algebraic operations complexity in \eqref{eq:S-non-strong-oracle-complexity-R}. Since only the outer loop requires gradient evaluations, the upper bound on $N_\epsilon$ follows directly from \eqref{eq:S-non-strong-conv}.
Moreover, since each inner iteration requires fewer than three matrix-vector multiplications,  the total number of matrix-vector multiplications across $N_\ep$ phases can be bounded as
\begin{align*}
  C_\ep \leq 3 \tsum_{t=1}^{N_\ep}S_t
  \leq 3N_\ep + \tfrac{3}{2}N_\ep(N_\ep+1)\Mbar\Delta
  = O\left\{N_\ep+\jim{\tfrac{\Mbar R}{\Llamr[c]}}
  +\tfrac{\Mbar}{\ep}\left(R\normsq{\xt[0]-\xstar}
  +\tfrac{\rlamr[c]^2}{R}\right)\right\}.
\end{align*}

\endproof

The next proof considers the strongly convex case.
\vgap\\
\textbf{Proof to Theorem \ref{thm:ANLA-S-str}, Corollary \ref{cor:S-str-opt} and Corollary \ref{cor:S-str-close}:}
It is straightforward to check that the outer-loop stepsize in \eqref{stp:S-str-outer} satisfies the condition \eqref{req:ANLA-s-outer}, and the adaptive inner-loop stepsize in \eqref{stp:S-str-init} and \eqref{stp:S-str-inner} satisfy the intra-phase condition \eqref{req:intra-phase}. Now we verify the inter-phase condition in \eqref{req:inter-phase}. Consider a fixed $t\geq 2$, we have
\begin{align*}
   & \wst[1][t] = \tfrac{\wt\delst[1]}{\sums \delst} = \tfrac{\wt \delst[1]}{\wt \MtU^2 \Gamt \Delta } = \tfrac{\delst[1]}{\MtU^2\Gamt \Delta} = \tfrac{\Wt}{\MtU \Gamt \Delta}, \\ &\wst[\St[t-1]][t-1]= \tfrac{\wt[t-1]\delst[\St[t-1]][t-1]}{\sums[\St[t-1]] \delst[s][t-1]} = \tfrac{\wt[t-1] \delst[\St[t-1]][t-1]}{\wtt \MtU[t-1]^2 \Gamt[t-1] \Delta } = \tfrac{\delst[\St[t-1]][t-1]}{\MtU[t-1]^2\Gamt[t-1] \Delta} = \tfrac{\Wt}{\MtU[t-1] \Gamt[t-1] \Delta}.
\end{align*}
Thus
\begin{align*}
  \wst[\St[t-1]][t-1] (\betast[\St-1][t-1] + \alpha/2 ) & = \tfrac{\Wt}{\MtU[t-1] \Gamt[t-1] \Delta} \tfrac{\alpha}{4}(\Wt[t-1] \Gamt[t-1] \MtU[t-1] + \St) \geq  \tfrac{\Wt}{\MtU[t-1] \Gamt[t-1] \Delta} \tfrac{\alpha}{4}(\delst[\St[t-1]][t-1]\Gamt[t-1])                     \\
                                                        & = \tfrac{\Wt}{\MtU[t-1] \Gamt[t-1] \Delta} \tfrac{\alpha}{4}(\MtU[t-1] \Wt \Gamt[t-1])= \tfrac{\alpha}{4} \tfrac{\Wt^2}{\Delta} = \tfrac{\Wt}{\MtU \Gamt \Delta} [\tfrac{\alpha}{4} \MtU \Gamt \Wt]= \wst[1]\betast[1].
\end{align*}
\begin{align*}
  \wst[\St[t-1]][t-1] \gamst[\St[t-1]][t-1] = \tfrac{\delst[\St[t-1]][t-1]}{\MtU[t-1]^2\Gamt[t-1] \Delta} (\tfrac{4}{\alpha}) \tfrac{\MtU[t-1]^2 \Gamt[t-1] }{\delst[\St[t-1]][t-1]} = \tfrac{4}{\alpha \Delta} = \tfrac{\delst[1]}{\MtU^2\Gamt \Delta} \tfrac{\MtU[t]^2 \Gamt[t] }{\delst[1][t]} = \wst[1][t] \gamst[1].
\end{align*}
\begin{align*}
  \betast[1] \gamst[\St[t-1]][t-1]  = (\tfrac{\alpha}{4} \MtU \Gamt \Wt) (\tfrac{4}{\alpha}) \tfrac{\MtU[t-1]^2 \Gamt[t-1] }{\MtU[t-1] \Wt} =  \MtU[t]\MtU[t-1] \Gamt \Gamt[t-1] \stackrel{(a)}{\geq} \tfrac{\MtU \Gamt}{\MtU[t-1] \Gamt[t-1]} \MtU[t-1]^2 = \tfrac{\wst[\St[t-1]][t-1]}{\wst[1]} \MtU[t-1]^2 = \rhost[1] \MtU[t-1]^2,
\end{align*}
where  the inequality in $(a)$ holds because we have $\Gamt[t-1] \geq 1$ as a consequence of its definition.
Thus all the requirements in Proposition \ref{pr:S-Q} are satisfied. We get from \eqref{eq:ANLA-S-Q-conv} that
\begin{align}\label{pfeq:thm-S-str-1}
  \begin{split}
    \sumwt Q(\zt; z) + \tfrac{\wt[N]\etat[N]}{2} \normsq{\xt[N] - x} & \leq \tfrac{\wst[1][1] \betast[1][1] + \wt[1]\etat[1]}{2} \normsq{\xt[0] - x}  + \tfrac{\wst[1][1]\gamst[1][1]}{2} \normsq{\lamst[0][1] - \lam}.
  \end{split}
\end{align}
Here a straightforward calculation of the stepsize choice in \eqref{stp:S-str-outer} and \eqref{stp:S-str-inner} leads to
\begin{align}\label{pfeq:thm-S-str-2}
  \wt[1]\etat[1] \leq 2L(\Lamr[c]), \wst[1][1] \betast[1][1] = 0,\text{ and } \wst[1][1] \gamst[1][1] = \jim{4/(\Delta \alpha)} = 4L(\Lamr[c])/R^2.
\end{align}
Similar arguments as that of Theorem \ref{thm:ANLA-S-non} imply the boundedness of $\MtU$, and that the ergodic average solution $\zbart$ defined according to \eqref{eq:ergo_solns} satisfies

\begin{equation}\label{thm5:Q-ave}
  Q(\zbart[N]; (\xstar; \lam, \nu, \pi))  \leq \tfrac{\Llamr[c]}{\sumwt} ( \tfrac{2\rlamr[c]^2}{R^2} +  \normsq{\xt[0] - \xstar}), \forall \lam \in  \Lamr[c], (\nu,\pi) \in [V, \Pi].  \end{equation}
Since $\sumwt \geq \max\{N(N+1)/2,
  \jim{\sqrt{2\kar[c]}[(1+ 1/\sqrt{2\kar[c]})^{N-4} -1]}\}$
(see \eqref{eq:wt_lower_bound}), we get the optimality gap and feasibility violation convergence bound in \eqref{eq:S-str-conv}. Moreover, since $\tfrac{\alpha}{2} \normsq{\xbart[N] - \xstar} \leq  Q(\zbart; (\xstar; \lamstar, \grad g(\xbart[N]), \grad f(\xbart[N])))$ (see \eqref{thmpf:Q-strong}), the convergence of $\xbart[N]$ to $\xstar$ in \eqref{eq:S-str-conv} also follows from \eqref{thm5:Q-ave}.
Next we establish the oracle and computation complexity bounds in Corollary \ref{cor:S-str-opt}. For any
$$N \geq \jim{\left\lceil
    (\sqrt{2\kar[c]}+1)
    \log\left(1+\frac{\sqrt{\Llamr[c]\alpha}
      (\normsq{\xt[0]-\xstar}+2\rlamr[c]^2/R^2)}
    {\sqrt{2}\,\ep}\right)
    \right\rceil+4}$$
we get
\begin{align*}
  Q(\zbart[N]; (\xstar; \lam, \nu, \pi)) \leq \ep, \forall \lam \in  \Lamr[c], (\nu,\pi) \in [V, \Pi], \\
  \Rightarrow \max\{F(\xbart[N]) - \Fstar, c \norm{[g(\xbart[N])]_+}\}  \leq \ep,
\end{align*}
so that $\xbart[N]$ is an $(\ep; \ep/c)$ solution. Therefore, the least number of phases $\Nep$ required for such a solution admits the upper bounded in the corollary statement.

Now we consider the corresponding number of matrix-vector multiplications $\Cep$ in the $\Nep$ phases, i.e., $\Cep = \sumt[\Nep] \St[t].$
\jim{Let
  \[
    A:=\tfrac{2\rlamr[c]^2}{R^2}+\normsq{\xt[0]-\xstar},
    \qquad
    T_\ep:=\max\left\{1,\frac{\Llamr[c]A}{\ep}\right\},
  \]
  and define $\Nep$ to be the first index for which $\sum_{t=1}^{\Nep}\wt\geq T_\ep$. Then $\xbart[\Nep]$ is an $(\ep;\ep/c)$-optimal solution by~\eqref{thm5:Q-ave}. Moreover, since $\kar[c]\geq1$, the stepsize choice gives $\wt/\wt[t-1]=1/\thetat\leq 2$. Hence, if $\Nep\geq2$, the minimality of $\Nep$ yields
  \[
    \sum_{t=1}^{\Nep}\wt
    \leq \sum_{t=1}^{\Nep-1}\wt+2\wt[\Nep-1]
    \leq 3\sum_{t=1}^{\Nep-1}\wt
    <3T_\ep.
  \]
  The same bound is immediate when $\Nep=1$.}
Calculating the sum of $\St$ directly is challenging, so we consider an easier quantity $\calRt[t] := [\St - 3]_+.$ An useful algebraic relation is
\begin{equation}\label{pfeq:Rt-lower-bound}
  \tsum_{i=1}^{t-1} \calRt[i] \leq \Mbar \Wt.\end{equation}
The result can be deduced by induction.
\jim{For $t=2$, Lemma~\ref{lm:alg-fact-St} gives
  \[
    \Wt[2]\Mbar\geq \Wt[2]\MtU[1]
    =\delst[\St[1]][1]=\frac{\St[1]}{\Gamt[1]}
    \geq \St[1]-2\geq\calRt[1].
  \]
  Assuming \eqref{pfeq:Rt-lower-bound} holds for some $t\geq2$, we have}
\begin{align*}
  \Wt[t+1] = \tfrac{\Wt \MtU + \tfrac{1}{\Gamt}(\St-1)}{\MtU} \stackrel{(a)}{\geq} \Wt[t] + \tfrac{[\St[t] - 3]_+}{\MtU} \stackrel{(b)}{\geq} \tfrac{1}{\Mbar}\{\tsum_{i=1}^{t-1} \calRt[i] + [\St[t] -3]_+\} = \jim{\tfrac{1}{\Mbar} \tsum_{i=1}^{t}\calRt[i]},
\end{align*}
where (a) follows from Lemma \ref{lm:alg-fact-St} and $\Gamt\geq1$, and (b) follows from the induction hypothesis and $\MtU\leq\Mbar$. Thus the principle of mathematical induction implies that \eqref{pfeq:Rt-lower-bound} is valid. Consequently, for $\calRt[\ep] = \sumt[\Nep] \calRt[t]$,  we have
\begin{align*}
  \calRt[\ep]^2/2 & \leq  \sums[\calRt[\ep]] s = \sumt[\Nep] \sums[\calRt] [(\tsum_{j=1}^{t-1} \calRt[j]) + s] \\
                  & \leq \sumt[\Nep] \sums[\calRt] [\Mbar \Wt + s]
  \leq \sumt[\Nep] \sums[\St-1] [\Mbar \Wt +  (s-1)]
  \leq \Mbar \sumt[\Nep] \sums[\St-1] [\Wt +  \tfrac{(s-1)}{\MtU}]                                             \\
                  & \stackrel{(a)}{\leq} \Mbar \sumt[\Nep] \MtU \wt \Delta
  \leq \Mbar^2\Delta\sumt[\Nep]\wt
  \jim{\leq 3\Mbar^2\Delta T_\ep}.
\end{align*}
where (a) follows from the minimality of $\St$ in \eqref{stp:St-strong}. Therefore,
\jim{
  \[
    \calRt[\ep]
    =O\left\{
    \frac{\Mbar R}{\sqrt{\alpha\Llamr[c]}}
    +\sqrt{\frac{1}{\alpha\ep}}\Mbar
    \bigl(R\norm{\xt[0]-\xstar}+\rlamr[c]\bigr)
    \right\}.
  \]}
Since $\sumt[\Nep] \St \leq \calRt[\ep] + 3 \Nep$, the
\jim{computation bound in Corollary~\ref{cor:S-str-opt} follows because its
  small-accuracy assumption implies $\ep\leq\Llamr[c]A$, so the
  $\ep$-independent initialization term is dominated. The oracle bound follows
  from $\log(1+x)\leq\log(2x)$ for $x\geq1$.}

Next, setting $c=1$, the oracle and computation complexity bounds in Corollary \ref{cor:S-str-close} can be derived similarly.

\endproof

\subsection{The Adaptive ACGD-S Method}\label{sec:adaptive-ACGD-S}

\jim{Throughout this subsection, we impose Assumption~\ref{ass:bounded_domain}. The diameter $D_X$ is used only in the complexity analysis and need not be known by the algorithm.}

Similar to the ACGD method, the proposed ACGD-S method's stepsize choice requires an upper bound to the aggregate Lipschitz smoothness constant $\Llamr[c]$ to ensure convergence, but, just as before, $\Llamr[c]$ is notoriously difficult to estimate because it  depends on the optimal dual multiplier. Similar to Subsection \ref{subsec:ANLA-search},  we propose a verifiable termination certificate called the PD-gap associated with a parameter $r$ to develop the adaptive ACGD-S method shown in Algorithm \ref{alg:adaptive-ACGD-S}



\begin{algorithm}[htb]
  \caption{The Adaptive ACGD-S Method}
  \label{alg:adaptive-ACGD-S}
  \begin{algorithmic}[1]

    \Require $\xt[-1]= \xundert[0] = \jim{\yst[0][1]}=\xt[0]\in X$, \jim{the scaling constant $c\geq1$, the strong convexity modulus $\alpha$, initial estimates $\Ltil_{f,0}\geq\alpha$ and $\Ltil_{g,0}\geq0$, $\Ltil_f=\Ltil_{f,0}$, $\Ltil_g=\Ltil_{g,0}$, $\rtil=2c$, and $\Ltilagg=\tilde L_{\mathrm{agg},0}:=\Ltil_{f,0}+2c\Ltil_{g,0}>0$}.

    \State Set $\pit[0] = \grad f(\xt[0]) $, $\nut[0] = \grad g(\xt[0])$, $\lamst[-1][1] = \lamst[0][1] = 0 $.
    \parState{\label{aACDG-S-line:compute-outer-stepsize} Compute stepsize parameter $\wt$, $\etat$, $\thetat$, and $\taut$ according to Theorem \ref{thm:ANLA-S-non} (Theorem \ref{thm:ANLA-S-str}) using $\Ltilagg$ and $R=1$  as proxies for problem parameters for the non-strongly convex (strongly convex) case.}
    \For{$t = 1,2,3,\ldots$}
    \parState{Set $\xundert \leftarrow (\taut \xundert[t-1] + \xtilt)/(1 + \taut)$ where  $\xtilt  = \xtt + \thetat (\xtt - \xt[t-2])$.}
    \parState{Set $\pit \leftarrow \grad f(\xundert)$ and $\nut \leftarrow \grad g(\xundert)$.}
    \parState{\label{aACDG-S-line:compute-inner-stepsize} Calculate inner loop iteration limit $\{\St\}$ stepsizes $\{\betast\}$ and $\{\gamst\}$,  and weights $\{\delst\}$ according to Theorem \ref{thm:ANLA-S-non} (Theorem \ref{thm:ANLA-S-str}) using $\Ltilagg$ and $R=1$ for the non-strongly convex (strongly convex) case.}
    \For{$s=1,2,..., \St$}
    \parState{Set $\htilst = \begin{cases} (\nut)^\top\lamst[0] + \rhost[1] (\nut[t-1])^\top (\lamst[0] - \lamst[-1]) & \text{if } s =1 \\
              (\nut)^\top\lamst[s-1] + \rhost (\nut)^\top (\lamst[s-1] - \lamst[s-2])    & \text{ o.w.}\end{cases}$}
    \parState{Set $\yst \leftarrow \argmin_{y \in X}\inner{\htilst + \pit}{y} + u(y) + \etat\normsq{y - \xtt}/2 +  \betast \normsq{y - \yst[s-1]}/2. $}
    \parState{Set $\lamst \leftarrow \argmax_{\lam \in \R^m_+}\inner{\lam}{\nut (\yst - \xundert) + g(\xundert)} - \gamst \normsq{\lam - \lamst[s-1]}/2.$}
    \EndFor
    \parState{Set $\lamst[0][t+1] = \lamst[\St][t]$, $\lamst[-1][t+1]= \lamst[\St-1][t]$, $\yst[0][t+1] = \yst[\St][t]$}
    \parState{Set $\xt = \sums \delst \yst /(\sums \delst)$ and $\lamtilt = \sums \delst \lamst /(\sums \delst)$.}
    \parState{\label{aACDG-S-line:ergodic-average}\jim{Compute the ergodic average $\xbar^t:= \tsum_{j=1}^t \wt[j] \xt[j] / (\tsum_{j=1}^t \wt[j])$.}}

    \parState{\label{aACDG-S-line:update-Ltil}Update the Lipschitz smoothness constant estimates:
      $$\Ltil_f\leftarrow \max\{\Ltil_f,  \tfrac{1}{2}\tfrac{\normsq{\grad f(\xunder^t)-\grad f(\xunder^{t-1})}}{f(\xunder^{t-1})-f(\xunder^t)-\inner{\grad f(\xunder^{t})}{\xunder^{t-1}-\xunder^{t}}}, \jim{\tfrac{1}{2}\tfrac{\normsq{\grad f(\xbar^t)-\grad f(\xunder^{t})}}{f(\xunder^t)-f(\xbar^t)-\inner{\grad f(\xbar^t)}{\xunder^{t}-\xbar^{t}}}}\},$$
      \begin{align*}\Ltil_g\leftarrow \max\left\{\Ltil_g, \tfrac{1}{2} \norm{[\tfrac{\normsq{\jim{\grad g_1(\xunder^t)-\grad g_1(\xunder^{t-1})}}}{g_1(\xunder^{t-1})-g_1(\xunder^t)-\inner{\grad g_1(\xunder^{t})}{\xunder^{t-1}-\xunder^{t}}},...,\tfrac{\normsq{\jim{\grad g_m(\xunder^t)-\grad g_m(\xunder^{t-1})}}}{g_m(\xunder^{t-1})-g_m(\xunder^t)-\inner{\grad g_m(\xunder^{t})}{\xunder^{t-1}-\xunder^{t}}}]},\right. \\
        \left.\jim{\tfrac{1}{2} \norm{[\tfrac{\normsq{\grad g_1(\xbar^t)-\grad g_1(\xunder^{t})}}{g_1(\xunder^t)-g_1(\xbar^t)-\inner{\grad g_1(\xbar^t)}{\xunder^{t}-\xbar^{t}}},...,\tfrac{\normsq{\grad g_m(\xbar^t)-\grad g_m(\xunder^{t})}}{g_m(\xunder^t)-g_m(\xbar^t)-\inner{\grad g_m(\xbar^t)}{\xunder^{t}-\xbar^{t}}}]}}\right\}.
      \end{align*}
    }
    \parState{\label{aACDG-S-line:calculate-FP-gap} Calculate the PD-gap $\DeltaPD(\rtil):=F(\xbar^t) + \rtil \norm{[g(\xbar^t)]_+} - \FunderPD^t$, where $\FunderPD^t$ is calculated using the lower linear approximation functions in \eqref{eq:ACGD-S-Funder}.}
    \If{\label{aACDG-S-line:smoothness-restart}$\Ltilagg<\Ltil_f+\rtil\Ltil_g$}
    \State Restart with $\Ltilagg\leftarrow\max\{2\Ltilagg,\Ltil_f+\rtil\Ltil_g\}$.
    \EndIf
    \While{\label{aACDG-S-line:restart-trigger}$\norm{[g(\xbar^t)]_+} > 2\DeltaPD(\rtil)/\rtil$}
    \State Set $\rtil^+\leftarrow2\rtil$.
    \If{$\Ltilagg<\Ltil_f+\rtil^+\Ltil_g$}
    \State Restart with $\Ltilagg\leftarrow\max\{2\Ltilagg,\Ltil_f+\rtil^+\Ltil_g\}$ and $\rtil\leftarrow\rtil^+$.
    \Else
    \State \jim{Set $\rtil\leftarrow\rtil^+$ and update $\DeltaPD(\rtil)$.}
    \EndIf
    \EndWhile
    \EndFor
  \end{algorithmic}
\end{algorithm}


\begin{definition}\label{def:primal-dual-gap-certificate}
  Given an evaluation point $\xbar$ and some lower linear approximation $l_f$ to $f$ and the lower linear approximations $l_{g_1}, l_{g_2}, \ldots, l_{g_m}$ to the non-negative scaled constraint functions, i.e., $\lambda_1 g_1, \lambda_2 g_2, \ldots, \lambda_m g_m$ with $\lambda_i \geq 0 \forall i \in [m]$, the PD-gap parameterized by $r$, denoted $\DeltaPD(r)$, is defined as:
  \begin{equation}\label{eq:PD-gap}
    \begin{split}
       & \DeltaPD(r) := F(\xbar) + r\norm{[g(\xbar)]_+} -\FunderPD                           \\
       & \text{ where } \FunderPD := \min_{x \in X} l_f(x) + u(x) + \sum_{i=1}^m l_{g_i}(x).
    \end{split}
  \end{equation}
\end{definition}

\jim{Under Assumption~\ref{ass:bounded_domain}, the minimum defining $\FunderPD$ is finite, so the PD-gap is well-defined.}

The PD-gap is more closely related to the primal dual iterates generated by the ACGD-S method, hence the name \textbf{P}rimal \textbf{D}ual gap.
In contrast to $\Funder$, $\FunderPD$ is obtained by minimizing a Lagrangian relaxation of the problem defining $\Funder$ in \eqref{eq:FP-gap}. Consequently, it is a weaker lower bound, satisfying $\FunderPD \leq \Funder$.
\jim{However, computing $\FunderPD$ is computationally cheaper than computing $\Funder$ because it requires only an $X$-projection of the form \eqref{eq:X-proj-def}, whereas computing $\Funder$ requires a solving a linearly constrained subproblem. Under our standing assumption that the $X$-projection is easy, it is typically less expensive. For example, when $u$ is quadratic and $X$ is a box or a Euclidean ball, the $X$-projection has a closed-form solution.}

Next, we show that the proposed PD-gap satisfies the requirements for a verifiable termination certificate stated in Property \ref{propy:verifiable-certificate-requirements} and can therefore be incorporated into the ACGD-S method in the same manner as the FP-gap is incorporated into the ACGD method.
As discussed above, given $l_f$, $\{l_{g_i}\}_{i=1}^m$, and the parameter $r$, the PD-gap can be computed without knowledge of any unknown problem parameters; hence, Property \ref{propy:verifiable-certificate-requirements}.a) is satisfied. Since the PD-gap is always at least as large as the corresponding FP-gap, Property \ref{propy:verifiable-certificate-requirements}.c) follows immediately.
\begin{lemma}\label{lm:PD-gap-to-suboptimality-feasibility}
  For the evaluation point $\xbar$, the PD-gap yields the following upper bounds on feasibility violation and optimality gap:
  \begin{enumerate}
    \item[a)] $F(\xbar) - \Fstar\leq \DeltaPD(r)$ for any $r\geq 0$.
    \item[b)] If $r \geq 2\jim{\rplus_c}$, then
          \begin{equation}\label{eq:PD-gap-to-feasibility-violation}
            \norm{[g(\xbar)]_+} \leq 2\DeltaPD(r)/r \jim{\leq\DeltaPD(r)/c}.
          \end{equation}
  \end{enumerate}
\end{lemma}

For Property \ref{propy:verifiable-certificate-requirements}.b), just like Lemma \ref{lm:FP-gap-from-Q}, we have the following lemma that shows the convergence of the PD-gap is implied by the convergence of the $Q$-gap associated with the dual variables $\Lambda\subset B^m_+(0;r)$ in the ACGD-S method.
\begin{lemma}\label{lm:PD-gap-from-Q}
  Given some primal dual iterates $\{(\xt[j],\tilde{\lambda}^j, \nut[j], \pit[j])\}_{j=1}^t$ and some nonnegative weights $\jim{\{\wt\}_{j=1}^t}$, let $\xbar^t$ be the weighted average of these iterates, i.e., $\xbar^t = \tsum_{j=1}^t \wt[j] \xt[j] /\tsum_{j=1}^t \wt[j]$. If the $Q$-gap associated with them is bounded by some $\Delta(r)$ for $\hat{\pi}^t = \grad f(\xbar^t)$, $\hat{\nu}^t = \grad g(\xbar^t)$,
  \begin{equation*}
    \max_{x \in X, \lambda \in B^m_+(0;r)}  \frac{1}{\tsum_{j=1}^t \wt[j]}\tsum_{j=1}^t \wt[j]\left(\La(\xt[j]; \lambda, \hat{\nu}^t, \hat{\pi}^t) - \La(x;\tilde{\lambda}^j,\nut[j],\pit[j])\right)\leq \Delta(r).
  \end{equation*}
  Then the following lower linear approximation functions $l_f$ and $l_{g_i}$ correspond to a PD-gap certificate associated with $\xbar^t$ where
  \begin{equation}\label{eq:ACGD-S-Funder}
    \begin{split}
      \text{ where } l_f(x;t) = \tfrac{1}{\tsum_{j=1}^{t} \wt[j]} \tsum_{j=1}^{t} \wt[j][\inner{\pit[j]}{x- \xundert[j]} + f(\xundert[j])].                           \\
      \text{ and } l_{g_i}(x;t) = \tfrac{1}{\tsum_{j=1}^{t} \wt[j] } \tsum_{j=1}^{t} \wt[j] \tilde{\lambda}^j_i[\inner{\nuit[j]}{x- \xundert[j]} + g_i(\xundert[j])]. \\
    \end{split}
  \end{equation}
  Importantly, the PD-gap associated with it satisfies $\DeltaPD(r) \leq \Delta(r)$.
\end{lemma}

\begin{proof}
  \jim{Let $W_t:=\tsum_{j=1}^t\wt[j]$. Since
  $\xbar^t=W_t^{-1}\tsum_{j=1}^t\wt[j]\xt[j]$, the convexity of $u$ and the
  conjugate identities for $\hat\pi^t=\grad f(\xbar^t)$ and
  $\hat\nu^t=\grad g(\xbar^t)$ give
  \begin{align*}
    F(\xbar^t)+r\norm{[g(\xbar^t)]_+}
     & =\max_{\lambda\in B_+^m(0;r)}
    \bigl\{F(\xbar^t)+\inner{\lambda}{g(\xbar^t)}\bigr\} \\
     & \leq\max_{\lambda\in B_+^m(0;r)}
    \frac{1}{W_t}\tsum_{j=1}^t\wt[j]
    \La(\xt[j];\lambda,\hat\nu^t,\hat\pi^t).
  \end{align*}
  On the other hand, the definitions in \eqref{eq:ACGD-S-Funder} imply
  \begin{align*}
    \FunderPD
     & =\min_{x\in X}\frac{1}{W_t}\tsum_{j=1}^t\wt[j]
    \left[\inner{\pit[j]}{x}-f^*(\pit[j])+u(x)
      +\tsum_{i=1}^m\tilde\lambda_i^j
    \bigl(\inner{\nut[j]_i}{x}-g_i^*(\nut[j]_i)\bigr)\right] \\
     & =\min_{x\in X}\frac{1}{W_t}\tsum_{j=1}^t\wt[j]
    \La(x;\tilde\lambda^j,\nut[j],\pit[j]).
  \end{align*}
  Subtracting these two relations and using the assumed $Q$-gap bound yields
  \begin{align*}
    \DeltaPD(r)
     & \leq\max_{x\in X,\,\lambda\in B_+^m(0;r)}
    \frac{1}{W_t}\tsum_{j=1}^t\wt[j]
    \left[\La(\xt[j];\lambda,\hat\nu^t,\hat\pi^t)
      -\La(x;\tilde\lambda^j,\nut[j],\pit[j])\right]
    \leq\Delta(r),
  \end{align*}
  which proves the claim.}
\end{proof}

\jim{Lemma~\ref{lm:PD-gap-from-Q} reduces the convergence analysis of the PD-gap
  to controlling the $Q$-gap over the current dual ball $B^m_+(0;\rtil)$. For ACGD-S, the choice
  of this ball requires more care than in Algorithm~\ref{alg:a-ACGD}.
  Algorithm~\ref{alg:a-ACGD} derives $\rtil$ from the current aggregate estimate
  through \eqref{eq:a-ANLA-r-lower-bound}, thereby selecting the largest dual
  radius supported by the current smoothness estimates. This radius can be much
  larger than that needed for the feasibility certificate, particularly when
  $\Ltil_g$ is small. Although this radius does not appear explicitly in the
  convergence bound for Algorithm~\ref{alg:a-ACGD}, it increases the complexity
  of the sliding mechanism. This motivates maintaining a separate radius estimate
  that is increased only when required by the feasibility certificate.}

\jim{To this end, Algorithm~\ref{alg:adaptive-ACGD-S} maintains $\rtil$ and
  $\Ltilagg$ as two separate nondecreasing quantities: $\rtil$ is doubled only
  when the feasibility test in Line~\ref{aACDG-S-line:restart-trigger} fails,
  whereas $\Ltilagg$ determines the stepsizes and must satisfy
  $\Ltilagg\geq\Ltil_f+\rtil\Ltil_g$. If the current $\Ltilagg$ supports an
  enlarged radius, the PD-gap is reevaluated without restarting; otherwise
  $\Ltilagg$ is increased and the method restarts. Thus a supported radius
  enlargement does not interrupt the sliding process, allowing the method to
  retain its accumulated weighted minorants and warm-started inner-loop states.}



\jim{The following proposition establishes the required PD-gap bound at the
  current certificate radius, including when that radius is increased within a
  run.}

\begin{proposition}\label{pr:a-ACGD-S-Q-convergence}
  \jim{Suppose Assumption~\ref{ass:bounded_domain} holds.} Consider the primal-dual iterates ${\zt:=(\xt; \lamtilt, \nut, \pit)}$ generated during a run of Algorithm~\ref{alg:adaptive-ACGD-S}. Let $\Ltilagg$ denote the fixed aggregate Lipschitz smoothness estimate used for this run, and let $\rtil_t$ denote the current value of the nondecreasing certificate radius. At iteration $t$, let $\Ltil_{f,t}$ and $\Ltil_{g,t}$ denote the current empirical smoothness estimates. If $\Ltilagg\geq\Ltil_{f,t}+\rtil_t\Ltil_{g,t}$, then the averaged solution $\xbar^t$ and the associated PD-certificate generated in Line~\ref{aACDG-S-line:calculate-FP-gap} satisfy
  \begin{equation}\label{eq:a-ACGD-S-Q-conv}
    \begin{split}
      \jim{\left(\tsum_{j=1}^t \wt[j]\right)\DeltaPD^{t}(\rtil_t)} \leq \max_{x \in X} 2\Ltilagg  [\normsq{\xt[0] - x} + 2(\rtil_t)^2],
    \end{split}
  \end{equation}
\end{proposition}

\begin{proof}
  \jim{Fix $t$ and set $\hat\pi^t:=\grad f(\xbar^t)$ and
    $\hat\nu^t:=\grad g(\xbar^t)$. Fix any reference point
    $z=(x;\lambda,\hat\nu^t,\hat\pi^t)$ with
    $x\in X$ and $\lambda\in B_+^m(0;\rtil_t)$. Retrace the proof of
    Proposition~\ref{pr:S-Q}. All intra-phase and inter-phase estimates are
    unchanged. Since $\norm{\lambda}\leq\rtil_t$, the hypothesis implies
    \[
      \Ltil_{f,t}+\norm{\lambda}\Ltil_{g,t}
      \leq\Ltil_{f,t}+\rtil_t\Ltil_{g,t}
      \leq\Ltilagg.
    \]
    For the outer-loop terms, the running-max definitions of
    $\Ltil_{f,t}$ and $\Ltil_{g,t}$ in
    Line~\ref{aACDG-S-line:update-Ltil}, together with the
    conjugate-Bregman identities, allow us to use the same Young-inequality
    argument as in the proof of Proposition~\ref{pr:a-Q-convergence}. In
    particular, for every $2\leq j\leq t$,
    \begin{align*}
       & \tau_j\left[\tsum_{i=1}^m\lambda_i
      U_{g_i^*}(\nu_i^j;\nu_i^{j-1})
      +U_{f^*}(\pi^j;\pi^{j-1})\right]             \\
       & \quad+\inner{\theta_j(\xt[j-1]-\xt[j-2])}
      {\pi^{j-1}-\pi^j+\tsum_{i=1}^m\lambda_i
      (\nu_i^{j-1}-\nu_i^j)}                       \\
       & \geq-\frac{\theta_j^2
        (\Ltil_{f,t}+\norm{\lambda}\Ltil_{g,t})}{2\tau_j}
      \normsq{\xt[j-1]-\xt[j-2]}                   \\
       & \geq-\frac{\theta_j\eta_{j-1}}{2}
      \normsq{\xt[j-1]-\xt[j-2]}.
    \end{align*}
    The last inequality follows from the preceding bound and
    $\eta_{j-1}\tau_j\geq\theta_j\Ltilagg$. For the remaining Bregman and
    inner-product terms at iteration $t$ involving $\hat\pi^t-\pi^t$ and
    $\hat\nu_i^t-\nu_i^t$, the corresponding argument in the proof of
    Proposition~\ref{pr:a-Q-convergence} gives the lower bound
    $-\eta_t\normsq{\xt[t]-\xt[t-1]}/2$, using
    $\eta_t(\tau_t+1)\geq\Ltilagg$. The displayed inequality and this lower
    bound remain valid if $\rtil$ is increased within the current run without
    triggering a restart because $\Ltil_{f,t}$ and $\Ltil_{g,t}$ are running
    maxima over all empirical ratios generated in the run, and the acceptance
    test verifies the same aggregate bound for the enlarged current radius
    $\rtil_t$. Consequently, the displayed inequality is valid
    for every $2\leq j\leq t$, including indices preceding such an increase,
    and the bound involving $\hat\pi^t$ and $\hat\nu^t$ is valid at iteration
    $t$.
    Finally, $\tau_1=0$ under both outer-loop stepsize choices, so the initial
    conjugate-Bregman term vanishes.

    Combining the displayed inequalities for $2\leq j\leq t$, the bound
    involving $\hat\pi^t$ and $\hat\nu^t$, and the unchanged intra-phase and
    inter-phase bounds gives
    \begin{align*}
      \sum_{i=1}^t \wt[i] Q(z^i;z)
      +\tfrac{\wt[t]}{2}(\etat[t]+\tfrac{\alpha}{2})
      \normsq{\xt[t]-x}
       & \leq \tfrac{\wst[1][1]\betast[1][1]+\wt[1]\etat[1]}{2}
      \normsq{\xt[0]-x}
      +\tfrac{\wst[1][1]\gamst[1][1]}{2}
      \normsq{\lamst[0][1]-\lambda}                             \\
       & \leq 2\Ltilagg\normsq{\xt[0]-x}
      +2\Ltilagg\normsq{\lambda}.
    \end{align*}
    It remains to verify the two coefficients on the first line of the
    preceding display. For $\alpha=0$, direct substitution of the stepsizes
    with $R=1$ gives
    \[
      \frac{\wst[1][1]\betast[1][1]+\wt[1]\etat[1]}{2}
      =\frac{3}{2}\Ltilagg,
      \qquad
      \frac{\wst[1][1]\gamst[1][1]}{2}
      =\frac{1}{2}\Ltilagg.
    \]
    For $\alpha>0$, the initialization and running-max update give
    $\Ltil_{f,t}\geq\alpha$, so the hypothesis implies
    $\Ltilagg\geq\alpha$. Hence the proxy condition number
    $\Ltilagg/\alpha$ is at least one, and the strongly convex stepsizes
    with $R=1$ give
    \[
      \frac{\wst[1][1]\betast[1][1]+\wt[1]\etat[1]}{2}
      =\Ltilagg,
      \qquad
      \frac{\wst[1][1]\gamst[1][1]}{2}
      =2\Ltilagg.
    \]
    Thus both coefficients are at most $2\Ltilagg$ in either case. Since
    $\lamst[0][1]=0$, the second inequality in the preceding display
    follows.

    Finally, drop the nonnegative terminal term and maximize over
    $x\in X$ and $\lambda\in B_+^m(0;\rtil_t)$. Lemma
    \ref{lm:PD-gap-from-Q} then yields  \eqref{eq:a-ACGD-S-Q-conv}.}
\end{proof}

\begin{theorem}\label{thm:aACGD-S}
  \jim{For a prescribed scaling constant $c\geq1$,} consider a smooth constrained optimization problem of the form \eqref{eq:prob} \jim{with a bounded domain satisfying Assumption \ref{ass:bounded_domain}}. \jim{Suppose the strong convexity modulus $\alpha$ is given,
    $\alpha\leq\Ltil_{f,0}\leq L_f$, $0\leq\Ltil_{g,0}\leq L_g$, and let
    \[
      \tilde L_{\mathrm{agg},0}:=\Ltil_{f,0}+2c\Ltil_{g,0}>0,
      \qquad \jim{\Lagg^c:=L_f+4\rplus_c L_g},
    \]
    and define
    \[
      K_r:=\ceil{\log_2\frac{\jim{\rplus_c}}{c}},\qquad
      K_L:=\ceil{\log_2\frac{2\jim{\Lagg^c}}
        {\tilde L_{\mathrm{agg},0}}}.
    \]
    Then $0<\tilde L_{\mathrm{agg},0}\leq\jim{\Lagg^c}$.}
  The adaptive ACGD-S method in Algorithm \ref{alg:adaptive-ACGD-S} finds an $(\ep,\ep/c)$-optimal solution for any $\ep>0$, where the final certificate radius satisfies \jim{$2c\leq\rtil<4\jim{\rplus_c}$}. It performs at most $K_r$ radius $\tilde{r}$ doublings and at most $K_L$ actual restarts. Across all runs, the method requires at most $N_\epsilon$ oracle evaluations and $C_\epsilon$ matrix-vector multiplications and projection operations in total.
  \begin{enumerate}
    \item[a)] In the non-strongly convex case with $\alpha=0$, let $\Mbar$ be an upper bound on $\jim{\max\{1,\norm{\grad g(\xunder^t)}\}}$ over all outer iterations $t$. Then
          \begin{align*}
            N_\ep & = \jim{\bigO\left\{
              \sqrt{\frac{\jim{\Lagg^c}}{\ep}}
            \sqrt{\DX^2+(\jim{\rplus_c})^2}+K_L+1\right\}}, \\
            C_\ep & = \jim{\bigO\left\{N_\ep+
              \frac{\Mbar}{\tilde L_{\mathrm{agg},0}}+
              \frac{(K_L+1)\Mbar}{\ep}
              \left[\DX^2+(\jim{\rplus_c})^2\right]\right\}}
          \end{align*}
    \item[b)] In the strongly convex case with $\alpha>0$, let $\Mbar$ be an upper bound on $\jim{\max\{1,\norm{\grad g(\xunder^t)}\}}$ over all outer iterations $t$. Then
          \begin{align*}
            N_\ep & = \jim{\bigO\left\{
              \sqrt{\frac{\jim{\Lagg^c}}{\alpha}}
              \log\left(1+\frac{\sqrt{\jim{\Lagg^c}\alpha}
            (\DX^2+(\jim{\rplus_c})^2)}{\ep}\right)+K_L+1\right\}}, \\
            C_\ep & = \jim{\bigO\left\{N_\ep+
              \frac{\Mbar}{\sqrt{\alpha\tilde L_{\mathrm{agg},0}}}+
              (K_L+1)\sqrt{\frac{1}{\alpha\ep}}\Mbar
              (\DX+\jim{\rplus_c})\right\}}.
          \end{align*}
  \end{enumerate}
\end{theorem}

\begin{proof}
  \jim{Under the restart convention stated after
    Algorithm~\ref{alg:adaptive-ACGD-S}, every actual restart begins a fresh run
    with the initialization required by
    Proposition~\ref{pr:a-ACGD-S-Q-convergence}. The running-max updates make the
    empirical estimates nondecreasing, and smoothness gives
    \[
      \alpha\leq\Ltil_{f,t}\leq L_f,
      \qquad 0\leq\Ltil_{g,t}\leq L_g.
    \]

    We first control the certificate radius. By
    Lemma~\ref{lm:PD-gap-to-suboptimality-feasibility}, the test in
    Line~\ref{aACDG-S-line:restart-trigger} cannot fail once
    $\rtil\geq2\jim{\rplus_c}$. Since the radius starts at $2c$, is
    retained across restarts, and is only doubled, it is doubled at most $K_r$
    times, and every radius encountered satisfies
    \[
      2c\leq\rtil<4\jim{\rplus_c}.
    \]
    A supported doubling reuses $F(\xbar^t)$, $g(\xbar^t)$, and
    $\FunderPD^t$ and therefore requires only scalar recomputation of the
    certificate.

    We next control the aggregate estimates. For every radius encountered,
    \[
      \Ltil_{f,t}+\rtil\Ltil_{g,t}
      \leq L_f+4\jim{\rplus_c}L_g
      =\jim{\Lagg^c}.
    \]
    The preceding bound applies to both the current radius $\rtil$ and the
    proposed radius $\rtil^+$. Hence, whenever either restart test is triggered,
    the empirical aggregate estimate in that test is at most
    $\jim{\Lagg^c}$, so $\Ltilagg<\jim{\Lagg^c}$. The update rule
    therefore gives
    \[
      2\Ltilagg\leq\Ltilagg^+
      \leq2\jim{\Lagg^c}.
    \]
    Consequently, there are at most $K_L$ actual restarts, every run uses
    $\Ltilagg\leq2\jim{\Lagg^c}$, and the total number of runs is at
    most $K_L+1$.

    Let $q\leq K_L$ be the number of actual restarts, let
    $\ell_k:=\Ltilagg^{(k)}$ be the aggregate estimate in run $k$, and set
    \[
      A_c:=\DX^2+(\jim{\rplus_c})^2.
    \]
    At every iteration that does not trigger a restart, the aggregate test and
    the radius loop ensure that the hypotheses of
    Proposition~\ref{pr:a-ACGD-S-Q-convergence} hold. Hence, for an absolute
    constant $C>0$,
    \[
      \left(\sum_{j=1}^t\wt[j]\right)\DeltaPD^t(\rtil_t)
      \leq C\ell_k A_c.
    \]
    Using the weight sequences in Theorems~\ref{thm:ANLA-S-non}
    and~\ref{thm:ANLA-S-str}, run $k$ therefore either restarts earlier or
    produces an iterate with $\DeltaPD(\rtil)\leq\ep$ within
    \[
      N_k=\begin{cases}
        \displaystyle
        \bigO\!\left(\sqrt{\frac{\ell_k A_c}{\ep}}+1\right),
         & \alpha=0, \\[2mm]
        \displaystyle
        \bigO\!\left(
        \sqrt{\frac{\ell_k}{\alpha}}
        \log\!\left(1+
        \frac{\sqrt{\ell_k\alpha}\,A_c}{\ep}\right)+1\right),
         & \alpha>0.
      \end{cases}
    \]
    Because $\ell_{k+1}\geq2\ell_k$ and
    $\ell_k\leq2\jim{\Lagg^c}$,
    \[
      \sum_{k=0}^{q}\sqrt{\ell_k}
      =\bigO\!\left(\sqrt{\jim{\Lagg^c}}\right).
    \]
    Summing the displayed per-run bounds, and using the monotonicity of the
    logarithm in the strongly convex case, gives the asserted value of $N_\ep$.

    It remains to count the inner-loop work. Each inner iteration uses a
    constant number of matrix-vector operations and projections, and each outer
    phase adds a constant amount of certificate work. Thus the inner-loop limits
    in Theorems~\ref{thm:ANLA-S-non} and~\ref{thm:ANLA-S-str}, with $R=1$, give
    the following per-run bounds:
    \[
      C_k=\begin{cases}
        \displaystyle
        \bigO\!\left(
        N_k+\frac{\Mbar}{\ell_k}
        +\frac{\Mbar A_c}{\ep}\right),
         & \alpha=0, \\[2mm]
        \displaystyle
        \bigO\!\left(
        N_k+\frac{\Mbar}{\sqrt{\alpha\ell_k}}
        +\frac{\Mbar(\DX+\jim{\rplus_c})}
        {\sqrt{\alpha\ep}}\right),
         & \alpha>0.
      \end{cases}
    \]
    The same geometric growth gives
    \[
      \sum_{k=0}^{q}\frac{1}{\ell_k}
      =\bigO\!\left(\frac{1}{\tilde L_{\mathrm{agg},0}}\right),
      \qquad
      \sum_{k=0}^{q}\frac{1}{\sqrt{\ell_k}}
      =\bigO\!\left(\frac{1}{\sqrt{\tilde L_{\mathrm{agg},0}}}\right).
    \]
    The accuracy-dependent final term in each case can occur in each of the at
    most $K_L+1$ runs. Summing the per-run bounds therefore yields the asserted
    value of $C_\ep$.

    Finally, whenever the post-radius-loop certificate satisfies
    $\DeltaPD(\rtil)\leq\ep$, Part~a) of
    Lemma~\ref{lm:PD-gap-to-suboptimality-feasibility} gives
    $F(\xbar)-\Fstar\leq\ep$. Moreover, exiting the radius loop certifies
    $\norm{[g(\xbar)]_+}\leq2\DeltaPD(\rtil)/\rtil$, and
    $\rtil\geq2c$. Therefore,
    \[
      F(\xbar)-\Fstar\leq\ep,
      \qquad
      \norm{[g(\xbar)]_+}
      \leq\frac{2\ep}{\rtil}
      \leq\frac{\ep}{c},
    \]
    proving the claimed accuracy.}
\end{proof}

Compared with the adaptive ACGD method in Subsection \ref{subsec:ANLA-search}, we have eliminated the \jim{constrained prox-mappings} required to compute the iterate $\xt$ in Line 6 of Algorithm \ref{alg:a-ACGD} and the termination criterion in Line 9, replacing them with matrix-vector multiplications and simple projections onto the feasible set $X$. These modifications enable the proposed method to scale to extremely large-scale problems with high dimensionality and many constraints.

\section{Numerical Experiments}\label{sec:numerical}

In this section, we conduct numerical experiments to evaluate the performance of the ACGD and ACGD-S methods on the following randomly generated quadratically constrained quadratic programming (QCQP) problems:
\begin{equation}
    \begin{split}
        \min_{x \in X} & f(x) = \frac{1}{2} x^\top Q x + c^\top x +\frac{\alpha}{2} \normsq{x}     \\
        s.t.           & g_i(x) = \frac{1}{2} x^\top P_i x + d_i^\top x \leq 0, \forall i \in [m],
    \end{split}
\end{equation}
where $Q$ and $P$ are randomly generated symmetric positive definite matrices, $c$ and $d$ are randomly generated vectors, and $X=[-10,10]^n$. The experiments are implemented in MATLAB and carried out on a MacBook Pro with an M3 Pro processor and 32 GB of memory. The reference optimal solution and objective value are computed using the Mosek solver.

\subsection{The Effect of Adaptive Stepsize Selection}

Recall that both the ACGD and ACGD-S methods crucially require an estimate of the aggregate Lipschitz smoothness constant, $\Ltilagg$. The true constant, denoted by $\Lagg$, is hard to determine in practice because it depends on the optimal dual variable $\lambda^*$. To assess how sensitive the proposed methods are to misspecification of this constant, we initialize $\Ltilagg$ using a scaled version of the true value:
$$\Ltilagg = \Lagg \times \text{scale factor},$$
where the $\text{scale factor}$ is chosen from $\{0.01, 0.05, 0.1, 0.5, 1, 3, 10\}$. For each method, we then calculate the number of iterations required to reach the desired relative optimality gap and constraint violation. The results are shown in Tables~\ref{tab:nonstrong} and~\ref{tab:strong}. A scale factor less than 1 indicates an underestimate of $\Lagg$, potentially leading to overly aggressive stepsizes and convergence to an infeasible solution, whereas a scale factor greater than 1 indicates an overestimate, resulting in overly conservative stepsizes.
\begin{figure}[ht]
    \centering
    \subfloat[Non-Strongly Convex Adaptive ACGD\label{fig:adaptive-acgd}]{%
        \includegraphics[width=0.8\textwidth]{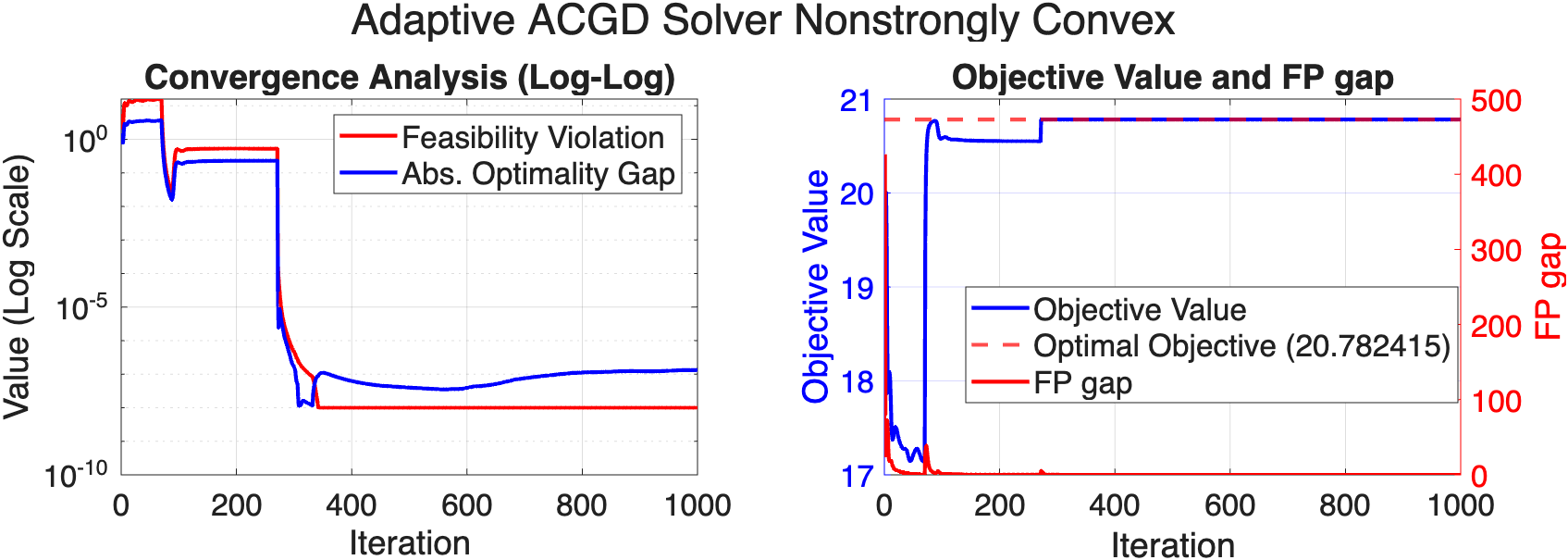}}
    \vspace{0.5cm}
    \subfloat[Strongly Convex Adaptive ACGD-S\label{fig:adaptive-acgd-s}]{%
        \includegraphics[width=0.8\textwidth]{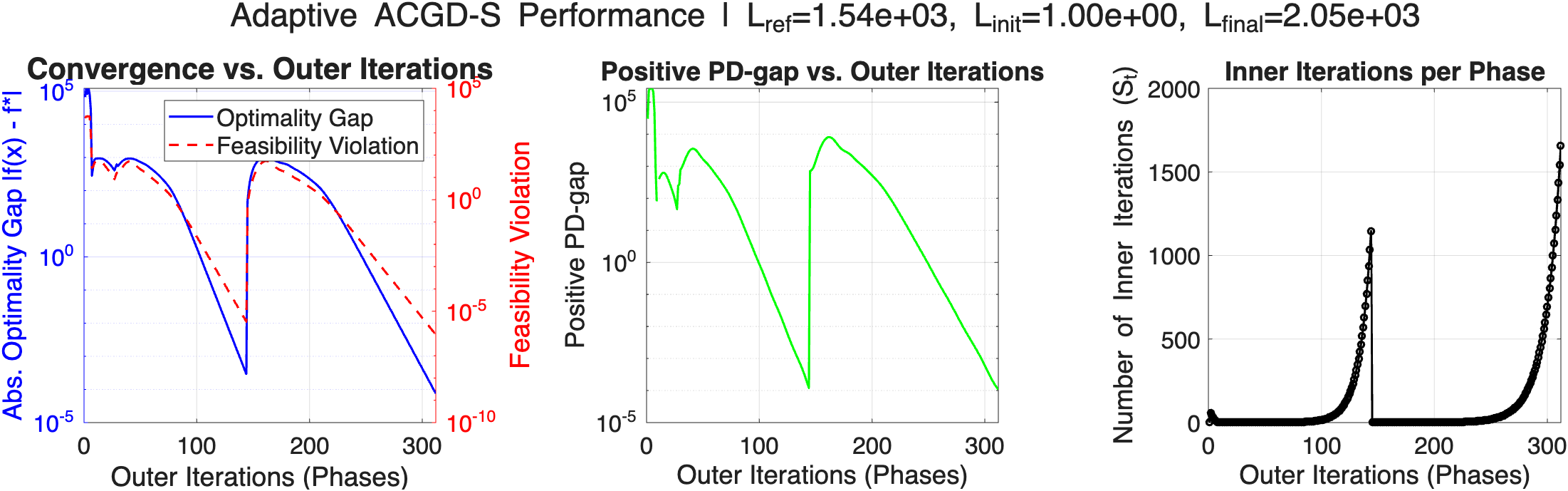}}

    \caption{Representative convergence plots for the Adaptive ACGD and Adaptive ACGD-S methods. For the non-strongly convex case, the FP-gap converges at an $\bigO(1/k^2)$ rate, while for the strongly convex case, the PD-gap converges at a linear rate, regardless of stepsize misspecification. By Lemma~\ref{lm:FP-gap-to-suboptimality-feasibility} and Lemma~\ref{lm:PD-gap-to-suboptimality-feasibility}, these convergence rates translate to corresponding rates for both the optimality gap and feasibility violation. When the feasibility violation exceeds a threshold, a restart is triggered to adjust the stepsize, which manifests as sudden jumps in the gap curves shown in these figures.}
    \label{fig:adaptive-comparison}
\end{figure}

\subsubsection*{Non-strongly Convex Case}

For the non-strongly convex case, we set $n = 30$, $m = 10$,  $\text{strong\_convexity} = 0.00\text{e+00}$, and $\Lagg = 3.2865\text{e+02}$. The termination criteria are $\text{relative\_gap} \leq 1.0\text{e-03}$, $\text{feasibility} \leq 1.0\text{e-03}$, and $\text{max\_iters} = 5000$.

\begin{table}[ht]
    \centering

    \caption{Number of iterations to reach $10^{-3}$ relative optimality gap and constraint violation---Non-Strongly Convex Case}
    \label{tab:nonstrong}
    \begin{tabular}{c|c|c|cc|cc}
        \hline
        Scaling & \jim{ACGD}        & \jim{Adaptive-ACGD} & \multicolumn{2}{c|}{\jim{ACGD-S}} & \multicolumn{2}{c}{\jim{Adaptive-ACGD-S}}                                         \\
        factor  & \jim{Outer Iters} & \jim{Outer Iters}   & \jim{Outer Iters}                 & \jim{Inner Iters}                         & \jim{Outer Iters} & \jim{Inner Iters} \\
        \hline
        10      & 24                & 24                  & $5000^{*}$                        & 55530                                     & $5000^{*}$        & 55530             \\
        3       & 13                & 13                  & 3190                              & 73518                                     & 3190              & 73518             \\
        1       & 6                 & 6                   & 1853                              & 73758                                     & 1853              & 73758             \\
        0.5     & 31                & 31                  & 1323                              & 74925                                     & 1323              & 74925             \\
        0.1     & $5000^{*}$        & 50                  & $5000^{*}$                        & 7843794                                   & 1542              & 88985             \\
        0.05    & $5000^{*}$        & 34                  & $5000^{*}$                        & 19226008                                  & 1558              & 90128             \\
        0.01    & $5000^{*}$        & 62                  & $5000^{*}$                        & 24739382                                  & 1118              & 76226             \\
        \hline
        \multicolumn{7}{l}{\footnotesize $^{*}$Timed out after 5000 outer iterations without meeting the termination criteria.}                                                   \\
    \end{tabular}
\end{table}

\subsubsection*{Strongly Convex Case}

For the strongly convex case, we set $n = 30$, $m = 10$, $\text{strong\_convexity} = 1.00\text{e+02}$, and $\Lagg = 3.1302\text{e+04}$. The termination criteria are $\text{relative\_gap} \leq 1.0\text{e-06}$, $\text{feasibility} \leq 1.0\text{e-06}$, and $\text{max\_iters} = 2000$.

\begin{table}[ht]
    \centering

    \caption{Number of iterations to reach $10^{-6}$ relative optimality gap and constraint violation - Strongly Convex Case}
    \label{tab:strong}
    \begin{tabular}{c|c|c|cc|cc}
        \hline
        Scaling & \jim{ACGD}        & \jim{Adaptive-ACGD} & \multicolumn{2}{c|}{\jim{ACGD-S}} & \multicolumn{2}{c}{\jim{Adaptive-ACGD-S}}                                         \\
        factor  & \jim{Outer Iters} & \jim{Outer Iters}   & \jim{Outer Iters}                 & \jim{Inner Iters}                         & \jim{Outer Iters} & \jim{Inner Iters} \\
        \hline
        10      & 316               & 316                 & $2000^{*}$                        & 43956                                     & $2000^{*}$        & 43956             \\
        3       & 171               & 171                 & 1164                              & 87345                                     & 1164              & 87345             \\
        1       & 93                & 93                  & 677                               & 86847                                     & 677               & 86847             \\
        0.5     & 118               & 118                 & 482                               & 86375                                     & 881               & 87172             \\
        0.1     & $2000^{*}$        & 136                 & $2000^{*}$                        & 9007720                                   & 885               & 87689             \\
        0.05    & $2000^{*}$        & 127                 & $2000^{*}$                        & 9335666                                   & 894               & 87705             \\
        0.01    & $2000^{*}$        & 84                  & $2000^{*}$                        & 9734320                                   & 1293              & 88801             \\
        \hline
        \multicolumn{7}{l}{\footnotesize $^{*}$Timed out after 2000 outer iterations without meeting the termination criteria.}                                                   \\
    \end{tabular}
\end{table}

\subsubsection*{Discussion}

The results warrant several comments. First, using an exact QP subproblem solver, the ACGD method converges quickly for both the non-strongly and strongly convex problems, provided the aggregate Lipschitz smoothness constant is correctly specified ($\text{scale factor} \geq 0.5$). However, when the scale factor falls below this threshold, the ACGD method becomes quite sensitive to its misspecification: for a $\text{scale factor} < 0.5$, the method fails to meet the termination criteria within the maximum number of iterations. \jim{In the non-strongly convex setting, underestimating $\Ltilagg$ can directly increase $S_t$, while the resulting aggressive outer steps can also increase $M_t$. These compounding effects lead to the extremely large inner-iteration counts reported in Table~\ref{tab:nonstrong}:}
\[
    \jim{S_t=\left\lceil\frac{tRM_t}{\Ltilagg}\right\rceil,\qquad
        M_t=\begin{cases}\norm{\nu^t},&\norm{\nu^t}>0,\\1,&\norm{\nu^t}=0,\end{cases}
        \qquad \nu^t=\grad g(\xunder^t).}
\]
In comparison, the Adaptive ACGD method is much more robust to stepsize misspecification. For all choices of the scale factor, the number of iterations required by the Adaptive ACGD method is comparable to the case where the aggregate Lipschitz smoothness constant is correctly specified.

Second, with the sliding subroutine to solve the inner problem inexactly, the ACGD-S method is even more sensitive to stepsize misspecification: there is only a narrow window ($\text{scale factor} \in [0.5, 3]$) for which the method can meet the termination criteria within the maximum number of iterations. Even when an overly conservative aggregate Lipschitz smoothness constant is used ($\text{scale factor} = 10$), the ACGD-S method fails to converge within the maximum number of iterations. This occurs because the number of inner iterations, \mbox{$S_t \propto 1/\Ltilagg$}, is too small, leading to inexact inner problem solutions and thus requiring many more outer iterations for convergence. In contrast, the Adaptive ACGD-S method is much more robust to stepsize misspecification. Whenever the initial guess of the aggregate Lipschitz smoothness constant is too small ($\text{scale factor} \leq 1$), the Adaptive ACGD-S method automatically restarts to adjust the stepsize. The number of iterations (both outer and inner) required by the Adaptive ACGD-S method is comparable to the case where the aggregate Lipschitz smoothness constant is correctly specified.

\subsection{The Effect of the Sliding Subroutine}

To assess the effect of the sliding subroutine, we now compare the performance of the Adaptive ACGD and Adaptive ACGD-S methods across different problem sizes. We focus on the more stable and practical adaptive methods and choose an initial guess of $1$ for the aggregate Lipschitz smoothness constant in all experiments. For each method and accuracy level, we record the first iteration at which both the feasibility violation and the relative optimality gap reach the desired accuracy. The operation counts and computation times reported in Tables~\ref{tab:adaptive-comparison} and~\ref{tab:adaptive-comparison-strong} are averaged over 5 random trials. The tables report gradient evaluations, matrix-vector operations, and computation time for different accuracy levels. 

\begin{table}[ht]
    \centering
    \setlength{\tabcolsep}{2pt}

    \caption{Performance comparison of Adaptive ACGD, Adaptive ACGD-S, and APDB across different problem sizes. All values are averaged over 5 random trials.}
    \label{tab:adaptive-comparison}
    \begin{tabular}{c|c|cc|ccc|cc}
        \hline
        \jim{Problem}                                                                                                                         & \jim{Accuracy}   & \multicolumn{2}{c|}{\jim{Adaptive-ACGD}} & \multicolumn{3}{c|}{\jim{Adaptive-ACGD-S}} & \multicolumn{2}{c}{\jim{APDB}}                                                                                       \\
        \jim{Size $(m,n)$}                                                                                                                    & \jim{Level}      & \jim{Grad Evals}                         & \jim{Time (s)}                             & \jim{Grad Evals}               & \jim{Matrix-Vec Ops} & \jim{Time (s)} & \jim{$2\times$ Grad Evals} & \jim{Time (s)} \\
                                                                                                                                              &                  &                                          &                                            &                                &                      &                & \jim{$=$ Matrix-Vec Ops}   &                \\
        \hline
        \multirow{3}{*}{\begin{tabular}{@{}c@{}}$(10, 20)$ \\ $L_{\text{agg}} = 2.14\text{e+01}$\end{tabular}}     & $1.0\text{e-01}$ & 2.4                                      & 0.03s                                      & 23.2                           & 191.6                & 0.02s          & 253.2                      & 0.03s          \\
                                                                                                                                              & $1.0\text{e-02}$ & 7.4                                      & 0.07s                                      & 63.4                           & 1246.8               & 0.06s          & 2975.6                     & 0.11s          \\
                                                                                                                                              & $1.0\text{e-03}$ & 22.2                                     & 0.23s                                      & 190.4                          & 11465.2              & 0.18s          & 29258.8                    & 0.87s          \\
        \hline
        \multirow{3}{*}{\begin{tabular}{@{}c@{}}$(100, 200)$ \\ $L_{\text{agg}} = 1.01\text{e+02}$\end{tabular}}   & $1.0\text{e-01}$ & 37.8                                     & 1.01s                                      & 26.0                           & 622.4                & 0.08s          & 9.6                        & 0.00s          \\
                                                                                                                                              & $1.0\text{e-02}$ & 47.2                                     & 1.26s                                      & 69.8                           & 2674.0               & 0.21s          & 5058.0                     & 0.48s          \\
                                                                                                                                              & $1.0\text{e-03}$ & 76.2                                     & 2.04s                                      & 208.0                          & 22738.0              & 0.62s          & 50047.6                    & 4.57s          \\
        \hline
        \multirow{3}{*}{\begin{tabular}{@{}c@{}}$(200, 500)$ \\ $L_{\text{agg}} = 1.43\text{e+02}$\end{tabular}}   & $1.0\text{e-01}$ & 38.3                                     & 3.41s                                      & 21.0                           & 1560.0               & 0.27s          & 4.0                        & 0.00s          \\
                                                                                                                                              & $1.0\text{e-02}$ & 59.0                                     & 5.25s                                      & 54.0                           & 3940.0               & 0.70s          & 7976.0                     & 6.34s          \\
                                                                                                                                              & $1.0\text{e-03}$ & 124.3                                    & 11.08s                                     & 158.3                          & 27372.6              & 2.04s          & 79200.8                    & 62.79s         \\
        \hline
        \multirow{3}{*}{\begin{tabular}{@{}c@{}}$(1000, 2000)$ \\ $L_{\text{agg}} = 3.76\text{e+02}$\end{tabular}} & $1.0\text{e-01}$ & 36.6                                     & 81.20s                                     & 19.2                           & 3540.4               & 4.76s          & 4.0                        & 0.03s          \\
                                                                                                                                              & $1.0\text{e-02}$ & 107.8                                    & 239.15s                                    & 47.4                           & 6268.8               & 11.75s         & 6296.0                     & 61.90s         \\
                                                                                                                                              & $1.0\text{e-03}$ & 136.2                                    & 302.17s                                    & 136.6                          & 32919.6              & 33.87s         & 62478.0                    & 617.84s        \\
        \hline
    \end{tabular}
\end{table}

\begin{table}[ht]
    \centering
    \setlength{\tabcolsep}{2pt}

    \caption{Performance comparison of Adaptive ACGD, Adaptive ACGD-S, and APDB across different problem sizes - Strongly Convex Case. All values are averaged over 5 random trials.}
    \label{tab:adaptive-comparison-strong}
    \begin{tabular}{c|c|cc|ccc|cc}
        \hline
        \jim{Problem}                                                                                                 & \jim{Accuracy}   & \multicolumn{2}{c|}{\jim{Adaptive-ACGD}} & \multicolumn{3}{c|}{\jim{Adaptive-ACGD-S}} & \multicolumn{2}{c}{\jim{APDB}}                                                                                       \\
        \jim{Size $(m,n)$}                                                                                            & \jim{Level}      & \jim{Grad Evals}                         & \jim{Time (s)}                             & \jim{Grad Evals}               & \jim{Matrix-Vec Ops} & \jim{Time (s)} & \jim{$2\times$ Grad Evals} & \jim{Time (s)} \\
                                                                                                                      &                  &                                          &                                            &                                &                      &                & \jim{$=$ Matrix-Vec Ops}   &                \\
        \hline
        \multirow{3}{*}{\begin{tabular}{@{}c@{}}$(10, 20)$ \\ $\kappa = 1.40\text{e+01}$\end{tabular}}     & $1.0\text{e-02}$ & 7.0                                      & 0.03s                                      & 110.6                          & 993.6                & 0.23s          & 30886.0                    & 0.90s          \\
                                                                                                                      & $1.0\text{e-04}$ & 18.8                                     & 0.09s                                      & 133.8                          & 4542.4               & 0.28s          & 49872.8                    & 1.44s          \\
                                                                                                                      & $1.0\text{e-06}$ & 31.6                                     & 0.15s                                      & 189.0                          & 43798.4              & 0.40s          & 92884.4                    & 2.67s          \\
        \hline
        \multirow{3}{*}{\begin{tabular}{@{}c@{}}$(100, 200)$ \\ $\kappa = 6.48\text{e+01}$\end{tabular}}   & $1.0\text{e-02}$ & 17.0                                     & 0.26s                                      & 205.4                          & 4608.0               & 0.39s          & 24026.4                    & 2.23s          \\
                                                                                                                      & $1.0\text{e-04}$ & 40.0                                     & 0.60s                                      & 249.4                          & 13048.0              & 0.47s          & 34008.0                    & 3.13s          \\
                                                                                                                      & $1.0\text{e-06}$ & 61.2                                     & 0.92s                                      & 293.2                          & 95568.0              & 0.55s          & 132602.8                   & 12.09s         \\
        \hline
        \multirow{3}{*}{\begin{tabular}{@{}c@{}}$(200, 500)$ \\ $\kappa = 1.19\text{e+02}$\end{tabular}}   & $1.0\text{e-02}$ & 33.2                                     & 2.72s                                      & 292.6                          & 17321.2              & 6.32s          & 33953.2                    & 26.94s         \\
                                                                                                                      & $1.0\text{e-04}$ & 75.4                                     & 6.19s                                      & 353.8                          & 31583.6              & 7.65s          & 49604.0                    & 39.17s         \\
                                                                                                                      & $1.0\text{e-06}$ & 115.0                                    & 9.44s                                      & 415.0                          & 172470.8             & 8.97s          & 204872.4                   & 160.77s        \\
        \hline
        \multirow{3}{*}{\begin{tabular}{@{}c@{}}$(1000, 2000)$ \\ $\kappa = 2.62\text{e+02}$\end{tabular}} & $1.0\text{e-02}$ & 40.3                                     & 139.97s                                    & 386.3                          & 43972.6              & 198.52s        & 27114.8                    & 244.57s        \\
                                                                                                                      & $1.0\text{e-04}$ & 84.0                                     & 291.42s                                    & 472.3                          & 63869.4              & 242.69s        & 48224.4                    & 432.51s        \\
                                                                                                                      & $1.0\text{e-06}$ & 118.0                                    & 409.53s                                    & 557.7                          & 258749.4             & 286.52s        & $205379.2^{*}$             & 1800.02s       \\
        \hline
        \multicolumn{9}{l}{\footnotesize $^{*}$APDB timed out after 30 minutes without reaching $10^{-6}$ in any of the five trials.}                                                                                                                                                                                                                   \\
    \end{tabular}
\end{table}

We make a few observations. First, both methods are highly efficient for solving constrained problems. They solve moderate-sized problems ($(m,n) = (10, 20)$ and $(100, 200)$) in under a second and larger-scale problems ($(m,n) = (1000, 2000)$) in just a few minutes. Notably, compared to~\cite{xu2020first}, both methods handle a large number of constraints ($m = 1000$) without a significant drop in performance.

Second, the Adaptive ACGD method, which uses an exact inner QP solver, \jim{generally} requires fewer oracle evaluations than the Adaptive ACGD-S method to reach the same accuracy. This suggests that the simpler Adaptive ACGD method may be preferable when oracle evaluations are computationally expensive.

Third, regarding computation time, the Adaptive ACGD-S method \jim{generally outperforms the Adaptive ACGD method for solving large-scale problems to a high accuracy}. This efficiency gain stems from the use of the sliding subroutine to solve the inner problem inexactly, which is substantially faster than the exact QP solver utilized by Adaptive ACGD. This advantage is particularly pronounced for large-scale problems, where solving the QP subproblem becomes computationally expensive.

\jim{Additionally, we test APD with backtracking line search (APDB)~\cite{aybat2021primal}. Each APDB backtracking trial evaluates one gradient and requires roughly two matrix-vector multiplications. As shown in the Adaptive ACGD-S and APDB columns of Tables~\ref{tab:adaptive-comparison} and~\ref{tab:adaptive-comparison-strong}, Adaptive ACGD-S uses a comparable number of matrix-vector multiplications while requiring significantly fewer expensive gradient evaluations. Empirically, this results in substantially shorter runtimes at higher accuracy levels.}

\section{Conclusion}\label{sec:conclusion}
To sum up, this paper proposes two efficient methods for large-scale function-constrained optimization. The simple ACGD method has the optimal oracle complexity, but it requires access to \jim{a constrained prox-mapping oracle}. The more complicated ACGD-S method has both the optimal oracle complexity and the optimal computation complexity. Lower complexity bounds are provided to demonstrate that the oracle complexity of both ACGD and ACGD-S is unimprovable for the general case of first-order methods. Together they provide a complete characterization of the difficulty of solving a smooth function-constrained optimization problem from both the oracle complexity and the computation complexity perspective.

\section{Statements and Declarations}
\subsection*{Funding}
This work is partially supported by the ONR grant N00014-20-1-2089.

\subsection*{Competing Interests}
The authors have no relevant financial or non-financial interests to disclose.

\subsection*{Availability of data and materials}
The numerical experiments use randomly generated problem instances as described in Section~\ref{sec:numerical}; no external datasets were used.
\medskip

\bibliographystyle{siam}
\bibliography{Reference}

\begin{thebibliography}{10}

\bibitem{Beck2017First}
{\sc A.~Beck}, {\em First-order methods in optimization}, vol.~25, SIAM, 2017.

\bibitem{bertsekas2009convex}
{\sc D.~Bertsekas}, {\em Convex optimization theory}, vol.~1, Athena Scientific, 2009.

\bibitem{boob2022stochastic}
{\sc D.~Boob, Q.~Deng, and G.~Lan}, {\em Stochastic first-order methods for convex and nonconvex functional constrained optimization}, Mathematical Programming, 197 (2023), pp.~215--279.

\bibitem{deng2026uniformly}
{\sc Q.~Deng, G.~Lan, and Z.~Lin}, {\em Uniformly optimal and parameter-free first-order methods for convex and function-constrained optimization}, INFORMS Journal on Computing,  (2026).

\bibitem{gandy2005portfolio}
{\sc R.~Gandy}, {\em Portfolio optimization with risk constraints}, PhD thesis, Universit{\"a}t Ulm, 2005.

\bibitem{aybat2021primal}
{\sc E.~Y. Hamedani and N.~S. Aybat}, {\em A primal-dual algorithm with line search for general convex-concave saddle point problems}, SIAM Journal on Optimization, 31 (2021), pp.~1299--1329.

\bibitem{lan2016gradient}
{\sc G.~Lan}, {\em Gradient sliding for composite optimization}, Mathematical Programming, 159 (2016), pp.~201--235.

\bibitem{LanBook}
\leavevmode\vrule height 2pt depth -1.6pt width 23pt, {\em First-order and stochastic Optimization Methods for Machine Learning}, Springer-Nature, 2020.

\bibitem{lan2013iteration}
{\sc G.~Lan and R.~D. Monteiro}, {\em Iteration-complexity of first-order penalty methods for convex programming}, Mathematical Programming, 138 (2013), pp.~115--139.

\bibitem{lanOuyang2016GradientSliding}
{\sc G.~Lan and Y.~Ouyang}, {\em Accelerated gradient sliding for structured convex optimization}, Computational Optimization and Applications, 82 (2022), pp.~361--394.

\bibitem{lan2023optimal}
{\sc G.~Lan, Y.~Ouyang, and Z.~Zhang}, {\em Optimal and parameter-free gradient minimization methods for convex and nonconvex optimization}, Mathematical Programming,  (2026).

\bibitem{lan2021graph}
{\sc G.~Lan, Y.~Ouyang, and Y.~Zhou}, {\em Graph topology invariant gradient and sampling complexity for decentralized and stochastic optimization}, SIAM Journal on Optimization, 33 (2023), pp.~1647--1675.

\bibitem{lan2022optimal}
{\sc G.~Lan and Z.~Zhang}, {\em Optimal methods for convex risk-averse distributed optimization}, SIAM Journal on Optimization, 33 (2023), pp.~1518--1557.

\bibitem{lin2018level}
{\sc Q.~Lin, S.~Nadarajah, and N.~Soheili}, {\em A level-set method for convex optimization with a feasible solution path}, SIAM Journal on Optimization, 28 (2018), pp.~3290--3311.

\bibitem{lin2025adaptive}
{\sc Q.~Lin, N.~Soheili, R.~Ma, and S.~Nadarajah}, {\em An adaptive parameter-free and projection-free restarting level set method for constrained convex optimization under the error bound condition}, Journal of Machine Learning Research, 26 (2025), pp.~1--45.

\bibitem{nemirovski2001lectures}
{\sc A.~Nemirovski}, {\em Lectures on modern convex optimization}, in Society for Industrial and Applied Mathematics (SIAM, Citeseer, 2001.

\bibitem{nemirovsky1991optimality}
{\sc A.~S. Nemirovsky}, {\em On optimality of krylov's information when solving linear operator equations}, Journal of Complexity, 7 (1991), pp.~121--130.

\bibitem{nemirovsky1983problem}
{\sc A.~S. Nemirovsky and D.~B. Yudin}, {\em Problem complexity and method efficiency in optimization.}, John Wiley UK/USA, 1983.

\bibitem{nesterov1983method}
{\sc Y.~Nesterov}, {\em A method of solving a convex programming problem with convergence rate {$O(1/k^2)$}}, in Doklady Akademii Nauk, vol.~269, Russian Academy of Sciences, 1983, pp.~543--547.

\bibitem{nesterov2003introductory}
\leavevmode\vrule height 2pt depth -1.6pt width 23pt, {\em Introductory lectures on convex optimization: A basic course}, vol.~87, Springer Science \& Business Media, 2003.

\bibitem{nesterov2013gradient}
{\sc Y.~Nesterov}, {\em Gradient methods for minimizing composite functions}, Mathematical programming, 140 (2013), pp.~125--161.

\bibitem{nesterov2018lectures}
{\sc Y.~Nesterov}, {\em Lectures on convex optimization}, vol.~137, Springer, 2018.

\bibitem{ouyang2021lower}
{\sc Y.~Ouyang and Y.~Xu}, {\em Lower complexity bounds of first-order methods for convex-concave bilinear saddle-point problems}, Mathematical Programming, 185 (2021), pp.~1--35.

\bibitem{rigollet2011neyman}
{\sc P.~Rigollet and X.~Tong}, {\em Neyman-pearson classification, convexity and stochastic constraints}, Journal of Machine Learning Research,  (2011).

\bibitem{xu2020first}
{\sc Y.~Xu}, {\em First-order methods for problems with {O(1)} functional constraints can have almost the same convergence rate as for unconstrained problems}, SIAM Journal on Optimization, 32 (2022), pp.~1759--1790.

\bibitem{yang2022data}
{\sc S.~Yang, X.~Li, and G.~Lan}, {\em Data-driven minimax optimization with expectation constraints}, Operations Research, 73 (2024), pp.~1345--1365.

\bibitem{zafar2017fairness}
{\sc M.~B. Zafar, I.~Valera, M.~G. Rogriguez, and K.~P. Gummadi}, {\em Fairness constraints: Mechanisms for fair classification}, in Artificial intelligence and statistics, PMLR, 2017, pp.~962--970.

\bibitem{zhang2019efficient}
{\sc Z.~Zhang, S.~Ahmed, and G.~Lan}, {\em Efficient algorithms for distributionally robust stochastic optimization with discrete scenario support}, SIAM Journal on Optimization, 31 (2021), pp.~1690--1721.

\bibitem{zhang2020optimal}
{\sc Z.~Zhang and G.~Lan}, {\em Optimal methods for convex nested stochastic composite optimization}, Mathematical Programming, 212 (2025), pp.~1--48.

\bibitem{zhang2025linearly}
{\sc Z.~Zhang and S.~Sra}, {\em Linearly convergent algorithms for nonsmooth problems with unknown smooth pieces}, arXiv preprint arXiv:2507.19465,  (2025).

\end{thebibliography}
\section{Appendix}

\def \barg {\hat{g}}
\begin{lemma}\label{lm:agg_strong_cvxity}
Let $\Lambda \subset \R^m_+$  and a convex vector-valued function $g:\R^n \rightarrow \R^m $ be given. If $\sumi \lam_i g_i$ is $L$-smooth for all $\lam \in \Lam$, i.e., $\norm{\sumi \lami {(\grad g_i(x) -  \grad g_i(\xbar)} )} \leq L \norm{x - \xbar} \forall x, \xbar \in \jim{\R^n}, \forall \lam \in \Lam$, the Bregman distance function generated by its (vector-valued) conjugate function $U_{g^*}$ satisfies 
$$ \inner{\lam}{U_g^*(\nu, \nubar)} \geq \normsq{\sumi \lami (\nui - \nubar_i)}/(2L)\ \forall \nu, \nubar \in \{\grad g(x): x \in \R^n\}.$$
\end{lemma}
\begin{proof}
Let  $\lam \in \Lam$,  $\nu = \grad g(x),$ and $\nubar = \grad g(\xbar) $ be given. Consider the function $\barg := \sumi \lami g_i $. Clearly,  $\barg$ is $L$-smooth such that the Bregman distance function generated by its conjugate satisfies $U_{\barg^*}(v, \bar{v}) \geq \normsq{v - \bar{v}}/(2 L).$ Since $\grad \barg(x) = \lambda^\top \grad g(x)$, we have 
\begin{align*}
 \lambda^\top U_{g^*}(\nu; \nubar) &\stackrel{(a)}{=} \sumi \lami U_{g^*_i}(\nu_i; \nubar_i)=\sumi \lami U_{g_i}(\xbar; x) \\
&= \sumi \lami [g_i(\xbar) - g_i(x)  - \inner{\grad g_i(x)}{\xbar - x}]\\
&= \barg(\xbar) - \barg(x) - \inner{\grad \barg(x)}{\xbar - x}\\
&\stackrel{(b)}{=} U_{\barg}(\xbar; x)
 = U_{\barg^*}(\lam^\top \nu; \lam^\top \nubar) \geq \normsq{\sumi \lami (\nui - \nubar_i)}/(2L),
\end{align*}
where (a) and (b) follows from the algebraic identity between Bregman distance functions generated by Fenchel conjugates $(h, h^*)$, $U_h(y;\bar y) = U_{h^*}(\grad h(\bar y); \grad h(y)).$
\end{proof}
\vgap
The next two lemmas provide some basic algebraic identities useful for deriving complexity bounds. 
\begin{lemma}\label{lm:al_fact}
Given an $x > 0$, the following algebraic relation is valid:
\begin{equation}\label{al-fact:wt}
h(y):= (1 + 1/x)^{y -3} \leq y\  \forall\  2x \geq  y\geq 1.
\end{equation}
\end{lemma}

\begin{proof}
First, we show the relation for $y=2x$, i.e. $h(2x)= (1 + 1/x)^{2x -3} \leq 2x\  \forall x >0.$
Let's consider two cases. If $ x \geq 4$, we have 
$$(1  + 1/x)^{2x -3} \leq [(1+ 1/x)^x]^2 \leq \exp(2) \leq 8 \leq 2x.$$
If $1.5<x < 4$, we have 
$$ [(1 + 1/x)^{x -1.5}]^2 \leq [1/(1 - \tfrac{x-1.5}{x})]^2 = (x/1.5)^2 = (x/2.25) \times x \leq 2 x.$$
\jim{Moreover, if $0.5\leq x\leq 1.5$, we have 
$$[(1 + 1/x)^{x -1.5}]^2 \leq 1 \leq 2x.$$
Thus we get $h(2x) \leq 2x \forall x\in [0.5, +\infty)$. }
Since $h(1) \leq 1$, the relation in \eqref{al-fact:wt} follows from the convexity of $h$ with respect to $y$.
\end{proof}
\vgap 
\begin{lemma}\label{lm:alg-fact-St}
Given non-negative parameters $\Delta_t$, $H$, and $h \in \{0,1\}$, suppose $\St = \min\{S\in \mathbb{N}_+: \sums[S] H + (s-h) \geq \Delta_t\}$, and $\Gamma$ is the non-negative root of $ \sums[\St]  [\Gamma H + (s-h)] = \Gamma^2\Delta_t$, then $\St$ satisfies 
$$\St / \Gamma \geq \St -2.$$
\end{lemma}

\begin{proof}
\jim{The statement clearly holds for $S_t\leq 2$, so it suffices to consider the case with $S_t\geq 3$.}

Suppose for the sake of contradiction that $\St / \Gamma < \St -2$. On the one hand, the definition of $\St$ implies that 
\begin{align*}
\Delta_t = H (\tfrac{\St}{\Gamma}) + \sums \tfrac{(s-h)}{\Gamma^2} = H (\tfrac{\St}{\Gamma}) +  \tfrac{1}{2}\tfrac{(\St-\jim{2}h+1)}{\Gamma} \tfrac{\St}{\Gamma}.
\end{align*}
On the other hand, the choice of $\St$ implies that $\sums[\St-1] [H + (s-h)] < \Delta_t$ such that 
$$\Delta_{t}> H(S_{t}-1)+\frac{1}{2}(S_{t}-1)(S_{t}-2h)$$

\jim{
Now we consider two cases. If $h=1$, then the contradiction assumption implies that $S_{t}-2>\frac{S_{t}}{\Gamma}>\frac{S_{t}-2h+1}{\Gamma}$
such that 
\[
\Delta_{t}> H(S_{t}-1)+\frac{1}{2}(S_{t}-1)(S_{t}-2h) >H\frac{S_{t}}{\Gamma}+\frac{1}{2}\frac{(S_{t}-2h+1)S_{t}}{\Gamma^{2}} .
\]

If $h=0$, the contradiction assumption implies that $1>\frac{1}{\Gamma}\frac{S_{t}}{S_{t}-2} $
and $S_{t}-1>\frac{S_{t}}{\Gamma}$ such that 
\begin{align*}
\Delta_{t} & >H(S_{t}-1)+\frac{1}{2}(S_{t}-1)(S_{t})\\
 & >H\frac{S_{t}}{\Gamma}+\frac{1}{2\Gamma^{2}}(S_{t}+1)(S_{t}-2)\frac{S_{t}^{2}}{(S_{t}-2)^{2}}\\
 & >H\frac{S_{t}}{\Gamma}+\half\frac{(S_{t}-2h+1)S_{t}}{\Gamma^{2}}.
\end{align*}

Therefore, both cases lead to the desired contradiction.

}
\end{proof}

\end{document}